\documentclass[11pt]{amsart}

\usepackage{amsfonts,amsmath,amssymb,amsthm}
\allowdisplaybreaks[4]

\usepackage[colorlinks,linkcolor=red,citecolor=green]{hyperref}
\usepackage{makecell}
\usepackage{braket}
\usepackage{graphics}
\usepackage{xcolor}
\usepackage{mathtools}
\usepackage{tikz}
\usetikzlibrary{arrows.meta}

\newcommand{\dd}{\mathrm{d}}

\newtheorem{theorem}{Theorem}
\newtheorem{lemma}[theorem]{Lemma}
\newtheorem{remark}[theorem]{Remark}
\newtheorem{proposition}[theorem]{Proposition}

\newtheorem{corollary}[theorem]{Corollary}

\numberwithin{theorem}{section}
\numberwithin{equation}{section}

\begin{document}
	
	\setlength{\baselineskip}{1.1\baselineskip}
	
	\title[Log-concavity of the first Dirichlet eigenfunction in hyperbolic space]{Log-concavity and level set curvature estimates for the first Dirichlet eigenfunction in hyperbolic space}
	
	\author{Wei Zhang}
	\address{School of Mathematics and Statistics\\
		Lanzhou University\\
		Lanzhou, 730000, Gansu Province, China.}
	\email{zhangw@lzu.edu.cn}
	
	\author{Qi Zhou}
	\address{School of Mathematics and Statistics\\
		Lanzhou University\\
		Lanzhou, 730000, Gansu Province, China.}
	\email{zhouqi2025@lzu.edu.cn, zhouqimath20@lzu.edu.cn}

	\date{\today}
	
	\maketitle
	
	\begin{abstract}
		In this paper, we study the first Dirichlet eigenfunction on bounded, smooth and strictly horo-convex domains in hyperbolic space $\mathbb{H}^n$. We prove its strict $\log$-concavity in arbitrary dimension and establish the Gauss curvature bound $K>1$ for its level sets. In dimensions two and three, we further obtain quantitative lower bounds for $K$ in terms of the minimum Gauss curvature of the boundary. In higher dimensions, the product of the two smallest principal curvatures of each level hypersurface exceeds one. Furthermore, we extend the strict $\log$-concavity result to bounded smooth horo-convex domains. Our approach combines maximum or minimum principles for suitable auxiliary functions, boundary convexity estimates and a domain deformation argument.
	\end{abstract}
	
	
	
	\section{Introduction}
	
	Given an elliptic boundary value problem on a bounded convex domain in a Riemannian manifold, does the solution inherit any convexity properties from the domain? This is a classical question in elliptic partial differential equations. There are two main aspects to such problems, namely the convexity of the solution (or of a suitable transformation of it) and the convexity of its level sets. A further goal is to obtain quantitative convexity estimates for solutions and in particular curvature estimates for their level sets in terms of the geometry of the boundary. The first Dirichlet eigenvalue problem and the torsion problem are two classical models for the study of these questions. For a bounded smooth domain $\Omega$ in a Riemannian manifold $(M^n,g)$, these problems are given by
	\begin{equation}\label{Equation-Dirichlet problem}
		\begin{cases}
			\Delta u=-\lambda u & \mathrm{in}\ \Omega,\\
			\hspace{0.4cm} u>0  & \mathrm{in}\ \Omega,\\
			\hspace{0.4cm} u=0 & \mathrm{on}\ \partial\Omega,
		\end{cases}
	\end{equation}
	and
	\begin{equation}\label{Equation-Torsion problem}
		\begin{cases}
			\Delta u=-2 & \mathrm{in}\ \Omega,\\
			\hspace{0.4cm} u=0 & \mathrm{on}\ \partial\Omega,
		\end{cases}
	\end{equation}
	respectively. Regarding the convexity of solutions to these two problems, we consider domains whose boundary principal curvatures are nonnegative. One studies the concavity of $\log u$ for the first Dirichlet eigenfunction $u$ satisfying \eqref{Equation-Dirichlet problem}. For the torsion function $u$ satisfying \eqref{Equation-Torsion problem}, the corresponding question concerns the convexity of $-\sqrt{u}$. In particular, both properties imply convexity of the superlevel sets of $u$.
	
	In the Euclidean plane, Makar-Limanov \cite{Makar-Limanov1971} proved in 1971 that $-\sqrt{u}$ is strictly convex when $u$ is the torsion function of a bounded convex domain. His proof is based on the superharmonicity of the $P$-function defined by $$P=2u\det\nabla^2u+2u_1u_2u_{12}-u_{11}u_2^2-u_{22}u_1^2.$$ For the first Dirichlet eigenvalue problem, Brascamp and Lieb \cite{Brascamp-Lieb1976} established $\log$-concavity of the first eigenfunction on bounded convex domains in $\mathbb{R}^n$. Acker, Payne and Philippin \cite{Acker-PP1981} subsequently adapted Makar-Limanov's idea to the two-dimensional eigenvalue problem to give another proof of $\log$-concavity and obtain pointwise curvature estimates for the level curves. The study of $\log$-concavity also has applications to fundamental gap estimates in spectral theory. In particular, Andrews and  Clutterbuck \cite{Andrews-Clutterbuck2011} established a sharp modulus of concavity for the logarithm of the first eigenfunction and used it to prove the fundamental gap conjecture for convex Euclidean domains. 
	
	Several general methods have been developed to study the convexity properties of solutions to elliptic equations. Korevaar \cite{Korevaar1983b} introduced the concavity maximum principle for quasilinear elliptic and parabolic equations by applying the classical maximum principle to a suitable concavity function. Kennington \cite{Kennington1985} further refined this method and established power concavity results for a broad class of boundary value problems; see also the monograph of Kawohl \cite{Kawohl1985}. An alternative approach developed by Alvarez, Lasry and Lions \cite{Alvarez-LL1997} uses convex envelopes together with the comparison principle and is sometimes referred to as the macroscopic convexity principle. The microscopic convexity principle, also known as the constant rank theorem, was introduced and developed through the work of Caffarelli and Friedman \cite{Caffarelli-Friedman1985}, Singer, Wong, Yau and Yau \cite{Singer-WYY1985} and Korevaar and Lewis \cite{Korevaar-Lewis1987}. In particular, these approaches provide different ways to establish the convexity of $-\log u$ and $-\sqrt{u}$ for \eqref{Equation-Dirichlet problem} and \eqref{Equation-Torsion problem}, respectively, on bounded convex domains in $\mathbb{R}^n$.
	
	A natural question is whether the above convexity results extend to domains  in a Riemannian manifold $M^n$. Convexity properties on space forms have also received considerable attention. In the following discussion, we focus on the first Dirichlet eigenvalue problem. For developments concerning the torsion problem, see our earlier work \cite{Zhang-Zhou2026} and the references therein. When $M^n$ is the sphere $\mathbb{S}^n$, Lee and Wang \cite{Lee-Wang1987} proved $\log$-concavity of the first Dirichlet eigenfunction on convex domains by adapting the continuity argument of Singer, Wong, Yau and Yau \cite{Singer-WYY1985}. They also obtained a lower bound for the fundamental gap. In his survey article, Korevaar \cite{Korevaar1987} reported unpublished joint work with Treibergs extending the constant rank theorem to spaces of nonnegative constant curvature and indicated its application to the first eigenvalue problem on the sphere. Further progress on the sphere was made by Seto, Wang and Wei \cite{Seto-WW2019}, who established improved $\log$-concavity estimates and sharp fundamental gap bounds for $n\geq 3$ under a diameter restriction. He, Wei and Zhang \cite{He-WZ2020} subsequently relaxed this restriction, while Dai, Seto and Wei \cite{Dai-SW2021} obtained the corresponding sharp fundamental gap bound for $n=2$ using a different one-dimensional comparison model.
	
	Beyond spaces of constant sectional curvature, Khan, Nguyen, Tuerkoen and Wei \cite{Khan-NTW2025} obtained strong $\log$-concavity estimates for the first eigenfunction on convex domains in positively curved surfaces under a condition involving the curvature and its derivatives, together with corresponding lower bounds for the fundamental gap. A classical example of a manifold with positive but nonconstant sectional curvature is complex projective space $\mathbb{CP}^n$, $n\geq 2$. In this setting, Xia \cite{Xia2026} showed strict $\log$-concavity of the first eigenfunction on smooth $1$-convex domains in $\mathbb{CP}^2$ by establishing quantitative Hessian estimates with applications to level set convexity and the fundamental gap. Under a symmetry assumption, Hsu \cite{Hsu2026} proved strict $\log$-concavity on smooth convex domains in $\mathbb{CP}^n$ for $n=2,3,4$. More recently, Cho \cite{Cho2026} established quantitative $\log$-concavity and fundamental gap estimates in arbitrary dimension under a stronger boundary convexity assumption without symmetry restrictions.
	
	The situation changes substantially in negative curvature. Shih \cite{Shih1989} constructed a geodesically convex domain in the hyperbolic plane $\mathbb{H}^2$ whose first Dirichlet eigenfunction is not $\log$-concave. Bourni, Clutterbuck, Nguyen, Stancu, Wei and Wheeler \cite{Bourni-CNSWW2022} subsequently proved that such counterexamples exist in $\mathbb{H}^n$ for every $n\geq 2$. These results show that geodesic convexity is insufficient to recover the corresponding Euclidean conclusions. A natural stronger geometric condition is horo-convexity. A domain $\Omega\subset\mathbb{H}^n$ is called horo-convex if every boundary point admits a supporting horoball containing $\Omega$. For a smooth domain, this is equivalent to requiring that all boundary principal curvatures are at least $1$. If they are all strictly greater than $1$, the domain is called strictly horo-convex. In this setting, Wei and Xiao \cite{Wei-Xiao2025} used a constant rank theorem to establish $\log$-concavity of the first eigenfunction on horo-convex domains of sufficiently small diameter in $\mathbb{H}^n$. In dimension two, Dai, Ennis, Nguyen and Wei \cite{Dai-ENW2026} removed the diameter restriction and further proved that every superlevel set is horo-convex or equivalently that every nontrivial level curve has geodesic curvature at least $1$. Their proof combines Killing fields with a nodal set argument.
	
	Motivated by these developments, we study $\log$-concavity of the first Dirichlet eigenfunction on horo-convex domains in hyperbolic space without a diameter restriction and obtain lower bounds for the curvature of its level sets.
	
	Before stating our results, we briefly review the $P$-function approach to convexity problems in Euclidean space. In 2012, Ma, Shi and Ye \cite{Ma-SY2012} extended the arguments of Makar-Limanov \cite{Makar-Limanov1971} and Acker, Payne and Philippin \cite{Acker-PP1981} to higher dimensions by constructing suitable auxiliary functions and establishing minimum principles for them. More precisely, for the first Dirichlet eigenvalue problem, they considered the auxiliary function $$\varphi=e^{-(n+1)v}\det\nabla^2v,$$ where $v=-\log u$. This function is related to the Gauss curvature of the graph of $v$. Under the assumption $\nabla^2v>0$, they proved that $\varphi$ is superharmonic modulo gradient terms with locally bounded coefficients. They then applied the minimum principle to obtain a lower bound for $\varphi$ in terms of the boundary Gauss curvature and the boundary gradient of $u$. They also used this estimate in a deformation argument to prove $\log$-concavity. Jia, Ma and Shi \cite{Jia-MS2023} further showed that the logarithm of this auxiliary function is superharmonic.
	
	Inspired by this approach, we establish maximum or minimum principles for suitable auxiliary functions in hyperbolic space. We combine these principles with boundary convexity estimates and a domain deformation argument to prove the following results.
	
	From now on, let $\Omega\subset\mathbb{H}^n$ be a bounded smooth domain and let $u$ satisfy \eqref{Equation-Dirichlet problem}. Write $v=-\log u$. For $0<t<\max_\Omega u$, denote the level sets of $u$ by $$\Sigma_t=\{x\in\Omega:u(x)=t\}.$$ At regular points, the principal curvatures of $\Sigma_t$ are taken with respect to the outward normal of the superlevel set $\{u>t\}$. We use $K$ to denote their product, namely the Gauss curvature. In dimension two, $K$ is the geodesic curvature of the level curve.
	
	Our first result concerns the two-dimensional case recently studied by Dai, Ennis, Nguyen and Wei \cite{Dai-ENW2026}. Here we provide a new proof using PDE methods and derive a quantitative lower bound for the Gauss curvature of the level sets in terms of the minimum boundary curvature.
	
	\begin{theorem}\label{Theorem-2-dimensional}
		Assume that $n=2$ and $\Omega$ is strictly horo-convex. Then, $\nabla^2v>0$ in $\Omega$ and the curvature of every level set $\Sigma_t$ satisfies $$K\geq\min_{\partial\Omega}\kappa>1.$$ In particular, every $\Sigma_t$ is strictly horo-convex. Here $\kappa$ denotes the curvature of $\partial\Omega$.
	\end{theorem}
	
	In fact, an approximation argument allows us to weaken the strict horo-convexity assumption on $\Omega$ to horo-convexity while retaining strict $\log$-concavity of the first eigenfunction.
	
	\begin{corollary}\label{Corollary-n=2}
		Assume that $n=2$ and $\Omega$ is horo-convex. Then, $\nabla^2v>0$ in $\Omega$.
	\end{corollary}
	
	We now turn to higher dimensions, where strict $\log$-concavity remains valid. In dimension three, the same Gauss curvature estimate holds as in the two-dimensional case with the minimum Gauss curvature of $\partial\Omega$ as a lower bound. For $n\geq 4$, the product of the two smallest principal curvatures is strictly larger than $1$, so the Gauss curvature also satisfies $K>1$. More precisely, we prove the following result.
	
	\begin{theorem}\label{Theorem-Higher-dimensional}
		Let $n\geq 3$ and assume that $\Omega$ is strictly horo-convex. Then, $\nabla^2v>0$ in $\Omega$. For $n=3$, the Gauss curvature of every level set $\Sigma_t$ satisfies $$K\geq\min_{\partial\Omega}(\kappa_1\kappa_2)>1,$$ where $\kappa_1,\kappa_2$ are the principal curvatures of $\partial\Omega$. For $n\geq 4$, the two smallest principal curvatures $\kappa_1\leq\kappa_2$ of $\Sigma_t$ satisfy $\kappa_1\kappa_2>1$. In particular, $K>1$ on every $\Sigma_t$.
	\end{theorem}
	
	The higher-dimensional counterpart of Corollary \ref{Corollary-n=2} is as follows.
	
	\begin{corollary}\label{Corollary}
		Let $n\geq 3$ and suppose that $\Omega$ is horo-convex. Then, $\nabla^2v>0$ in $\Omega$.
	\end{corollary}
	
	\begin{remark}
		Similar results for the torsion problem \eqref{Equation-Torsion problem} can also be established using similar methods. More precisely, $-\sqrt{u}$ can be shown to be strictly convex on bounded smooth strictly horo-convex domains in $\mathbb{H}^n$ without a diameter restriction. This removes the small diameter assumption in our earlier work \cite{Zhang-Zhou2026}.
	\end{remark}
	
	We now outline the proof. Our argument follows the auxiliary function approach developed by Ma, Shi and Ye \cite{Ma-SY2012}. The key is to construct suitable auxiliary functions and establish maximum or minimum principles for them. We then combine these principles with a domain deformation argument to prove strict $\log$-concavity. For our choice of auxiliary functions $\phi$, a boundary convexity estimate is also essential to keep $\nabla^2v$ positive definite near the boundary throughout the deformation. The ball case establishes $\nabla^2v>0$ and $\phi<1$. Smooth convergence of the eigenfunctions together with the boundary estimates ensures that these inequalities persist on domains sufficiently close to a ball. Proposition \ref{Proposition-Domain deformation} provides the starting point for the continuity argument.
	
	In dimension two, our choice of $\phi$ is motivated by the construction of Dai, Ennis, Nguyen and Wei \cite{Dai-ENW2026}. Its boundary value is $\kappa^{-2}$, where $\kappa$ denotes the curvature of $\partial\Omega$. We establish a maximum principle for $\phi$ and obtain $\phi\leq\max_{\partial\Omega}\kappa^{-2}$. Together with the strict eigenvalue bound $\lambda>\frac{1}{4}$, this estimate rules out the first degeneration of $\nabla^2v$ along the deformation. On each regular level curve, the comparison $K^{-2}\leq\phi$ establishes the curvature estimate.
	
	For the three-dimensional case, the auxiliary function has boundary value $(\kappa_1\kappa_2)^{-1}$, where $\kappa_1,\kappa_2$ are the principal curvatures of $\partial\Omega$. Geometrically, it is related to the area of a parallelogram spanned by an orthonormal tangent pair, with the area measured in the metric defined by the inverse Hessian; see Section \ref{Subsection-Auxiliary functions} for details. A similar maximum principle controls $\phi$ by its boundary values. The strict bound $\lambda>1$ allows us to exclude both possible ranks of a degenerate Hessian. The comparison $K^{-1}\leq\phi$ completes the Gauss curvature estimate.
	
	For $n\geq 4$, we extend this construction by maximizing the same area quantity over all two-dimensional tangent planes of $\Sigma_t$. To control $\phi$ along the deformation, we establish a minimum principle for $e^{-v}(1-\phi)$. This provides a uniform upper bound for $\phi$ strictly below $1$, which rules out Hessian degeneration and proves strict $\log$-concavity. The level set curvature estimate  follows from the same bound.
	
	The rest of this paper is organized as follows. In Section \ref{Section-Preliminaries}, we collect the necessary preliminaries, including basic geometric facts, auxiliary functions, boundary convexity estimates and domain deformation. Sections \ref{Section-n=2}, \ref{Section-n=3} and \ref{Section-n} contain the proofs of our main results in dimensions two, three and $n\geq 4$, respectively.

	\section{Preliminaries}\label{Section-Preliminaries}
	
	This section begins with the geometry of level sets and the hyperbolic identities used below. We then introduce the auxiliary functions and explain their geometric meaning. Finally, we prove the ball and boundary estimates and describe the domain deformation used in the continuity argument.
	
	\subsection{Geometry of level sets}\label{Subsection-Geometry of level sets}
	
	We recall the formulas for the second fundamental form and the curvature of level sets; see Ma and Zhang \cite{Ma-Zhang2014}. Let $g$ denote the hyperbolic metric on $\mathbb{H}^n$ and let $\nabla$ be its Levi-Civita connection. We write $\Braket{\cdot,\cdot}$ and $|\cdot|$ for the inner product and norm induced by $g$, respectively. Unless otherwise stated, Latin indices range from $1$ to $n$ and Greek indices range from $1$ to $n-1$.
	
	For a smooth function $v$, we denote the components of its first and second covariant derivatives in a local orthonormal frame $\{e_1,e_2,\cdots,e_n\}$ by $$v_i=\nabla_{e_i}v,\quad v_{ij}=\nabla^2v(e_i,e_j).$$ By the definition of the Hessian,
	\begin{equation}\label{Equality-Hessian of v}
		v_{ij}=\nabla_{e_i}(\nabla_{e_j}v)-\Braket{\nabla_{e_i}e_j,\nabla v}.
	\end{equation}
	
	Let $u$ be the first Dirichlet eigenfunction satisfying \eqref{Equation-Dirichlet problem} and set $v=-\log u$. For $0<t<\max_\Omega u$, write $$\Sigma_t=\{x\in\Omega:u(x)=t\}=\{x\in\Omega:v(x)=-\log t\}.$$ If $t$ is a regular value of $u$, meaning that $\nabla u(x)\neq 0$ for every $x\in\Sigma_t$, then $\Sigma_t$ is a smooth hypersurface by the implicit function theorem.
	
	Since $\nabla u$ points into the superlevel set $\{u>t\}$ along $\Sigma_t$, the outward unit normal to $\Sigma_t$ is $$e_n=-\frac{\nabla u}{|\nabla u|}=\frac{\nabla v}{|\nabla v|}.$$ Along $\Sigma_t$, we complete $e_n$ to a local orthonormal frame by choosing $e_1,e_2,\cdots,e_{n-1}$ tangent to $\Sigma_t$. As $v$ is constant on $\Sigma_t$, we have
	\begin{equation}\label{Equality-Gradient of v}
		v_\alpha=0,\quad v_n=\Braket{\nabla v,e_n}=|\nabla v|>0.
	\end{equation}
	The second fundamental form of $\Sigma_t$ with respect to $e_n$ is defined by
	\begin{equation}\label{Formula-Second fundamental form}
		h_{\alpha\beta}=\Braket{\nabla_{e_\alpha}e_n,e_\beta}=-\Braket{\nabla_{e_\alpha}e_\beta,e_n}.
	\end{equation}
	Combining \eqref{Equality-Hessian of v}-\eqref{Formula-Second fundamental form}, we obtain $$v_{\alpha\beta}=\nabla_{e_\alpha}(\nabla_{e_\beta}v)-\Braket{\nabla_{e_\alpha}e_\beta,\nabla v}=-|\nabla v|\Braket{\nabla_{e_\alpha}e_\beta,e_n}=|\nabla v|h_{\alpha\beta}.$$ Consequently,
	\begin{equation}\label{Equality-Second fundamental form}
		h_{\alpha\beta}=\frac{v_{\alpha\beta}}{|\nabla v|}.
	\end{equation}
	
	Let $\kappa_1,\kappa_2,\cdots,\kappa_{n-1}$ denote the principal curvatures of $\Sigma_t$, namely the eigenvalues of $(h_{\alpha\beta})$. The level hypersurface is strictly convex if and only if $(h_{\alpha\beta})$ is positive definite. In particular, \eqref{Equality-Second fundamental form} shows that if $\nabla^2v>0$, then every regular level set of $v$ is strictly convex. We use $K$ to denote the product of the principal curvatures, that is, the Gauss curvature of $\Sigma_t$. It follows from \eqref{Equality-Second fundamental form} that
	\begin{equation}\label{Equality-Gauss curvature of level set}
		K=\prod^{n-1}_{\alpha=1}\kappa_\alpha=\det(h_{\alpha\beta})=\frac{\det(v_{\alpha\beta})}{|\nabla v|^{n-1}}.
	\end{equation}
	
	\subsection{Hyperbolic geometry}\label{Subsection-Hyperbolic geometry}
	
	We next record the covariant derivative identities in hyperbolic space $\mathbb{H}^n$ that will be used in the proofs below. In a local orthonormal frame, we use the curvature convention $$R_{ijkl}=-(\delta_{ik}\delta_{jl}-\delta_{il}\delta_{jk}).$$ For a smooth function $u$, we use subscripts to denote covariant derivatives, as in Section \ref{Subsection-Geometry of level sets}. The Ricci commutation identity gives
	\begin{equation}\label{Formula-Ricci commutation}
		u_{ijk}-u_{ikj}=\delta_{ij}u_k-\delta_{ik}u_j.
	\end{equation}
	In particular, for $i\neq j$,
	\begin{equation}\label{Formula-Third derivative commutation}
		u_{jii}=u_{iij}-u_j.
	\end{equation}
	For fourth-order covariant derivatives, we have the commutation formula
	\begin{align*}
		u_{ijkl}-u_{ijlk}
		=&\sum_{m=1}^n(R_{mikl}u_{mj}+R_{mjkl}u_{im})\\
		=&~\delta_{ik}u_{jl}-\delta_{il}u_{jk}+\delta_{jk}u_{il}-\delta_{jl}u_{ik}.
	\end{align*}
	Combining this with \eqref{Formula-Ricci commutation}, we obtain for each fixed $k$,	
	\begin{equation}\label{Formula-Fourth derivative commutation}
		u_{ijkk}=u_{kkij}-2u_{ij}+2\delta_{ij}u_{kk}-\delta_{ik}u_{jk}+\delta_{jk}u_{ik}.
	\end{equation}
	
	For later use, we also recall the radial formulas in geodesic polar coordinates centered at $o\in\mathbb{H}^n$. Writing $r$ for the geodesic distance from $o$, the metric takes the form $$g=\dd r^2+\sinh^2r\dd z^2,$$ where $\dd z^2$ denotes the standard metric on $\mathbb{S}^{n-1}$. For a radial function $u=u(r)$, its Laplacian is
	\begin{equation}\label{Formula-Radial Laplacian}
		\Delta u=u''+(n-1)\coth ru',
	\end{equation}
	and its Hessian is given by
	\begin{equation}\label{Formula-Radial Hessian}
		\nabla^2u=u''\dd r^2+u'\coth r(g-\dd r^2),
	\end{equation}
	where $\dd r^2=\dd r\otimes \dd r$. In an orthonormal frame with $e_n=\partial_r$ and $e_\alpha$ tangent to the geodesic sphere centered at $o$, formula \eqref{Formula-Radial Hessian} gives
	\begin{equation}\label{Equality-Radial Hessian}
		u_{\alpha\beta}=u'\coth r\delta_{\alpha\beta},\quad u_{\alpha n}=0,\quad u_{nn}=u'',
	\end{equation}
	where $1\leq\alpha,\beta\leq n-1$. Thus the radial and tangential eigenvalues of $\nabla^2u$ are $u''$ and $u'\coth r$, respectively, with the latter having multiplicity $n-1$. If $u$ is smooth at $o$, then $u'(0)=0$ and $\nabla^2u(o)=u''(0)g$.
	
	We conclude this subsection with a basic lower bound for the first Dirichlet eigenvalue, which will play a key role in the subsequent arguments.
	\begin{lemma}\label{Lemma-Eigenvalue lower bound}
		Suppose that $u$ satisfies \eqref{Equation-Dirichlet problem} on a bounded, smooth domain $\Omega\subset\mathbb{H}^n$, $n\geq 2$. Then, its first eigenvalue satisfies $$\lambda>\frac{(n-1)^2}{4}.$$
	\end{lemma}
	\begin{proof}
		Choose geodesic polar coordinates with center $o\notin\bar{\Omega}$, so that $r$ is smooth near $\bar{\Omega}$ and $\nabla r=\partial_r$ has unit length. By \eqref{Formula-Radial Laplacian}, we have $$\Delta r=(n-1)\coth r>n-1.$$ Integrating with respect to the volume measure induced by $g$, we use the boundary condition and the Cauchy-Schwarz inequality to obtain
		\begin{align*}
			(n-1)\int_\Omega u^2
			&~<\int_\Omega u^2\Delta r=-2\int_\Omega u\Braket{\nabla u,\nabla r}\\
			&~\leq2\bigg(\int_\Omega u^2\bigg)^{1/2}\bigg(\int_\Omega|\nabla u|^2\bigg)^{1/2}.
		\end{align*}
		On the other hand, multiplying \eqref{Equation-Dirichlet problem} by $u$ and integrating over $\Omega$ yields $$\int_\Omega|\nabla u|^2=\lambda\int_\Omega u^2.$$ Consequently, $$(n-1)\int_\Omega u^2<2\sqrt{\lambda}\int_\Omega u^2.$$ Since $u\not\equiv 0$, the asserted estimate follows.
	\end{proof}
	
	\subsection{Auxiliary functions}\label{Subsection-Auxiliary functions}
	
	We introduce the auxiliary functions used in the proofs of our main results and explain their geometric meaning. Their boundary limits reflect the curvature of $\partial\Omega.$ Throughout, we work in the region where $\nabla^2v>0$.
	
	We first consider the case $n=2$ and define the auxiliary function
	\begin{equation}\label{Auxiliary function-2-dim}
		\phi=\frac{|(\nabla^2v)\nabla v|^2}{(\det\nabla^2v)^2}.
	\end{equation}
	In terms of components, this takes the form $$\phi=\frac{(v_{11}v_1+v_{12}v_2)^2+(v_{12}v_1+v_{22}v_2)^2}{(v_{11}v_{22}-v_{12}^2)^2},$$ which is smooth wherever $\nabla^2v>0$ and vanishes at critical points of $v$. Its relation to the level set curvature can be seen as follows. At a regular point, let $K$ be the curvature of the level curve and choose $e_1$ tangent to the level curve and $e_2=\nabla v/|\nabla v|$. In this frame, $v_1=0$ and $v_2=|\nabla v|$, and hence \eqref{Equality-Second fundamental form} gives 
	\begin{equation}\label{Estimate-Auxiliary function-2-dim}
		\phi=\frac{v_2^2(v_{12}^2+v_{22}^2)}{(v_{11}v_{22}-v_{12}^2)^2}\geq\frac{v_2^2}{v_{11}^2}=\frac{1}{K^2}.
	\end{equation}
	
	Our choice of $\phi$ is also motivated by the work of Dai, Ennis, Nguyen and Wei \cite{Dai-ENW2026} on the horo-convexity of level curves. Let $J=\big(\begin{smallmatrix}
		0 & -1\\
		1 & 0
	\end{smallmatrix}\big)$ denote rotation by $\pi/2$ in a local orthonormal frame. We can equivalently write $$\phi=|(\nabla^2v)^{-1}J\nabla v|^2,$$ using the two-dimensional identity $(\nabla^2v)^{-1}J=J\nabla^2v/\det\nabla^2v$ and the fact that $J$ preserves lengths. In their argument, after proving $\nabla^2v>0$, they use the vector $(\nabla^2v)^{-1}J\nabla v$ to construct suitable Killing fields. If its length exceeds $1$, the construction produces a function satisfying the same eigenvalue equation that vanishes together with its gradient at an interior point. However, horo-convexity ensures that the function has exactly two boundary zeros, which by the planar structure of nodal sets rules out such an interior zero. This contradiction yields $$|(\nabla^2v)^{-1}J\nabla v|\leq 1,$$ and consequently $K\geq 1$, proving the horo-convexity of the level curves.
	
	In our approach, we apply a maximum principle to show that $\phi$ attains its maximum on $\partial\Omega$. Combining this with \eqref{Estimate-Auxiliary function-2-dim} and the boundary identity $\phi=\kappa^{-2}$ from Proposition \ref{Proposition-Boundary behavior}, with $\kappa$ denoting the curvature of $\partial\Omega$, we conclude that the curvature of every regular level curve is bounded below by the minimum curvature of $\partial\Omega$, that is, $$K\geq\min_{\partial\Omega}\kappa.$$
	
	We next consider the case $n=3$, where the regular level sets are surfaces, and define
	\begin{equation}\label{Auxiliary function-3-dim}
		\phi=\frac{(\nabla v)^T(\nabla^2v)\nabla v}{\det\nabla^2v}=\frac{\sum_{i,j=1}^3v_iv_jv_{ij}}{\det(v_{ij})}.
	\end{equation}
	Here, $\nabla v$ is regarded as a column vector and $(\nabla v)^T$ denotes its transpose. As in the two-dimensional case, choose an orthonormal frame with $e_1,e_2$ tangent to a regular level surface and $e_3=\nabla v/|\nabla v|$. The positive definiteness of the Hessian together with \eqref{Equality-Gauss curvature of level set} implies
	\begin{equation}\label{Estimate-Auxiliary function-3-dim}
		\phi=\frac{|\nabla v|^2v_{33}}{\det(v_{ij})}\geq\frac{|\nabla v|^2}{v_{11}v_{22}-v_{12}^2}=\frac{1}{K},
	\end{equation}
	where $K$ is the Gauss curvature of the level surface. The same maximum principle argument combined with the boundary identity $\phi=(\kappa_1\kappa_2)^{-1}$ established in Section \ref{Subsection-The ball case and boundary convexity estimate} yields $$K\geq\min_{\partial\Omega}(\kappa_1\kappa_2),$$ where $\kappa_1,\kappa_2$ denote the principal curvatures of $\partial\Omega$.
	
	The auxiliary function $\phi$ defined in \eqref{Auxiliary function-3-dim} also has a geometric interpretation in terms of tangent areas measured by the inverse Hessian. At each point $p$, the positive definite inverse matrix $(v^{ij})=(v_{ij})^{-1}$ defines an inner product on the tangent space to $\mathbb{H}^3$ at $p$, given by $\sum_{i,j=1}^3v^{ij}X_iY_j$, where $X_i,Y_j$ are the components of tangent vectors $X,Y\in T_p\mathbb{H}^3$ in an orthonormal frame. Let $\xi,\eta$ be local tangent vector fields along $\Sigma_t$ orthonormal with respect to $g$ and orthogonal to $\nabla v$. In terms of their components $\xi_i,\eta_i$, \eqref{Auxiliary function-3-dim} can then be expressed as
	\begin{equation}\label{Auxiliary function-3-dim'}
		\phi=|\nabla v|^2\left[\bigg(\sum_{i,j=1}^3v^{ij}\xi_i\xi_j\bigg)\bigg(\sum_{i,j=1}^3v^{ij}\eta_i\eta_j\bigg)-\bigg(\sum_{i,j=1}^3v^{ij}\xi_i\eta_j\bigg)^2\right].
	\end{equation}
	The expression in brackets is the squared area of the parallelogram spanned by $\xi,\eta$ measured with respect to the metric defined by the inverse Hessian.
	
	The area interpretation of the three-dimensional auxiliary function in \eqref{Auxiliary function-3-dim'} motivates the extension of the construction to higher dimensions. For $n\geq 4$, we maximize the same area quantity over all two-dimensional planes tangent to $\Sigma_t$ and define
	\begin{equation}\label{Auxiliary function-Higher dimensional}
		\phi=|\nabla v|^2\max_{\substack{|\xi|=|\eta|=1,\Braket{\xi,\eta}=0\\ \Braket{\xi,\nabla v}=\Braket{\eta,\nabla v}=0}}\left[\bigg(\sum_{i,j=1}^nv^{ij}\xi_i\xi_j\bigg)\bigg(\sum_{i,j=1}^nv^{ij}\eta_i\eta_j\bigg)-\bigg(\sum_{i,j=1}^nv^{ij}\xi_i\eta_j\bigg)^2\right].
	\end{equation}
	As in the three-dimensional case, this construction provides a curvature estimate. Let $0<\kappa_1\leq\cdots\leq \kappa_{n-1}$ be the principal curvatures of $\Sigma_t$. Choosing $\xi,\eta$ in the principal directions corresponding to $\kappa_1,\kappa_2$, respectively, and using \eqref{Equality-Second fundamental form}, we obtain $$\phi\geq\frac{1}{\kappa_1\kappa_2}.$$
	
	Following the strategy used in dimensions two and three, we combine a minimum principle for $e^{-v}(1-\phi)$ with a domain deformation to establish $\phi<1$ in $\Omega$, using the boundary identity $$\phi=\frac{1}{\kappa_1\kappa_2}<1\quad\mathrm{on}\ \partial\Omega,$$ established in Section \ref{Subsection-The ball case and boundary convexity estimate}, where $\kappa_1,\kappa_2>1$ are the two smallest principal curvatures of $\partial\Omega$. Together with the preceding comparison, this yields $$\kappa_1\kappa_2>1$$ on every regular level hypersurface.
	
	Finally, we define $\phi=0$ at critical points of $v$, so that $\phi$ is continuous wherever $\nabla^2v>0$. Since the maximizing tangent plane may vary from point to point, $\phi$ need not be smooth. For the local calculations below, choose a smooth local orthonormal frame near a regular point $x_0$ such that the tangential projections of the first two directions span a maximizing plane at $x_0$. We orthonormalize these projections to obtain smooth tangent vector fields $\xi,\eta$. Substituting into the expression in \eqref{Auxiliary function-Higher dimensional}, we calculate $$\tilde{\phi}=|\nabla v|^2\frac{\det\big(\sum_{i,j}v^{ij}\big(\delta_{ia}-\frac{v_iv_a}{|\nabla v|^2}\big)\big(\delta_{jb}-\frac{v_jv_b}{|\nabla v|^2}\big)\big)_{a,b=1,2}}{\det\big(\delta_{ab}-\frac{v_av_b}{|\nabla v|^2}\big)_{a,b=1,2}}.$$ Moreover,
	\begin{equation}\label{Auxiliary function-Higher dimensional'}
		\tilde{\phi}=\frac{|\nabla v|^4}{\sum_{i=3}^nv_i^2}\det\bigg(v^{ab}-\frac{v_b}{|\nabla v|^2}\sum_iv^{ai}v_i-\frac{v_a}{|\nabla v|^2}\sum_iv^{bi}v_i+\frac{v_av_b}{|\nabla v|^4}\sum_{i,j}v^{ij}v_iv_j\bigg)_{a,b=1,2}.
	\end{equation}
	By construction, $\tilde{\phi}$ is smooth near $x_0$, agrees with $\phi$ at $x_0$ and is bounded above by $\phi$ nearby.

	\subsection{The ball case and boundary convexity estimate}\label{Subsection-The ball case and boundary convexity estimate}
	
	To prepare for the proofs of the main theorems, we establish two convexity results for the first Dirichlet eigenfunctions, one on geodesic balls and the other near the boundary of horo-convex domains. The proofs are similar to those in \cite{Wei-Xiao2025} and \cite{Zhang-Zhou2026}. We also compute the auxiliary functions introduced in Section \ref{Subsection-Auxiliary functions} on geodesic balls and study their boundary behavior on horo-convex domains.
	
	Let $B_R(o)\subset\mathbb{H}^n$ be the geodesic ball centered at $o$ with radius $R$. We first prove Theorems \ref{Theorem-2-dimensional} and \ref{Theorem-Higher-dimensional} for geodesic balls, providing the results needed for the subsequent deformation argument.
	\begin{proposition}\label{Proposition-Ball case}
		Let $u$ be the first Dirichlet eigenfunction satisfying \eqref{Equation-Dirichlet problem} on $B_R(o)$. Then the function $v=-\log u$ satisfies $\nabla^2v>0$ in $B_R(o)$. Moreover, every level set $\Sigma_t$, $0<t<\max_{B_R(o)}u$, is a geodesic sphere centered at $o$ with all principal curvatures equal to $\coth r$, where $r$ is its radius.
	\end{proposition}
	\begin{proof}
		By separation of variables, the first Dirichlet eigenfunction $u$ on $B_R(o)$ is radial. Thus $v=-\log u$ depends only on $r$, the geodesic distance from $o$. By the radial Laplacian formula \eqref{Formula-Radial Laplacian}, the equation for $v$ becomes
		\begin{equation}\label{Equation-Radial-v}
			v''+(n-1)\coth rv'=\lambda+(v')^2\quad 0<r<R.
		\end{equation}
		Smoothness at the center gives $v'(0)=0$. Taking the limit $r\rightarrow 0$ in \eqref{Equation-Radial-v} and using $$\lim_{r\rightarrow 0}\coth rv'(r)=\lim_{r\rightarrow 0}(r\coth r)\frac{v'(r)}{r}=v''(0),$$ we obtain $v''(0)=\lambda/n>0$.
		
		Multiplying \eqref{Equation-Radial-v} by $\sinh^{n-1}r$ and integrating from $0$ to $r$, we have $$\sinh^{n-1}rv'(r)=\int_0^r\sinh^{n-1}s\big(\lambda+(v'(s))^2\big)\dd s>0.$$ Thus $v'>0$ for $0<r<R$. As noted in Section \ref{Subsection-Hyperbolic geometry}, the radial Hessian formula \eqref{Formula-Radial Hessian} shows that $\nabla^2v$ has $n-1$ tangential eigenvalues equal to $v'\coth r$ and one radial eigenvalue $v''$. It remains to prove that $v''>0$.
		
		Arguing by contradiction, suppose that $v''$ vanishes somewhere in $(0,R)$. Denoting its first zero by $r_0$, one has $$v''(r)>0\quad\mathrm{for}\ 0\leq r<r_0,\quad v''(r_0)=0,\quad v'''(r_0)\leq 0.$$ On the other hand, differentiating \eqref{Equation-Radial-v} gives
		\begin{equation}\label{Equality-Differentiating radial equation}
			v'''=\big(2v'-(n-1)\coth r\big)v''+\frac{n-1}{\sinh^2r}v'.
		\end{equation}
		Evaluating \eqref{Equality-Differentiating radial equation} at $r_0$, we obtain $$v'''(r_0)=\frac{n-1}{\sinh^2r_0}v'(r_0)>0,$$ a contradiction. Hence $v''>0$ on $[0,R)$, and therefore $\nabla^2v>0$ in $B_R(o)$.
		
		For the second assertion, radiality and strict monotonicity of $v$ yield $\Sigma_t=\partial B_{r_t}(o)$ for some $r_t\in(0,R)$. By \eqref{Formula-Radial Hessian}, in an orthonormal frame tangent to $\Sigma_t$, we have $$h_{\alpha\beta}=\frac{v_{\alpha\beta}}{|\nabla v|}=\frac{v'(r_t)\coth r_t\delta_{\alpha\beta}}{v'(r_t)}=\coth r_t\delta_{\alpha\beta}.$$ Therefore, all its principal curvatures equal $\coth r_t$, completing the proof.
	\end{proof}
	
	On $B_R(o)$, we can explicitly evaluate the auxiliary functions introduced in Section \ref{Subsection-Auxiliary functions}. Applying the radial Hessian formula \eqref{Equality-Radial Hessian} to $v$, we obtain, in an orthonormal frame with $e_n=\partial_r$ and $e_1,e_2,\cdots,e_{n-1}$ tangent to the geodesic sphere centered at $o$, $$v_n=v',\quad v_{\alpha\beta}=v'\coth r\delta_{\alpha\beta},\quad v_{\alpha n}=0,\quad v_{nn}=v''$$ for $0<r<R$. Substituting these identities into \eqref{Auxiliary function-2-dim} and \eqref{Auxiliary function-3-dim} gives $\phi=\tanh^2r$ when $n=2,3$. For $n\geq 4$, we have $v^{\alpha\beta}=(v'\coth r)^{-1}\delta_{\alpha\beta}$, and hence, for every orthonormal tangent pair $\xi,\eta$, $$\sum_{i,j=1}^nv^{ij}\xi_i\xi_j=\sum_{i,j=1}^nv^{ij}\eta_i\eta_j=\frac{1}{v'\coth r},\quad \sum_{i,j=1}^nv^{ij}\xi_i\eta_j=0.$$ The definition of $\phi$ in \eqref{Auxiliary function-Higher dimensional} consequently yields $\phi=\tanh^2r$. With $\phi(o)=0$ at the center, we obtain in every dimension
	\begin{equation}\label{Estimate-Auxiliary function-ball case}
		0\leq\phi<\tanh^2R<1\quad\mathrm{in}\ B_R(o).
	\end{equation}

	We next show that $v=-\log u$ is strictly convex near the boundary of horo-convex domains.
	
	\begin{proposition}\label{Proposition-Boundary convexity estimate}
		Let $\Omega$ be a bounded, smooth horo-convex domain in $\mathbb{H}^n$ for $n\geq 2$, and let $u$ be the first Dirichlet eigenfunction satisfying \eqref{Equation-Dirichlet problem}. Then, for $v=-\log u$, there exists a neighborhood $\mathcal{N}$ of $\partial\Omega$ such that $\nabla^2v>0$ in $\mathcal{N}\cap\Omega$.
	\end{proposition}
	\begin{proof}
		Let $\nu$ denote the unit outward normal to $\partial\Omega$. By the maximum principle and Hopf's lemma, we have $u>0$ in $\Omega$ and $\nabla_\nu u<0$ on $\partial\Omega$. Since $u=0$ on $\partial\Omega$, all its tangential derivatives vanish there, and hence $$\nabla u=(\nabla_\nu u)\nu=-|\nabla u|\nu\quad\mathrm{on}\ \partial\Omega.$$
		
		Choose an orthonormal frame $\{e_1,e_2,\cdots,e_{n-1},\nu\}$ along $\partial\Omega$ and denote its second fundamental form by $h_{\alpha\beta}$. On $\partial\Omega$, we have $$u_{\alpha\beta}=e_\alpha(e_\beta u)-\Braket{\nabla_{e_\alpha}e_\beta,\nabla u}=h_{\alpha\beta}\nabla_\nu u=-|\nabla u|h_{\alpha\beta}.$$ Thus, for every unit tangential vector $e$, horo-convexity implies
		\begin{equation}\label{Estimate-Boundary-Hessian of u}
			-\nabla^2u(e,e)=|\nabla u|h(e,e)\geq|\nabla u|>0.
		\end{equation}
		
		Since $v=-\log u$, direct computation in $\Omega$ yields $$\nabla v=-\frac{\nabla u}{u},\quad \nabla^2v=\frac{\nabla u\otimes\nabla u}{u^2}-\frac{\nabla^2u}{u}.$$ Therefore, for any unit vector $e$,
		\begin{equation}\label{Equality-Boundary-Hessian of v}
			\nabla^2v(e,e)=\frac{|\nabla_eu|^2}{u^2}-\frac{\nabla^2u(e,e)}{u}.
		\end{equation}
		Near $\partial\Omega$, write $e=e^\top+e^\bot$, where $e^\bot=\Braket{e,\nabla u/|\nabla u|}\nabla u/|\nabla u|$. Then $$|\nabla_eu|^2=\Braket{e,\nabla u}^2=\left\langle e,|\nabla u|\frac{\nabla u}{|\nabla u|}\right\rangle^2=|\nabla u|^2|e^\bot|^2,$$ and substituting into \eqref{Equality-Boundary-Hessian of v}, we obtain $$\nabla^2v(e,e)=\frac{|\nabla u|^2|e^\bot|^2}{u^2}-\frac{\nabla^2u(e,e)}{u}.$$ By \eqref{Estimate-Boundary-Hessian of u}, compactness and continuity, there exists $\eta>0$ such that $\nabla^2v(e,e)>0$ near $\partial\Omega$ for every unit vector $e$ satisfying $|e^\bot|\leq\eta$. Since $u\rightarrow 0$ while $|\nabla u|$ is bounded away from zero by Hopf's lemma and $\nabla^2u$ is bounded near $\partial\Omega$, we have $\nabla^2v(e,e)\rightarrow+\infty$ uniformly for $|e^\bot|\geq\eta$ as $x\rightarrow\partial\Omega$. Hence $\nabla^2v>0$ in $\mathcal{N}\cap\Omega$ for a sufficiently small neighborhood $\mathcal{N}$ of $\partial\Omega$, which completes the proof.
	\end{proof}
	
	The preceding boundary convexity estimate guarantees that the auxiliary functions introduced in Section \ref{Subsection-Auxiliary functions} are well defined near $\partial\Omega$. We now determine their boundary values.
	\begin{proposition}\label{Proposition-Boundary behavior}
		Let $\Omega\subset\mathbb{H}^n$, $n\geq 2$, be a bounded, smooth and horo-convex domain, and let $u$ be its first Dirichlet eigenfunction satisfying \eqref{Equation-Dirichlet problem}. Set $v=-\log u$ and let $\phi$ be the corresponding auxiliary function defined in \eqref{Auxiliary function-2-dim}, \eqref{Auxiliary function-3-dim} or \eqref{Auxiliary function-Higher dimensional}. Then $\phi$ extends continuously to $\partial\Omega$. In dimension two,
		\begin{equation}\label{Boundary behavior-2-dim}
			\lim_{x\rightarrow p}\phi(x)=\frac{1}{\kappa(p)^2},\quad p\in\partial\Omega,
		\end{equation}
		where $\kappa$ is the boundary curvature. For $n\geq 3$,
		\begin{equation}\label{Boundary behavior-higher dim}
			\lim_{x\rightarrow p}\phi(x)=\frac{1}{\kappa_1(p)\kappa_2(p)},\quad p\in\partial\Omega,
		\end{equation}
		where $\kappa_1(p)\leq\kappa_2(p)$ are the two smallest principal curvatures of $\partial\Omega$ at $p$.
	\end{proposition}
	\begin{proof}
		Fix $p\in\partial\Omega$. By Hopf's lemma and continuity, $\nabla u\neq 0$ near $p$. Choose a local orthonormal frame with $e_n=-\nabla u/|\nabla u|$ and $e_1,e_2,\cdots,e_{n-1}$ tangent to the level hypersurface of $u$, so that $u_\alpha=0$ and $u_n=-|\nabla u|$. At $p$, choose the tangential vectors in the principal directions of $\partial\Omega$, ordered so that $$1\leq\kappa_1(p)\leq\kappa_2(p)\leq\cdots\leq\kappa_{n-1}(p).$$
		
		On $\partial\Omega$, $e_n$ is the outward unit normal, and \eqref{Estimate-Boundary-Hessian of u} implies
		\begin{equation}\label{Equality-Boundary-tangential Hessian of u}
			-u_{\alpha\beta}(p)=|\nabla u(p)|h_{\alpha\beta}(p)=|\nabla u(p)|\kappa_\alpha(p)\delta_{\alpha\beta}.
		\end{equation}
		Since $v=-\log u$ and $u_\alpha=0$ in the chosen frame, $$v_{\alpha\beta}=-\frac{u_{\alpha\beta}}{u},\quad v_{\alpha n}=-\frac{u_{\alpha n}}{u},\quad v_{nn}=\frac{u_n^2-uu_{nn}}{u^2}.$$ In dimension two, $v_1=0$ and $v_2=-u_2/u$, so \eqref{Auxiliary function-2-dim} becomes $$\phi=\frac{v_2^2(v_{12}^2+v_{22}^2)}{(v_{11}v_{22}-v_{12}^2)^2}=\frac{u_2^2\big(u^2u_{12}^2+(u_2^2-uu_{22})^2\big)}{\big(u_{11}(u_2^2-uu_{22})+uu_{12}^2\big)^2}.$$ As $x\rightarrow p$, we have $u\rightarrow 0$, $u_2(p)<0$ and $u_{11}(p)=\kappa(p)u_2(p)$. Therefore, $$\lim_{x\rightarrow p}\phi(x)=\frac{u_2(p)^2}{u_{11}(p)^2}=\frac{1}{\kappa(p)^2}.$$
		
		We now consider $n\geq 3$. A direct calculation from $\sum_{k=1}^nv_{ik}v^{kj}=\delta_{ij}$ and $|\nabla v|=|\nabla u|/u$ shows that $$(|\nabla v|v^{\alpha\beta})_{\alpha,\beta=1}^{n-1}=|\nabla u|\bigg(\bigg(-u_{\alpha\beta}-\frac{uu_{\alpha n}u_{\beta n}}{u_n^2-uu_{nn}}\bigg)_{\alpha,\beta=1}^{n-1}\bigg)^{-1}.$$ Taking $x\rightarrow p$ and recalling that $u(p)=0$ and $u_n(p)\neq 0$, we conclude from \eqref{Equality-Boundary-tangential Hessian of u} that $$\lim_{x\rightarrow p}(|\nabla v|v^{\alpha\beta})_{\alpha,\beta=1}^{n-1}=\mathrm{diag}\bigg(\frac{1}{\kappa_1(p)},\frac{1}{\kappa_2(p)},\cdots,\frac{1}{\kappa_{n-1}(p)}\bigg).$$ In particular, each component satisfies
		\begin{equation}\label{Equality-Boundary-scaled inverse Hessian}
			\lim_{x\rightarrow p}(|\nabla v|v^{\alpha\beta})=\frac{\delta_{\alpha\beta}}{\kappa_\alpha(p)}.
		\end{equation}
		In dimension three, choosing $\xi=e_1$ and $\eta=e_2$ reduces \eqref{Auxiliary function-3-dim'} to $$\phi=(|\nabla v|v^{11})(|\nabla v|v^{22})-(|\nabla v|v^{12})^2.$$ It follows from \eqref{Equality-Boundary-scaled inverse Hessian} that $$\lim_{x\rightarrow p}\phi(x)=\frac{1}{\kappa_1(p)\kappa_2(p)}.$$ For $n\geq 4$, the convergence in \eqref{Equality-Boundary-scaled inverse Hessian} is uniform for all orthogonal unit tangent vectors $\xi$ and $\eta$. Hence, by \eqref{Auxiliary function-Higher dimensional}, $$\lim_{x\rightarrow p}\phi(x)=\max_{1\leq\alpha<\beta\leq n-1}\frac{1}{\kappa_\alpha(p)\kappa_\beta(p)}=\frac{1}{\kappa_1(p)\kappa_2(p)},$$ where the maximum is attained in the first two principal directions. Since the boundary values depend continuously on $p$, $\phi$ extends continuously to $\partial\Omega$. This completes the proof.
	\end{proof}

	\subsection{Domain deformation}\label{Subsection-Domain deformation}
	
	To establish convexity, we use a continuity argument that relies on a deformation of the domain. We therefore need a smooth family of strictly horo-convex domains connecting the given domain to a geodesic ball. Following \cite{Wei-Xiao2025}, we describe this deformation and establish the inequalities needed to begin the continuity argument.
	
	Let $\Omega\subset\mathbb{H}^n$, $n\geq 2$, be a bounded, smooth and strictly horo-convex domain. For $n\geq 3$, strict horo-convexity ensures that $\partial\Omega$ is an admissible initial hypersurface for the locally constrained curvature flow of Hu, Li and Wei \cite{Hu-LW2022}. The flow produces a smooth family of strictly horo-convex domains $\{\Omega_t\}_{t\geq 0}$ with $\Omega_0=\Omega$, whose boundaries converge smoothly to a geodesic sphere of radius $R>0$ centered at $o\in\mathbb{H}^n$, $$\partial\Omega_t\rightarrow\partial B_R(o)\quad \mathrm{as}\ t\rightarrow\infty.$$ For $n=2$, the corresponding family follows from the dimension reduction argument in \cite{Wei-Xiao2025}.
	
	Along this deformation, let $u_t$ be the first Dirichlet eigenfunction satisfying \eqref{Equation-Dirichlet problem} on $\Omega_t$, normalized by $\int_{\Omega_t}u_t^2=1$, where integration is with respect to the volume measure induced by $g$. Set $v_t=-\log u_t$ and let $\phi_t$ be the corresponding auxiliary function from Section \ref{Subsection-Auxiliary functions}, defined wherever $\nabla^2v_t>0$. By standard elliptic theory, the normalized eigenfunctions $u_t$ depend smoothly on $t$ and converge to the normalized first Dirichlet eigenfunction $u_\infty$ on $B_R(o)$ with smooth convergence up to the boundary. In particular, $v_t$ converges locally smoothly to $v_\infty=-\log u_\infty$ in the interior.
	
	The convergence above together with the ball result and the boundary estimates allows us to begin the continuity argument through the following proposition.
	\begin{proposition}\label{Proposition-Domain deformation}
		There exists $T>0$ such that $\nabla^2v_t>0$ and $\phi_t<1$ in $\Omega_t$ for every $t\geq T$.
	\end{proposition}
	\begin{proof}
		Proposition \ref{Proposition-Ball case} and \eqref{Estimate-Auxiliary function-ball case} establish $\nabla^2v_\infty>0$ and $\phi_\infty<1$ on $B_R(o)$, where $\phi_\infty$ is the auxiliary function corresponding to $v_\infty$. Since $v_t$ converges locally smoothly to $v_\infty$, the same conclusions hold on each compact subset of $B_R(o)$ for all sufficiently large $t$.
		
		Meanwhile, smooth convergence of the domains and eigenfunctions provides uniform derivative bounds for $u_t$ and a positive lower bound for $|\nabla u_t|$ on $\partial\Omega_t$ for all sufficiently large $t$. These bounds allow us to apply the boundary convexity argument for $v_t=-\log u_t$ in Proposition \ref{Proposition-Boundary convexity estimate} uniformly, so that $\nabla^2v_t>0$ in a boundary neighborhood of uniform width. Moreover, \eqref{Boundary behavior-2-dim} and \eqref{Boundary behavior-higher dim} together with the convergence of the boundary principal curvatures to $\coth R>1$ imply that the boundary values of $\phi_t$ converge uniformly to $\tanh^2R<1$. The uniform bounds above and the formulas in the proof of Proposition \ref{Proposition-Boundary behavior} ensure that $\phi_t$ approaches its boundary values uniformly in $t$. Hence, $\phi_t<1$ in a boundary neighborhood of uniform width.
		
		Finally, the interior and boundary conclusions together show that there exists $T>0$ such that $\nabla^2v_t>0$ and $\phi_t<1$ throughout $\Omega_t$ for every $t\geq T$, completing the proof.
	\end{proof}
	
	\section{The two-dimensional case}\label{Section-n=2}
	
	We begin with the two-dimensional case and prove Theorem \ref{Theorem-2-dimensional}, illustrating the general strategy for our main results. Let $u$ be the first Dirichlet eigenfunction satisfying \eqref{Equation-Dirichlet problem} on a bounded, smooth and strictly horo-convex domain $\Omega\subset\mathbb{H}^2$. The function $v=-\log u$ solves
	\begin{equation}\label{Equation-v}
		\Delta v=\lambda+|\nabla v|^2\quad\mathrm{in}\ \Omega.
	\end{equation}
	We use the auxiliary function $\phi$ defined in \eqref{Auxiliary function-2-dim} and write $\kappa$ for the boundary curvature. The proof proceeds by first establishing a maximum principle for $\phi$ under the assumptions $\nabla^2v>0$ and $\phi<1$. We then apply the domain deformation from Section \ref{Subsection-Domain deformation} and rule out the first degeneration of the Hessian, thereby establishing strict $\log$-concavity. Moreover, the level set curvature estimate follows from this maximum principle.
	
	\subsection{A maximum principle}
	
	The following proposition bounds $\phi$ in terms of the boundary curvature.
	
	\begin{proposition}\label{Proposition-Maximum principle-2-dim}
		Suppose that $\nabla^2v>0$ and $\phi<1$ in $\Omega$. Then
		\begin{equation}\label{Estimate-Boundary-2-dim}
			\phi\leq \max_{\partial\Omega}\kappa^{-2}<1.
		\end{equation}
	\end{proposition}
	\begin{proof}
		In view of \eqref{Boundary behavior-2-dim}, $\phi$ extends continuously to $\bar{\Omega}$ with positive boundary values $\kappa^{-2}$. To obtain the desired estimate \eqref{Estimate-Boundary-2-dim}, it suffices to show that $\phi$ attains its maximum on $\partial\Omega$. Suppose to the contrary that \eqref{Estimate-Boundary-2-dim} fails. Then, $\phi$ attains its maximum at some $x_0\in\Omega$, where $0<\phi(x_0)<1$, $\nabla\phi(x_0)=0$ and $\Delta\phi(x_0)\leq 0$.
		
		Note that $\nabla|\nabla v|^2=2(\nabla^2v)\nabla v$, so $\phi(x_0)>0$ implies $\nabla|\nabla v|^2(x_0)\neq 0$. Choose orthonormal coordinates at $x_0$ whose second coordinate vector points in the direction of $\nabla|\nabla v|^2$. All the following calculations are carried out at $x_0$, and all indices range from $1$ to $2$. By the choice of coordinates,
		\begin{equation}\label{2-dim-Coordinate condition}
			v_1v_{11}+v_2v_{12}=0,
		\end{equation}
		so the positive definiteness of $\nabla^2v$ yields
		\begin{equation}\label{2-dim-Coordinate condition'}
			v_1=-\frac{v_{12}}{v_{11}}v_2,\quad v_2>0.
		\end{equation}
		Consequently,
		\begin{equation}\label{2-dim-Coordinate condition''}
			v_1v_{12}+v_2v_{22}=\frac{v_2}{v_{11}}(v_{11}v_{22}-v_{12}^2).
		\end{equation}
		From \eqref{2-dim-Coordinate condition} and \eqref{2-dim-Coordinate condition''}, the expression for $\phi$ reduces to
		\begin{equation}\label{2-dim-Condition-phi}
			\phi=\frac{v_2^2}{v_{11}^2}.
		\end{equation}
		
		Differentiating the definition of $\phi$ in \eqref{Auxiliary function-2-dim} and using $\nabla\phi(x_0)=0$, we find
		\begin{equation}\label{2-dim-1st condition}
			\begin{aligned}
				0=\frac{\phi_k}{2\phi}
				=&~\frac{v_1v_{12k}+v_2v_{22k}+v_{12}v_{1k}+v_{22}v_{2k}}{v_1v_{12}+v_2v_{22}}\\
				&-\frac{v_{22}v_{11k}+v_{11}v_{22k}-2v_{12}v_{12k}}{v_{11}v_{22}-v_{12}^2}.
			\end{aligned}
		\end{equation}
		The relations \eqref{2-dim-Coordinate condition'} and \eqref{2-dim-Coordinate condition''} simplify \eqref{2-dim-1st condition} to
		\begin{equation}\label{2-dim-1st condition'}
			v_{22}v_{11k}-v_{12}v_{12k}=\frac{v_{11}}{v_2}(v_{12}v_{1k}+v_{22}v_{2k}).
		\end{equation}
		Taking derivatives once more, summing over $k$ and using $\phi_k=0$, we have
		\begin{equation}\label{2-dim-Equality-F1}
			\begin{aligned}
				\frac{\Delta\phi}{2\phi}
				=&-\frac{1}{v_{11}v_{22}-v_{12}^2}\sum_k(v_{22}v_{11kk}+v_{11}v_{22kk}-2v_{12}v_{12kk}+2v_{11k}v_{22k}-2v_{12k}^2)\\
				&+\frac{1}{(v_1v_{12}+v_2v_{22})^2}\sum_k(v_1v_{11k}+v_2v_{12k}+v_{11}v_{1k}+v_{12}v_{2k})^2\\
				&+\frac{1}{v_1v_{12}+v_2v_{22}}\sum_k(v_1v_{12kk}+v_2v_{22kk}+2v_{12k}v_{1k}+2v_{22k}v_{2k}\\
				&\hspace{3.4cm}+v_{12}v_{1kk}+v_{22}v_{2kk}).
			\end{aligned}
		\end{equation}
		
		To handle the fourth-order derivative terms in \eqref{2-dim-Equality-F1}, we differentiate the equation \eqref{Equation-v} and use \eqref{Formula-Third derivative commutation} to write
		\begin{equation}\label{2-dim-3rd derivative term}
			\sum_kv_{ikk}=(\Delta v)_i-v_i=2\sum_kv_kv_{ki}-v_i.
		\end{equation}
		The second derivatives of \eqref{Equation-v} together with \eqref{Formula-Ricci commutation} yield
		\begin{align*}
			(\Delta v)_{ij}
			=&~2\sum_kv_{ki}v_{kj}+2\sum_kv_kv_{kij}\\
			=&~2\sum_kv_{ki}v_{kj}+2\sum_kv_kv_{ijk}+2v_iv_j-2\delta_{ij}\sum_kv_k^2.
		\end{align*}
		It follows from \eqref{Formula-Ricci commutation} and \eqref{Formula-Fourth derivative commutation} that
		\begin{equation}\label{2-dim-4th derivative term}
			\begin{aligned}
				\sum_kv_{ijkk}
				=&~(\Delta v)_{ij}-4v_{ij}+2\delta_{ij}\Delta v\\
				=&~2\sum_kv_kv_{ijk}+2\sum_kv_{ki}v_{kj}+2v_iv_j-4v_{ij}+2\lambda\delta_{ij}.
			\end{aligned}
		\end{equation}
		Using the relations \eqref{2-dim-Coordinate condition'} and \eqref{2-dim-Coordinate condition''}, we calculate the fourth-order derivative terms in \eqref{2-dim-Equality-F1} as
		\begin{align*}
			&-\frac{1}{v_{11}v_{22}-v_{12}^2}\sum_k(v_{22}v_{11kk}+v_{11}v_{22kk}-2v_{12}v_{12kk})\\
			&+\frac{1}{v_1v_{12}+v_2v_{22}}\sum_k(v_1v_{12kk}+v_2v_{22kk})\\
			=&-\frac{1}{v_{11}v_{22}-v_{12}^2}\sum_k(v_{22}v_{11kk}-v_{12}v_{12kk})\\
			=&-\frac{2}{v_{11}v_{22}-v_{12}^2}\sum_kv_k(v_{22}v_{11k}-v_{12}v_{12k})\\
			&-\frac{2}{v_{11}v_{22}-v_{12}^2}\sum_k(v_{22}v_{1k}^2-v_{12}v_{1k}v_{2k})\\
			&-\frac{2(v_1^2v_{22}-v_1v_2v_{12})}{v_{11}v_{22}-v_{12}^2}+4-\frac{2\lambda v_{22}}{v_{11}v_{22}-v_{12}^2},
		\end{align*}
		where \eqref{2-dim-4th derivative term} has been used in the last equality. Applying \eqref{2-dim-Coordinate condition''}, we rewrite \eqref{2-dim-Equality-F1} as
		\begin{equation}\label{2-dim-Equality-F2}
			\begin{aligned}
				\frac{\Delta\phi}{2\phi}
				=&-\frac{2}{v_{11}v_{22}-v_{12}^2}\sum_kv_k(v_{22}v_{11k}-v_{12}v_{12k})\\
				&-\frac{2}{v_{11}v_{22}-v_{12}^2}\sum_k(v_{11k}v_{22k}-v_{12k}^2)\\
				&+\frac{v_{11}^2}{v_2^2(v_{11}v_{22}-v_{12}^2)^2}\sum_k(v_1v_{11k}+v_2v_{12k}+v_{11}v_{1k}+v_{12}v_{2k})^2\\
				&+\frac{2v_{11}}{v_2(v_{11}v_{22}-v_{12}^2)}\sum_k(v_{1k}v_{12k}+v_{2k}v_{22k})\\
				&+\frac{v_{11}}{v_2(v_{11}v_{22}-v_{12}^2)}\sum_k(v_{12}v_{1kk}+v_{22}v_{2kk})\\
				&-\frac{2}{v_{11}v_{22}-v_{12}^2}\sum_k(v_{22}v_{1k}^2-v_{12}v_{1k}v_{2k})\\
				&-\frac{2(v_1^2v_{22}-v_1v_2v_{12})}{v_{11}v_{22}-v_{12}^2}+4-\frac{2\lambda v_{22}}{v_{11}v_{22}-v_{12}^2}.
			\end{aligned}
		\end{equation}
		
		We next deal with the third-order derivative terms in \eqref{2-dim-Equality-F2}. The relations \eqref{2-dim-1st condition'} and \eqref{2-dim-3rd derivative term} imply
		\begin{align*}
			&-\frac{2}{v_{11}v_{22}-v_{12}^2}\sum_kv_k(v_{22}v_{11k}-v_{12}v_{12k})\\
			=&-\frac{2v_{11}}{v_2(v_{11}v_{22}-v_{12}^2)}\sum_kv_k(v_{12}v_{1k}+v_{22}v_{2k})
		\end{align*}
		and
		\begin{align*}
			&~\frac{v_{11}}{v_2(v_{11}v_{22}-v_{12}^2)}\sum_k(v_{12}v_{1kk}+v_{22}v_{2kk})\\
			=&~\frac{2v_{11}}{v_2(v_{11}v_{22}-v_{12}^2)}\sum_kv_k(v_{12}v_{1k}+v_{22}v_{2k})-\frac{v_{11}(v_1v_{12}+v_2v_{22})}{v_2(v_{11}v_{22}-v_{12}^2)},
		\end{align*}
		respectively. Thus,  \eqref{2-dim-Coordinate condition''} shows that
		\begin{equation}\label{2-dim-3rd derivative term1}
			\begin{aligned}
				&-\frac{2}{v_{11}v_{22}-v_{12}^2}\sum_kv_k(v_{22}v_{11k}-v_{12}v_{12k})\\
				&+\frac{v_{11}}{v_2(v_{11}v_{22}-v_{12}^2)}\sum_k(v_{12}v_{1kk}+v_{22}v_{2kk})\\
				=&-\frac{v_{11}(v_1v_{12}+v_2v_{22})}{v_2(v_{11}v_{22}-v_{12}^2)}\\
				=&-1.
			\end{aligned}
		\end{equation}
		For the quadratic terms in the third derivatives in \eqref{2-dim-Equality-F2}, we use \eqref{2-dim-1st condition'} and the positive definiteness of $\nabla^2v$ to express $v_{11k}$ in terms of $v_{12k}$, namely $$v_{11k}=\frac{v_{12}}{v_{22}}v_{12k}+\frac{v_{11}}{v_2v_{22}}(v_{12}v_{1k}+v_{22}v_{2k}).$$ It follows that
		\begin{equation}\label{2-dim-3rd derivative term2}
			v_1v_{11k}+v_2v_{12k}+v_{11}v_{1k}+v_{12}v_{2k}=\frac{v_{11}v_{22}-v_{12}^2}{v_{22}}\bigg(v_{1k}+\frac{v_2}{v_{11}}v_{12k}\bigg)
		\end{equation}
		by \eqref{2-dim-Coordinate condition'} and
		\begin{equation}\label{2-dim-3rd derivative term3}
			\begin{aligned}
				&-v_{11k}v_{22k}+v_{12k}^2+\frac{v_{11}}{v_2}(v_{1k}v_{12k}+v_{2k}v_{22k})\\
				=&~\frac{v_{11}}{v_2v_{22}}(v_{22}v_{12k}-v_{12}v_{22k})\bigg(v_{1k}+\frac{v_2}{v_{11}}v_{12k}\bigg).
			\end{aligned}
		\end{equation}
		Combining \eqref{2-dim-3rd derivative term1}-\eqref{2-dim-3rd derivative term3} with \eqref{2-dim-Equality-F2}, we conclude from \eqref{2-dim-Condition-phi} that
		\begin{equation}\label{2-dim-Equality-F3}
			\begin{aligned}
				\frac{\Delta\phi}{2\phi}
				=&~\frac{1}{\phi v_{22}^2}\sum_k\bigg(v_{1k}+\frac{v_2}{v_{11}}v_{12k}\bigg)^2\\
				&+\frac{2v_{11}}{v_2v_{22}(v_{11}v_{22}-v_{12}^2)}\sum_k(v_{22}v_{12k}-v_{12}v_{22k})\bigg(v_{1k}+\frac{v_2}{v_{11}}v_{12k}\bigg)\\
				&-\frac{2}{v_{11}v_{22}-v_{12}^2}\sum_k(v_{22}v_{1k}^2-v_{12}v_{1k}v_{2k})\\
				&-\frac{2\lambda v_{22}}{v_{11}v_{22}-v_{12}^2}-\frac{2(v_1^2v_{22}-v_1v_2v_{12})}{v_{11}v_{22}-v_{12}^2}+3.
			\end{aligned}
		\end{equation}
		
		We first express the third derivatives in \eqref{2-dim-Equality-F3} in terms of $v_{121}$ and $v_{122}$. Differentiating \eqref{Equation-v} with respect to $x_1$ and using \eqref{2-dim-Coordinate condition}, we obtain $$v_{111}+v_{221}=0.$$ For differentiation in the $x_2$ direction, \eqref{2-dim-Coordinate condition''} leads to $$v_{112}+v_{222}=\frac{2v_2}{v_{11}}(v_{11}v_{22}-v_{12}^2).$$ Applying \eqref{Formula-Third derivative commutation} and substituting the expression for $v_1$ from \eqref{2-dim-Coordinate condition'}, we rewrite these identities as
		\begin{equation}\label{2-dim-3rd derivative identities}
			\begin{aligned}
				v_{112}&=v_{121}+v_2,\\[8pt]
				v_{111}&=-v_{221}=-v_{122}+\frac{v_2v_{12}}{v_{11}},\\
				v_{222}&=\frac{2v_2}{v_{11}}(v_{11}v_{22}-v_{12}^2)-v_{121}-v_2.
			\end{aligned}
		\end{equation}
		We now determine $v_{121}$ and $v_{122}$ from \eqref{2-dim-1st condition'}. It follows directly from \eqref{2-dim-3rd derivative identities} and \eqref{2-dim-Condition-phi} that
		\begin{align*}
			v_{12}v_{121}+v_{22}v_{122}=&-\frac{v_{11}v_{12}}{v_2}\big(v_{11}+(1-\phi)v_{22}\big),\\
			v_{22}v_{121}-v_{12}v_{122}=&~\frac{v_{11}}{v_2}(v_{12}^2+v_{22}^2-\phi v_{11}v_{22}).
		\end{align*}
		Since $v_{12}^2+v_{22}^2>0$, the resulting expressions are
		\begin{align*}
			v_{121}=&~\frac{v_{11}}{v_2}\bigg(v_{22}-v_{11}+(1-\phi)\frac{v_{22}(v_{11}v_{22}-v_{12}^2)}{v_{12}^2+v_{22}^2}\bigg)
			\shortintertext{and}
			v_{122}=&-\frac{v_{11}v_{12}}{v_2}\bigg(2-\phi+(1-\phi)\frac{v_{11}v_{22}-v_{12}^2}{v_{12}^2+v_{22}^2}\bigg).
		\end{align*}
		Taking $k=1,2$ in the first two terms on the right-hand side of \eqref{2-dim-Equality-F3}, we find
		\begin{align*}
			\frac{1}{v_{22}}\bigg(v_{11}+\frac{v_2}{v_{11}}v_{121}\bigg)=&~1+(1-\phi)\frac{v_{11}v_{22}-v_{12}^2}{v_{12}^2+v_{22}^2},\\
			\frac{1}{v_{22}}\bigg(v_{12}+\frac{v_2}{v_{11}}v_{122}\bigg)=&-(1-\phi)\frac{v_{12}(v_{11}+v_{22})}{v_{12}^2+v_{22}^2},\\
			v_{22}v_{121}-v_{12}v_{221}=&~\frac{v_{11}}{v_2}\big(v_{12}^2+v_{22}^2-\phi(v_{11}v_{22}-v_{12}^2)\big)
			\shortintertext{and}
			v_{22}v_{122}-v_{12}v_{222}=&-(1-\phi)\frac{v_{11}v_{12}(v_{11}+v_{22})}{v_2}-\frac{2v_2v_{12}(v_{11}v_{22}-v_{12}^2)}{v_{11}}.
		\end{align*}
		Substituting these identities into \eqref{2-dim-Equality-F3}, we arrive at
		\begin{equation}\label{2-dim-Equality-F4}
			\begin{aligned}
				\frac{\Delta\phi}{2\phi}
				=&~\frac{1}{\phi}\bigg(\bigg(1+(1-\phi)\frac{v_{11}v_{22}-v_{12}^2}{v_{12}^2+v_{22}^2}\bigg)^2+(1-\phi)^2\frac{v_{12}^2(v_{11}+v_{22})^2}{(v_{12}^2+v_{22}^2)^2}\bigg)\\
				&+\frac{2}{\phi(v_{11}v_{22}-v_{12}^2)}\bigg(\big(v_{12}^2+v_{22}^2-\phi(v_{11}v_{22}-v_{12}^2)\big)\\
				&\hspace{3cm}\cdot\bigg(1+(1-\phi)\frac{v_{11}v_{22}-v_{12}^2}{v_{12}^2+v_{22}^2}\bigg)+(1-\phi)^2\frac{v_{12}^2(v_{11}+v_{22})^2}{v_{12}^2+v_{22}^2}\bigg)\\
				&+4(1-\phi)\frac{v_{12}^2(v_{11}+v_{22})}{v_{12}^2+v_{22}^2}-\frac{2}{v_{11}v_{22}-v_{12}^2}\sum_k(v_{22}v_{1k}^2-v_{12}v_{1k}v_{2k})\\
				&-\frac{2\lambda v_{22}}{v_{11}v_{22}-v_{12}^2}-2\frac{v_1^2v_{22}-v_1v_2v_{12}}{v_{11}v_{22}-v_{12}^2}+3.
			\end{aligned}
		\end{equation}
		
		Finally, we compute the lower-order terms in \eqref{2-dim-Equality-F4}. Direct calculation shows that
		\begin{equation}\label{2-dim-Lower terms1}
			\begin{aligned}
				&~\bigg(1+(1-\phi)\frac{v_{11}v_{22}-v_{12}^2}{v_{12}^2+v_{22}^2}\bigg)^2+(1-\phi)^2\frac{v_{12}^2(v_{11}+v_{22})^2}{(v_{12}^2+v_{22}^2)^2}\\
				=&~1+2(1-\phi)\frac{v_{11}v_{22}-v_{12}^2}{v_{12}^2+v_{22}^2}+(1-\phi)^2\frac{v_{11}^2+v_{12}^2}{v_{12}^2+v_{22}^2},
			\end{aligned}
		\end{equation}
		\begin{equation}\label{2-dim-Lower terms2}
			\begin{aligned}
				&~\big(v_{12}^2+v_{22}^2-\phi(v_{11}v_{22}-v_{12}^2)\big)\bigg(1+(1-\phi)\frac{v_{11}v_{22}-v_{12}^2}{v_{12}^2+v_{22}^2}\bigg)\\
				&+(1-\phi)^2\frac{v_{12}^2(v_{11}+v_{22})^2}{v_{12}^2+v_{22}^2}\\
				=&~v_{12}^2+v_{22}^2+(1-2\phi)(v_{11}v_{22}-v_{12}^2)+(1-\phi)^2(v_{11}^2+v_{12}^2)\\
				&-(1-\phi)\frac{(v_{11}v_{22}-v_{12}^2)^2}{v_{12}^2+v_{22}^2},
			\end{aligned}
		\end{equation}
		\begin{equation}\label{2-dim-Lower terms3}
			\frac{2}{v_{11}v_{22}-v_{12}^2}\sum_k(v_{22}v_{1k}^2-v_{12}v_{1k}v_{2k})=2v_{11}
		\end{equation}
		and
		\begin{equation}\label{2-dim-Lower terms4}
			\begin{aligned}
				-\frac{2\lambda v_{22}}{v_{11}v_{22}-v_{12}^2}-2\frac{v_1^2v_{22}-v_1v_2v_{12}}{v_{11}v_{22}-v_{12}^2}
				=&-\frac{2v_{22}(v_{11}+v_{22})}{v_{11}v_{22}-v_{12}^2}+2\phi v_{11}\\
				=&~2\phi v_{11}-2-2\frac{v_{12}^2+v_{22}^2}{v_{11}v_{22}-v_{12}^2},
			\end{aligned}
		\end{equation}
		where we have used \eqref{Equation-v}, \eqref{2-dim-Coordinate condition'} and \eqref{2-dim-Condition-phi}. Hence, \eqref{2-dim-Lower terms1}-\eqref{2-dim-Lower terms4} reduce \eqref{2-dim-Equality-F4} to
		\begin{equation}\label{2-dim-Equality-F5}
			\begin{aligned}
				\frac{\Delta\phi}{2\phi}
				=&~\frac{1-\phi}{\phi}\bigg(3+2\frac{v_{12}^2+v_{22}^2}{v_{11}v_{22}-v_{12}^2}-2\phi v_{11}\bigg)\\
				&+\frac{(1-\phi)^2}{\phi}(v_{11}^2+v_{12}^2)\bigg(\frac{2}{v_{11}v_{22}-v_{12}^2}+\frac{1}{v_{12}^2+v_{22}^2}\bigg)\\
				&+4(1-\phi)\frac{v_{12}^2(v_{11}+v_{22})}{v_{12}^2+v_{22}^2}.
			\end{aligned}
		\end{equation}
		Since $\nabla^2v>0$ and $0<\phi<1$, the last two terms in \eqref{2-dim-Equality-F5} are nonnegative. For the first term, we apply \eqref{Equation-v}, \eqref{2-dim-Coordinate condition'} and \eqref{2-dim-Condition-phi} once more to obtain $$3+2\frac{v_{12}^2+v_{22}^2}{v_{11}v_{22}-v_{12}^2}-2\phi v_{11}=1+\frac{2}{v_{11}}\bigg(\lambda+\phi v_{12}^2+\frac{v_{12}^2(v_{11}+v_{22})}{v_{11}v_{22}-v_{12}^2}\bigg)>0.$$ Therefore, $\Delta\phi(x_0)>0$, contrary to $\Delta\phi(x_0)\leq 0$. This completes the proof.	
	\end{proof}
	
	\subsection{Proof of Theorem \ref{Theorem-2-dimensional}}\label{Subsection-Proof of Theorem 1.1}
	
	We now combine the maximum principle in Proposition \ref{Proposition-Maximum principle-2-dim} with the eigenvalue bound $\lambda>\frac{1}{4}$ to rule out Hessian degeneration along the domain deformation from Section \ref{Subsection-Domain deformation} and prove Theorem \ref{Theorem-2-dimensional}.
	
	\begin{proof}[Proof of Theorem \ref{Theorem-2-dimensional}]
		Consider the domain deformation $\{\Omega_t\}_{t\geq 0}$ from Section \ref{Subsection-Domain deformation}. Let $u_t$ be the first Dirichlet eigenfunction on $\Omega_t$ with eigenvalue $\lambda_t$ and set $v_t=-\log u_t$. Denote the curvature of the boundary $\partial\Omega_t$ by $\kappa_t$ and define $\phi_t$ by \eqref{Auxiliary function-2-dim} with $v$ replaced by $v_t$. According to Proposition \ref{Proposition-Domain deformation}, there exists $T>0$ such that $\nabla^2v_t>0$ and $\phi_t<1$ throughout $\Omega_t$ whenever $t\geq T$. Let $t_\ast\in[0,T)$ be the smallest time for which both conditions hold for all $t>t_\ast$. We first verify that both conditions remain valid at $t=t_\ast$ and then show that $t_\ast=0$.
		
		By continuity, $\nabla^2v_{t_\ast}\geq 0$ in $\Omega_{t_\ast}$. For the computations below we write $v=v_{t_\ast}$ and $\lambda=\lambda_{t_\ast}$. Combining $\phi_t<1$ with \eqref{Auxiliary function-2-dim} and letting $t\rightarrow t_\ast$ leads to
		\begin{equation}\label{2-dim-Inequality-phi-t_ast}
			(v_{11}v_1+v_{12}v_2)^2+(v_{12}v_1+v_{22}v_2)^2\leq(v_{11}v_{22}-v_{12}^2)^2.
		\end{equation}
		Suppose that $\nabla^2v$ degenerates at some $x_0\in\Omega_{t_\ast}$. Choose orthonormal coordinates at $x_0$ such that $\nabla^2v(x_0)$ is diagonal. Since $\nabla^2v\geq 0$ and \eqref{Equation-v} guarantees $\Delta v>0$, we may assume $v_{11}>0$ and $v_{12}=v_{22}=0$ at $x_0$. Then, \eqref{2-dim-Inequality-phi-t_ast} forces $v_1=0$. All computations below are performed at $x_0$.
		
		Since $v_{22}$ attains the local minimum $0$ at $x_0$, it follows that $v_{222}=0$. Differentiating equation \eqref{Equation-v} with respect to $x_2$, we obtain $$v_{112}=-v_{222}+2v_1v_{12}+2v_2v_{22}=0.$$ The commutation formula \eqref{Formula-Third derivative commutation} then gives $$v_{121}=v_{211}=v_{112}-v_2=-v_2.$$ On the other hand, \eqref{2-dim-Inequality-phi-t_ast} and the positive semi-definiteness of $\nabla^2v$ provide the estimate $$|v_{11}v_1+v_{12}v_2|\leq v_{11}v_{22}-v_{12}^2.$$ Both $v_{11}v_{22}-v_{12}^2$ and $v_{11}v_{22}-v_{12}^2+v_{11}v_1+v_{12}v_2$ attain the local minimum $0$ at $x_0$ and hence $$\nabla(v_{11}v_1+v_{12}v_2)(x_0)=0.$$ In particular, the relations $v_{12}=v_1=0$ and $v_{121}=-v_2$ simplify the derivative in the $x_1$ direction to $$0=v_2v_{121}+v_{11}^2=v_{11}^2-v_2^2.$$ Furthermore, equation \eqref{Equation-v} can now be written as $$\lambda=v_{11}-v_2^2=v_{11}(1-v_{11})\leq\frac{1}{4},$$ which contradicts Lemma \ref{Lemma-Eigenvalue lower bound}. Hence, $\nabla^2v_{t_\ast}>0$ throughout $\Omega_{t_\ast}$.
		
		Now that $\phi_{t_\ast}$ is well defined, letting $t\rightarrow t_\ast$ in the estimate of Proposition \ref{Proposition-Maximum principle-2-dim} for $t>t_\ast$ implies $$\phi_{t_\ast}\leq\max_{\partial\Omega_{t_\ast}}\kappa_{t_\ast}^{-2}<1.$$ It follows that both conditions still hold at $t=t_\ast$.
		
		We now show that $\nabla^2v_t>0$ and $\phi_t<1$ in $\Omega_t$ for all $t\geq 0$. Suppose for contradiction that $t_\ast>0$. By continuity, both conditions continue to hold on compact subsets of $\Omega_{t_\ast}$ for $t$ sufficiently close to $t_\ast$. Near the boundary, Proposition \ref{Proposition-Boundary convexity estimate} and \eqref{Boundary behavior-2-dim} ensure that they hold uniformly for $t$ close to $t_\ast$. Hence, there exists $\varepsilon\in(0,t_\ast)$ such that $\nabla^2v_t>0$ and $\phi_t<1$ throughout $\Omega_t$ for every $t\in[t_\ast-\varepsilon,t_\ast]$. This contradicts the definition of $t_\ast$. Thus, $t_\ast=0$ and $\nabla^2v>0$ in $\Omega$.
		
		It remains to prove the assertion about level sets. Since $\Omega$ is strictly horo-convex and $v$ is strictly convex with $v\rightarrow+\infty$ at the boundary, $v$ has exactly one critical point. This point is its unique minimum point and the unique maximum point of $u$. Therefore, each level set $\Sigma_t$ with $0<t<\max_\Omega u$ is a smooth closed curve and is strictly convex by \eqref{Equality-Second fundamental form}. Furthermore, \eqref{Estimate-Auxiliary function-2-dim} and \eqref{Estimate-Boundary-2-dim} imply $$\frac{1}{K^2}\leq\phi\leq\max_{\partial\Omega}\kappa^{-2}$$ and hence $$K\geq\min_{\partial\Omega}\kappa>1.$$ This completes the proof.
	\end{proof}
	
	At the end of this section, we prove Corollary \ref{Corollary-n=2}, where the domain is assumed to be horo-convex rather than strictly horo-convex.
	
	\begin{proof}[Proof of Corollary \ref{Corollary-n=2}]
		As in Wei and Xiao \cite[Section 4.1]{Wei-Xiao2025}, we first approximate $\Omega$ by smooth horo-convex domains from which their flow can start. Their Theorem 3.1 together with the dimension reduction in their Section 4.1 ensures strict horo-convexity at positive times. Choosing sufficiently small positive times, we obtain smooth strictly horo-convex domains $\Omega_\varepsilon$ with $\partial\Omega_\varepsilon\rightarrow\partial\Omega$ in $C^2$.
		
		Let $u_\varepsilon$ be the first Dirichlet eigenfunction on $\Omega_\varepsilon$ with the same normalization as $u$ and set $v_\varepsilon=-\log u_\varepsilon$. Standard elliptic theory ensures that $u_\varepsilon\rightarrow u$ and $v_\varepsilon\rightarrow v$ locally smoothly in $\Omega$. Theorem \ref{Theorem-2-dimensional} and the estimates in its proof therefore imply $$\nabla^2v\geq 0\quad\mathrm{in}\ \Omega$$ and $$(v_{11}v_1+v_{12}v_2)^2+(v_{12}v_1+v_{22}v_2)^2\leq(v_{11}v_{22}-v_{12}^2)^2.$$ Since $\lambda>\frac{1}{4}$ by Lemma \ref{Lemma-Eigenvalue lower bound}, the rank-one degeneration argument in Section \ref{Subsection-Proof of Theorem 1.1} applies. Hence, $$\nabla^2v>0\quad\mathrm{in}\ \Omega.$$
	\end{proof}
	
	\section{The three-dimensional case}\label{Section-n=3}
	
	We turn to the proof of Theorem \ref{Theorem-Higher-dimensional} in dimension three. Let $u$ be the first Dirichlet eigenfunction on a bounded, smooth and strictly horo-convex domain $\Omega\subset\mathbb{H}^3$ satisfying \eqref{Equation-Dirichlet problem} and set $v=-\log u$. The proof follows the same strategy as in dimension two. Throughout this section, the auxiliary function $\phi$ is defined in \eqref{Auxiliary function-3-dim} and $\kappa_1,\kappa_2$ denote the principal curvatures of $\partial\Omega$. The Latin indices $i,j,k,l,m$ range over $1,2,3$ and the Greek index $\alpha$ over $1,2$.
	
	\subsection{A maximum principle}
	
	We first bound $\phi$ by its boundary value.
	
	\begin{proposition}\label{Proposition-Maximum principle-3-dim}
		Suppose that $\nabla^2v>0$ and $\phi<1$ in $\Omega$. Then
		\begin{equation}\label{Estimate-Boundary-3-dim}
			\phi\leq\max_{\partial\Omega}(\kappa_1\kappa_2)^{-1}<1.
		\end{equation}
	\end{proposition}
	\begin{proof}
		According to \eqref{Boundary behavior-higher dim}, $\phi$ extends continuously to $\bar{\Omega}$ with positive boundary values $(\kappa_1\kappa_2)^{-1}$. If \eqref{Estimate-Boundary-3-dim} fails, $\phi$ reaches its maximum at some $x_0\in\Omega$. At $x_0$, $0<\phi(x_0)<1$, $\nabla\phi(x_0)=0$ and $\Delta\phi(x_0)\leq 0$. Since $\phi(x_0)>0$ and $\nabla^2v>0$, we have $\nabla|\nabla v|^2(x_0)\neq 0$. Choose normal orthonormal coordinates at $x_0$ so that the third coordinate vector points in this direction and the Hessian is diagonal on its orthogonal plane. In these coordinates, $v_{12}=0$ and the first two components of $\nabla|\nabla v|^2$ vanish. Thus, $v_1v_{11}+v_3v_{13}=0$ and $v_2v_{22}+v_3v_{23}=0$. Since $v_{11}>0$ and $v_{22}>0$, these identities imply
		\begin{equation}\label{3-dim-Coordinate condition}
			v_1=-\frac{v_{13}}{v_{11}}v_3\quad\mathrm{and}\ v_2=-\frac{v_{23}}{v_{22}}v_3.
		\end{equation}
		The positive definiteness of the Hessian also ensures that $$v_{33}-\frac{v_{13}^2}{v_{11}}-\frac{v_{23}^2}{v_{22}}=\frac{\det(v_{ij})}{v_{11}v_{22}}>0.$$ By the choice of frame and \eqref{3-dim-Coordinate condition}, $$0<v_1v_{13}+v_2v_{23}+v_3v_{33}=v_3\bigg(v_{33}-\frac{v_{13}^2}{v_{11}}-\frac{v_{23}^2}{v_{22}}\bigg),$$ so $v_3>0$ by the preceding inequality.
		
		All calculations below are performed at $x_0$. The relations \eqref{3-dim-Coordinate condition} show that
		\begin{equation}\label{3-dim-Identity}
			\sum_{i,j}v_iv_jv_{ij}=v_3^2\bigg(v_{33}-\frac{v_{13}^2}{v_{11}}-\frac{v_{23}^2}{v_{22}}\bigg),
		\end{equation} 
		and reduce $\phi$ to
		\begin{equation}\label{3-dim-Coordinate condition'}
			\phi=\frac{v_3^2}{v_{11}v_{22}}.
		\end{equation}
		Differentiating \eqref{Auxiliary function-3-dim} and applying $\nabla\phi(x_0)=0$, we have
		\begin{equation}\label{3-dim-1st condition}
			\begin{aligned}
				0=\frac{\phi_k}{\phi}
				=&~\frac{\sum_{i,j}v_iv_jv_{ijk}+2\sum_{i,j}v_iv_{ij}v_{jk}}{\sum_{i,j}v_iv_jv_{ij}}-\sum_{i,j}v^{ij}v_{ijk}\\
				=&~2\frac{v_{3k}}{v_3}-\frac{v_{11k}}{v_{11}}-\frac{v_{22k}}{v_{22}}.
			\end{aligned}
		\end{equation}
		Here the derivative of the determinant is computed from the identity $$\frac{(\det(v_{ij}))_k}{\det(v_{ij})}=\sum_{i,j}v^{ij}v_{ijk}.$$ Differentiating \eqref{Auxiliary function-3-dim} once more, summing over $k$ and using $\nabla\phi=0$, we arrive at
		\begin{equation}\label{3-dim-Equality-F1}
			\begin{aligned}
				\frac{\Delta\phi}{\phi}
				=&~\frac{1}{\sum_{i,j}v_iv_jv_{ij}}\sum_{i,j,k}\big(v_iv_jv_{ijkk}+4v_iv_{jk}v_{ijk}+2v_{ij}v_{ik}v_{jk}+2v_iv_{ij}v_{jkk}\big)\\
				&-\sum_{i,j,k}v^{ij}v_{ijkk}+\sum_{i,j,l,m,k}v^{ij}v^{lm}\big(v_{ilk}v_{jmk}-v_{ijk}v_{lmk}\big).
			\end{aligned}
		\end{equation}
		
		We first handle the fourth-order derivative terms in \eqref{3-dim-Equality-F1}. Taking the derivative of equation \eqref{Equation-v} and invoking \eqref{Formula-Third derivative commutation} and \eqref{Formula-Ricci commutation}, we obtain
		\begin{equation}\label{3-dim-3rd derivative term}
			\sum_kv_{ikk}=(\Delta v)_i-2v_i=2\sum_kv_kv_{ki}-2v_i,
		\end{equation}
		and
		\begin{align*}
			(\Delta v)_{ij}
			=&~2\sum_kv_{ki}v_{kj}+2\sum_kv_kv_{kij}\\
			=&~2\sum_kv_{ki}v_{kj}+2\sum_kv_kv_{ijk}+2v_iv_j-2\delta_{ij}\sum_kv_k^2.
		\end{align*}
		Furthermore, \eqref{Formula-Fourth derivative commutation} implies
		\begin{equation}\label{3-dim-4th derivative term}
			\begin{aligned}
				\sum_kv_{ijkk}
				=&~(\Delta v)_{ij}-6v_{ij}+2\delta_{ij}\Delta v\\
				=&~2\sum_kv_kv_{ijk}+2\sum_kv_{ki}v_{kj}+2v_iv_j-6v_{ij}+2\lambda\delta_{ij}.
			\end{aligned}
		\end{equation}
		To simplify the coefficients of the fourth-order derivatives, we use \eqref{3-dim-Identity} and the explicit expressions for the inverse Hessian to write
		\begin{equation}\label{3-dim-Inverse Hessian identity}
			\frac{v_iv_j}{\sum_{l,m}v_lv_mv_{lm}}-v^{ij}=-\frac{\delta_{i1}\delta_{j1}}{v_{11}}-\frac{\delta_{i2}\delta_{j2}}{v_{22}}.
		\end{equation}
		Together with \eqref{3-dim-4th derivative term}, this identity reduces the fourth-order derivative terms in \eqref{3-dim-Equality-F1} to
		\begin{align*}
			\sum_{i,j}\bigg(\frac{v_iv_j}{\sum_{l,m}v_lv_mv_{lm}}-v^{ij}\bigg)\sum_kv_{ijkk}
			=&-\frac{1}{v_{11}}\sum_kv_{11kk}-\frac{1}{v_{22}}\sum_kv_{22kk}\\
			=&-2\sum_kv_k\bigg(\frac{v_{11k}}{v_{11}}+\frac{v_{22k}}{v_{22}}\bigg)-2\bigg(v_{11}+v_{22}+\frac{v_{13}^2}{v_{11}}+\frac{v_{23}^2}{v_{22}}\bigg)\\
			&-2\bigg(\frac{v_1^2}{v_{11}}+\frac{v_2^2}{v_{22}}\bigg)+12-2\lambda\bigg(\frac{1}{v_{11}}+\frac{1}{v_{22}}\bigg),
		\end{align*}
		where $v_{12}=0$ has been used in the last equality. Hence, the preceding identity and \eqref{3-dim-3rd derivative term} allow us to rewrite \eqref{3-dim-Equality-F1} as
		\begin{equation}\label{3-dim-Equality-F2}
			\begin{aligned}
				\frac{\Delta\phi}{\phi}
				=&-2\sum_kv_k\bigg(\frac{v_{11k}}{v_{11}}+\frac{v_{22k}}{v_{22}}\bigg)-2\bigg(v_{11}+v_{22}+\frac{v_{13}^2}{v_{11}}+\frac{v_{23}^2}{v_{22}}\bigg)\\
				&-2\bigg(\frac{v_1^2}{v_{11}}+\frac{v_2^2}{v_{22}}\bigg)+8-2\lambda\bigg(\frac{1}{v_{11}}+\frac{1}{v_{22}}\bigg)+\frac{4}{\sum_{i,j}v_iv_jv_{ij}}\sum_k\bigg(\sum_iv_iv_{ik}\bigg)^2\\
				&+\frac{1}{\sum_{i,j}v_iv_jv_{ij}}\bigg(4\sum_{i,j,k}v_iv_{jk}v_{ijk}+2\sum_{i,j,k}v_{ij}v_{ik}v_{jk}\bigg)\\
				&+\sum_{i,j,l,m,k}v^{ij}v^{lm}\big(v_{ilk}v_{jmk}-v_{ijk}v_{lmk}\big).
			\end{aligned}
		\end{equation}
		
		We next deal with the third-order derivative terms in \eqref{3-dim-Equality-F2} and combine them with the lower-order terms. For the first term and the term containing $(\sum_iv_iv_{ik})^2$, \eqref{3-dim-Coordinate condition} and \eqref{3-dim-1st condition} imply
		\begin{align*}
			&-2\sum_kv_k\bigg(\frac{v_{11k}}{v_{11}}+\frac{v_{22k}}{v_{22}}\bigg)+\frac{4}{\sum_{i,j}v_iv_jv_{ij}}\sum_k\bigg(\sum_iv_iv_{ik}\bigg)^2\\
			=&-\frac{4}{v_3}\sum_kv_kv_{3k}+4\bigg(v_{33}-\frac{v_{13}^2}{v_{11}}-\frac{v_{23}^2}{v_{22}}\bigg)=0.
		\end{align*}
		By \eqref{3-dim-Inverse Hessian identity}, the last term expands as
		\begin{align*}
			&\sum_{i,j,l,m,k}v^{ij}v^{lm}(v_{ilk}v_{jmk}-v_{ijk}v_{lmk})\\
			=&\sum_{\alpha,\beta=1}^2\sum_k\frac{v_{\alpha\beta k}^2}{v_{\alpha\alpha}v_{\beta\beta}}+\frac{2}{\sum_{i,j}v_iv_jv_{ij}}\sum_{\alpha=1}^2\sum_k\frac{1}{v_{\alpha\alpha}}\sum_{i,j}v_iv_jv_{\alpha ik}v_{\alpha jk}\\
			&+\frac{1}{\big(\sum_{i,j}v_iv_jv_{ij}\big)^2}\sum_k\sum_{i,j,l,m}v_iv_jv_lv_mv_{ilk}v_{jmk}-\sum_k\bigg(\frac{v_{11k}}{v_{11}}+\frac{v_{22k}}{v_{22}}+\frac{\sum_{i,j}v_iv_jv_{ijk}}{\sum_{i,j}v_iv_jv_{ij}}\bigg)^2\\
			=&\sum_k\bigg(\frac{v_{11k}^2}{v_{11}^2}+\frac{v_{22k}^2}{v_{22}^2}+\frac{2v_{12k}^2}{v_{11}v_{22}}\bigg)+\frac{2}{\sum_{i,j}v_iv_jv_{ij}}\sum_{\alpha=1}^2\sum_k\frac{1}{v_{\alpha\alpha}}\bigg(\sum_iv_iv_{\alpha ik}\bigg)^2\\
			&+\frac{1}{\big(\sum_{i,j}v_iv_jv_{ij}\big)^2}\sum_k\bigg(\sum_{i,j}v_iv_jv_{ijk}\bigg)^2-\sum_k\bigg(\frac{v_{11k}}{v_{11}}+\frac{v_{22k}}{v_{22}}\bigg)^2\\
			&-\frac{2}{\sum_{i,j}v_iv_jv_{ij}}\sum_k\bigg(\frac{v_{11k}}{v_{11}}+\frac{v_{22k}}{v_{22}}\bigg)\sum_{i,j}v_iv_jv_{ijk}-\frac{1}{\big(\sum_{i,j}v_iv_jv_{ij}\big)^2}\sum_k\bigg(\sum_{i,j}v_iv_jv_{ijk}\bigg)^2\\
			=&\sum_k\bigg(\frac{v_{11k}^2}{v_{11}^2}+\frac{v_{22k}^2}{v_{22}^2}+\frac{2v_{12k}^2}{v_{11}v_{22}}-\bigg(\frac{v_{11k}}{v_{11}}+\frac{v_{22k}}{v_{22}}\bigg)^2\bigg)+\frac{2}{\sum_{i,j}v_iv_jv_{ij}}\sum_{\alpha=1}^2\sum_k\frac{1}{v_{\alpha\alpha}}\bigg(\sum_iv_iv_{\alpha ik}\bigg)^2\\
			&-\frac{2}{\sum_{i,j}v_iv_jv_{ij}}\sum_k\bigg(\frac{v_{11k}}{v_{11}}+\frac{v_{22k}}{v_{22}}\bigg)\sum_{i,j}v_iv_jv_{ijk}.
		\end{align*}
		We combine the final term above with the term containing $\sum_{i,j,k}v_iv_{jk}v_{ijk}$ in \eqref{3-dim-Equality-F2} to obtain
		\begin{align*}
			&~\frac{2}{\sum_{i,j}v_iv_jv_{ij}}\sum_k\bigg(2\sum_{i,j}v_iv_{jk}v_{ijk}-\bigg(\frac{v_{11k}}{v_{11}}+\frac{v_{22k}}{v_{22}}\bigg)\sum_{i,j}v_iv_jv_{ijk}\bigg)\\
			=&\frac{4}{\sum_{i,j}v_iv_jv_{ij}}\sum_{\alpha=1}^2\sum_k\frac{1}{v_{\alpha\alpha}}\bigg(\sum_iv_iv_{\alpha ik}\bigg)\bigg(\sum_iv_{\alpha i}v_{ik}\bigg),
		\end{align*} 
		where we have used \eqref{3-dim-Coordinate condition} and \eqref{3-dim-1st condition}. The term involving $\sum_{i,j,k}v_{ij}v_{ik}v_{jk}$ in \eqref{3-dim-Equality-F2} can be rewritten with \eqref{3-dim-Inverse Hessian identity} and \eqref{3-dim-1st condition} as
		\begin{align*}
			\frac{2}{\sum_{i,j}v_iv_jv_{ij}}\sum_{i,j,k}v_{ij}v_{ik}v_{jk}
			=&~\frac{2}{\sum_{i,j}v_iv_jv_{ij}}\sum_{\alpha=1}^2\sum_k\frac{1}{v_{\alpha\alpha}}\bigg(\sum_iv_{\alpha i}v_{ik}\bigg)^2+\frac{2}{v_3^2}\sum_kv_{3k}^2\\
			=&~\frac{2}{\sum_{i,j}v_iv_jv_{ij}}\sum_{\alpha=1}^2\sum_k\frac{1}{v_{\alpha\alpha}}\bigg(\sum_iv_{\alpha i}v_{ik}\bigg)^2+\frac{1}{2}\sum_k\bigg(\frac{v_{11k}}{v_{11}}+\frac{v_{22k}}{v_{22}}\bigg)^2.
		\end{align*}
		Consequently, \eqref{3-dim-Equality-F2} becomes
		\begin{equation}\label{3-dim-Equality-F3}
			\begin{aligned}
				\frac{\Delta\phi}{\phi}
				=&~8-2\bigg(v_{11}+v_{22}+\frac{v_{13}^2}{v_{11}}+\frac{v_{23}^2}{v_{22}}\bigg)-2\bigg(\frac{v_1^2}{v_{11}}+\frac{v_2^2}{v_{22}}\bigg)-2\lambda\bigg(\frac{1}{v_{11}}+\frac{1}{v_{22}}\bigg)\\
				&+\sum_k\bigg(\frac{v_{11k}^2}{v_{11}^2}+\frac{v_{22k}^2}{v_{22}^2}+\frac{2v_{12k}^2}{v_{11}v_{22}}-\frac{1}{2}\bigg(\frac{v_{11k}}{v_{11}}+\frac{v_{22k}}{v_{22}}\bigg)^2\bigg)\\
				&+\frac{2}{\sum_{i,j}v_iv_jv_{ij}}\sum_{\alpha=1}^2\sum_k\frac{1}{v_{\alpha\alpha}}\\
				&\hspace{0.2cm}\bigg(\bigg(\sum_iv_iv_{\alpha ik}\bigg)^2+2\bigg(\sum_iv_iv_{\alpha ik}\bigg)\bigg(\sum_iv_{\alpha i}v_{ik}\bigg)+\bigg(\sum_iv_{\alpha i}v_{ik}\bigg)^2\bigg)\\
				=&~8-2\bigg(v_{11}+v_{22}+\frac{v_{13}^2}{v_{11}}+\frac{v_{23}^2}{v_{22}}\bigg)-2\bigg(\frac{v_1^2}{v_{11}}+\frac{v_2^2}{v_{22}}\bigg)-2\lambda\bigg(\frac{1}{v_{11}}+\frac{1}{v_{22}}\bigg)\\
				&+\frac{1}{2}\sum_k\bigg(\frac{v_{11k}}{v_{11}}-\frac{v_{22k}}{v_{22}}\bigg)^2+\frac{2}{v_{11}v_{22}}\sum_kv_{12k}^2\\
				&+\frac{2}{\sum_{i,j}v_iv_jv_{ij}}\sum_{\alpha=1}^2\sum_k\frac{1}{v_{\alpha\alpha}}\bigg(\sum_iv_iv_{\alpha ik}+\sum_iv_{\alpha i}v_{ik}\bigg)^2.
			\end{aligned}
		\end{equation}
		
		We continue to estimate the third-order derivative terms in \eqref{3-dim-Equality-F3}. Recalling \eqref{3-dim-1st condition}, we have
		\begin{equation}\label{3-dim-3rd derivative relation}
			v_{22k}=v_{22}\bigg(\frac{2v_{3k}}{v_3}-\frac{v_{11k}}{v_{11}}\bigg)
		\end{equation}
		and hence
		\begin{equation}\label{3-dim-3rd derivative identity1}
			\frac{1}{2}\sum_k\bigg(\frac{v_{11k}}{v_{11}}-\frac{v_{22k}}{v_{22}}\bigg)^2=\frac{2}{v_3^2}\sum_k\bigg(\frac{v_3v_{11k}}{v_{11}}-v_{3k}\bigg)^2.
		\end{equation}
		From the commutation formulas \eqref{Formula-Ricci commutation} and \eqref{Formula-Third derivative commutation} together with \eqref{3-dim-Coordinate condition} and \eqref{3-dim-1st condition}, we obtain $$v_{121}=v_{112}+\frac{v_3v_{23}}{v_{22}}$$ and $$v_{122}=v_{22}\bigg(\frac{2v_{13}}{v_3}-\frac{v_{111}}{v_{11}}\bigg)+\frac{v_3v_{13}}{v_{11}}.$$ It follows from these identities and \eqref{3-dim-Coordinate condition'} that
		\begin{equation}\label{3-dim-3rd derivative identity2}
			\begin{aligned}
				\frac{2}{v_{11}v_{22}}\sum_kv_{12k}^2
				=&~\frac{2}{v_3^2}\phi(v_{121}^2+v_{122}^2+v_{123}^2)\\
				=&~\frac{2}{v_3^2}\bigg(\frac{v_{11}}{v_{22}}\bigg(\frac{v_3v_{112}}{v_{11}}+\phi v_{23}\bigg)^2+\frac{v_{22}}{v_{11}}\bigg(\frac{v_3v_{111}}{v_{11}}-(2+\phi)v_{13}\bigg)^2+\phi v_{123}^2\bigg).
			\end{aligned}
		\end{equation}
		Since $\nabla^2v>0$, we can keep only the terms with $k=\alpha$ in the last term of \eqref{3-dim-Equality-F3}. Using \eqref{3-dim-Identity} and applying Young's inequality, we conclude that
		\begin{equation}\label{3-dm-Estimate-3rd derivative}
			\begin{aligned}
				&~\frac{2}{\sum_{i,j}v_iv_jv_{ij}}\sum_{\alpha=1}^2\sum_k\frac{1}{v_{\alpha\alpha}}\bigg(\sum_iv_iv_{\alpha ik}+\sum_iv_{\alpha i}v_{ik}\bigg)^2\\
				\geq&~\frac{2}{v_3^2\big(v_{33}-\frac{v_{13}^2}{v_{11}}-\frac{v_{23}^2}{v_{22}}\big)}\sum_{\alpha=1}^2\frac{1}{v_{\alpha\alpha}}\bigg(\sum_iv_iv_{\alpha i\alpha}+\sum_iv_{\alpha i}^2\bigg)^2\\
				\geq&~\frac{4}{v_3^2}\sum_{\alpha=1}^2\bigg(\sum_iv_iv_{\alpha i\alpha}+\sum_iv_{\alpha i}^2\bigg)-\frac{2(v_{11}+v_{22})}{v_3^2}\bigg(v_{33}-\frac{v_{13}^2}{v_{11}}-\frac{v_{23}^2}{v_{22}}\bigg).
			\end{aligned}
		\end{equation}
		To calculate the first term on the right-hand side of \eqref{3-dm-Estimate-3rd derivative}, we write the third-order terms using $v_{111}$, $v_{112}$ and $v_{113}$. Applying \eqref{3-dim-3rd derivative relation} with $k=2,3$ and \eqref{Formula-Third derivative commutation}, we obtain the following identities
		\begin{align*}
			v_{222}=&~v_{22}\bigg(\frac{2v_{23}}{v_3}-\frac{v_{112}}{v_{11}}\bigg),\\
			v_{232}=&~v_{22}\bigg(\frac{2v_{33}}{v_3}-\frac{v_{113}}{v_{11}}\bigg)-v_3,
			\shortintertext{and}
			v_{131}=&~v_{113}-v_3.
		\end{align*}
		These identities combined with the expressions for $v_{121}$ and $v_{122}$ lead to
		\begin{equation}\label{3-dim-3rd derivative identity3}
			\begin{aligned}
				&\sum_{\alpha=1}^2\bigg(\sum_iv_iv_{\alpha i\alpha}+\sum_iv_{\alpha i}^2\bigg)\\
				=&~v_1(v_{111}+v_{122})+v_2(v_{121}+v_{222})+v_3(v_{131}+v_{232})+v_{11}^2+v_{22}^2+v_{13}^2+v_{23}^2\\
				=&-\frac{v_3v_{13}}{v_{11}}\bigg(v_{111}+v_{22}\bigg(\frac{2v_{13}}{v_3}-\frac{v_{111}}{v_{11}}\bigg)+\frac{v_3v_{13}}{v_{11}}\bigg)\\
				&-\frac{v_3v_{23}}{v_{22}}\bigg(v_{112}+\frac{v_3v_{23}}{v_{22}}+v_{22}\bigg(\frac{2v_{23}}{v_3}-\frac{v_{112}}{v_{11}}\bigg)\bigg)\\
				&+v_3\bigg(v_{113}+v_{22}\bigg(\frac{2v_{33}}{v_3}-\frac{v_{113}}{v_{11}}\bigg)-2v_3\bigg)+v_{11}^2+v_{22}^2+v_{13}^2+v_{23}^2\\
				=&-\frac{v_3v_{13}(v_{11}-v_{22})}{v_{11}^2}v_{111}-\frac{v_3v_{23}(v_{11}-v_{22})}{v_{11}v_{22}}v_{112}+\frac{v_3(v_{11}-v_{22})}{v_{11}}v_{113}\\
				&+2v_{22}v_{33}+v_{11}^2+v_{22}^2-2v_3^2+\bigg(1-\frac{v_{22}}{v_{11}}(2+\phi)\bigg)v_{13}^2-\bigg(1+\frac{v_{11}}{v_{22}}\phi\bigg)v_{23}^2,
			\end{aligned}
		\end{equation}
		where \eqref{3-dim-Coordinate condition} has been used in the second equality and \eqref{3-dim-Coordinate condition'} in the last equality.
		
		Combining \eqref{3-dim-3rd derivative identity1}-\eqref{3-dim-3rd derivative identity3} and collecting terms in $\frac{v_3v_{111}}{v_{11}}$, $\frac{v_3v_{112}}{v_{11}}$ and $\frac{v_3v_{113}}{v_{11}}$, we obtain the following estimate from \eqref{3-dim-Equality-F3}
		\begin{equation}\label{3-dim-Inequality-F1}
			\begin{aligned}
				\frac{\Delta\phi}{\phi}
				\geq&~8-2\bigg(v_{11}+v_{22}+\frac{v_{13}^2}{v_{11}}+\frac{v_{23}^2}{v_{22}}\bigg)-2\bigg(\frac{v_1^2}{v_{11}}+\frac{v_2^2}{v_{22}}\bigg)-2\lambda\bigg(\frac{1}{v_{11}}+\frac{1}{v_{22}}\bigg)\\
				&+\frac{2}{v_3^2}\bigg(\frac{v_{11}+v_{22}}{v_{11}}\bigg(\frac{v_3v_{111}}{v_{11}}\bigg)^2+\frac{v_{11}+v_{22}}{v_{22}}\bigg(\frac{v_3v_{112}}{v_{11}}\bigg)^2+\bigg(\frac{v_3v_{113}}{v_{11}}\bigg)^2\\
				&\hspace{1.2cm}-\frac{2(2v_{11}+(1+\phi)v_{22})v_{13}}{v_{11}}\frac{v_3v_{111}}{v_{11}}-\frac{2v_{11}(1-\phi)v_{23}}{v_{22}}\frac{v_3v_{112}}{v_{11}}\\
				&\hspace{1.2cm}+2(v_{11}-v_{22}-v_{33})\frac{v_3v_{113}}{v_{11}}+\phi v_{123}^2+2v_{11}^2+2v_{22}^2+v_{33}^2-4v_3^2\\
				&\hspace{1.2cm}+(3v_{22}-v_{11})v_{33}+\bigg(4+\frac{v_{22}}{v_{11}}(1+\phi)^2\bigg)v_{13}^2+\frac{v_{11}}{v_{22}}(1-\phi)^2v_{23}^2\bigg)\\
				=&~8-2\bigg(v_{11}+v_{22}+\frac{v_{13}^2}{v_{11}}+\frac{v_{23}^2}{v_{22}}\bigg)-2\bigg(\frac{v_1^2}{v_{11}}+\frac{v_2^2}{v_{22}}\bigg)-2\lambda\bigg(\frac{1}{v_{11}}+\frac{1}{v_{22}}\bigg)\\
				&+\frac{2}{v_3^2}\bigg(\frac{v_{11}+v_{22}}{v_{11}}\bigg(\frac{v_3v_{111}}{v_{11}}-\frac{(2v_{11}+(1+\phi)v_{22})v_{13}}{v_{11}+v_{22}}\bigg)^2\\
				&\hspace{1.2cm}+\frac{v_{11}+v_{22}}{v_{22}}\bigg(\frac{v_3v_{112}}{v_{11}}-\frac{v_{11}(1-\phi)v_{23}}{v_{11}+v_{22}}\bigg)^2\\
				&\hspace{1.2cm}+\frac{v_{11}v_{22}(1-\phi)^2}{v_{11}+v_{22}}\bigg(\frac{v_{13}^2}{v_{11}}+\frac{v_{23}^2}{v_{22}}\bigg)+\bigg(\frac{v_3v_{113}}{v_{11}}+v_{11}-v_{22}-v_{33}\bigg)^2\\
				&\hspace{1.2cm}+\phi v_{123}^2+(v_{11}+v_{22})v_{33}+(v_{11}+v_{22})^2-4v_3^2\bigg)\\
				\geq&~\frac{2(v_{11}+v_{22})}{v_3^2}\Delta v-2\bigg(v_{11}+v_{22}+\frac{v_{13}^2}{v_{11}}+\frac{v_{23}^2}{v_{22}}\bigg)-2\bigg(\frac{v_1^2}{v_{11}}+\frac{v_2^2}{v_{22}}\bigg)\\
				&-2\lambda\bigg(\frac{1}{v_{11}}+\frac{1}{v_{22}}\bigg),
			\end{aligned}
		\end{equation}
		where the last inequality follows from $v_{11},v_{22}>0$ and $\phi>0$.
		
		Finally, equation \eqref{Equation-v} and \eqref{3-dim-Coordinate condition} yield $$\Delta v=\lambda+v_3^2+\frac{v_3^2v_{13}^2}{v_{11}^2}+\frac{v_3^2v_{23}^2}{v_{22}^2}.$$ Substituting this identity into \eqref{3-dim-Inequality-F1} and using \eqref{3-dim-Coordinate condition} and \eqref{3-dim-Coordinate condition'}, we conclude that
		\begin{align*}
			\frac{\Delta\phi}{\phi}
			\geq&~\frac{2\lambda}{v_3^2}(v_{11}+v_{22})-2\lambda\bigg(\frac{1}{v_{11}}+\frac{1}{v_{22}}\bigg)\\
			&+2\frac{v_{22}v_{13}^2}{v_{11}^2}+2\frac{v_{11}v_{23}^2}{v_{22}^2}-2\frac{v_3^2v_{13}^2}{v_{11}^3}-2\frac{v_3^2v_{23}^2}{v_{22}^3}\\
			=&~\frac{2\lambda}{v_3^2}(v_{11}+v_{22})(1-\phi)+2\frac{v_{13}^2}{v_{11}^2}\bigg(v_{22}-\frac{v_3^2}{v_{11}}\bigg)+2\frac{v_{23}^2}{v_{22}^2}\bigg(v_{11}-\frac{v_3^2}{v_{22}}\bigg)\\
			=&~2(1-\phi)\bigg(\frac{\lambda(v_{11}+v_{22})}{v_3^2}+\frac{v_{22}v_{13}^2}{v_{11}^2}+\frac{v_{11}v_{23}^2}{v_{22}^2}\bigg)\\
			>&~0,
		\end{align*}
		since $\nabla^2v>0$, $0<\phi<1$ and $\lambda>0$. This contradicts $\Delta\phi(x_0)\leq 0$ and completes the proof. 
	\end{proof}
	
	\subsection{Proof of Theorem \ref{Theorem-Higher-dimensional} for $n=3$}\label{Subsection-Proof of Theorem 1.2-n=3}
	
	With Proposition \ref{Proposition-Maximum principle-3-dim} established, we turn to the proof of Theorem \ref{Theorem-Higher-dimensional} for $n=3$. We use the domain deformation as in Section \ref{Subsection-Proof of Theorem 1.1} and show that neither of the two possible types of Hessian degeneration can occur.
	\begin{proof}
		Let $\{\Omega_t\}_{t\geq 0}$ and $v_t$ be as in Section \ref{Subsection-Domain deformation}. Define $\phi_t$ by \eqref{Auxiliary function-3-dim} with $v$ replaced by $v_t$. Denote the first Dirichlet eigenvalue of $\Omega_t$ by $\lambda_t$ and its boundary principal curvatures by $\kappa_{1,t}$, $\kappa_{2,t}$. Choose $T>0$ as in Proposition \ref{Proposition-Domain deformation} so that $\nabla^2v_t>0$ and $\phi_t<1$ in $\Omega_t$ whenever $t\geq T$. We take $t_\ast\in[0,T)$ to be the smallest time such that both conditions hold for every $t>t_\ast$. As in the proof of Theorem \ref{Theorem-2-dimensional}, we complete the argument by verifying both conditions at $t=t_\ast$ and showing that $t_\ast=0$.
		
		For the following calculations, write $v=v_{t_\ast}$ and $\lambda=\lambda_{t_\ast}$. By continuity, $\nabla^2v\geq 0$. Taking $t\rightarrow t_\ast$ in \eqref{Estimate-Boundary-3-dim}, we obtain
		\begin{equation}\label{3-dim-Inequality-phi-t_ast}
			\sum_{i,j}v_iv_jv_{ij}\leq\bigg(\max_{\partial\Omega_{t_\ast}}(\kappa_{1,t_\ast}\kappa_{2,t_\ast})^{-1}\bigg)\det(v_{ij}).
		\end{equation}
		Assume that $\nabla^2v$ is degenerate at some $x_0\in\Omega_{t_\ast}$. Since $\nabla^2v\geq 0$ and equation \eqref{Equation-v} gives $\Delta v>0$, its rank at $x_0$ is either two or one. Moreover, \eqref{3-dim-Inequality-phi-t_ast} shows that both $\det(v_{ij})$ and $$\bigg(\max_{\partial\Omega_{t_\ast}}(\kappa_{1,t_\ast}\kappa_{2,t_\ast})^{-1}\bigg)\det(v_{ij})-\sum_{i,j}v_iv_jv_{ij}$$ attain their minimum value $0$ at $x_0$. Their Laplacians are therefore nonnegative there.
		
		Choose orthonormal coordinates at $x_0$ such that $\nabla^2v(x_0)$ is diagonal. All calculations below are performed at $x_0$.
		
		\textbf{Case 1.} $\nabla^2v$ has rank two. We may assume that $v_{11}>0$, $v_{22}>0$ and $v_{33}=0$. Then, \eqref{3-dim-Inequality-phi-t_ast} forces $v_1=v_2=0$. To calculate $\Delta\det(v_{ij})$, we note that $v_{33}$ attains the local minimum $0$ at $x_0$. Thus, $v_{33k}=0$. Differentiating \eqref{Equation-v} with respect to $x_3$ yields $$v_{113}+v_{223}=-v_{333}+2\sum_iv_iv_{i3}=0.$$ The commutation formulas \eqref{Formula-Ricci commutation} and \eqref{Formula-Third derivative commutation} therefore imply
		\begin{align*}
			v_{131}=&~v_{113}-v_3,\\
			v_{232}=&~-v_{113}-v_3,\\
			v_{133}=&~v_{233}=0,\\
			\intertext{and}
			v_{132}=&~v_{231}=v_{123}.
		\end{align*}
		Moreover, \eqref{3-dim-4th derivative term} reduces to $$\sum_kv_{33kk}=2\lambda+2v_3^2=2(v_{11}+v_{22}).$$ It follows that
		\begin{equation}\label{3-dim-Laplacian-det-rank two}
			\begin{aligned}
				\Delta\det(v_{ij})
				=&~v_{11}v_{22}\sum_kv_{33kk}-2v_{22}\sum_kv_{13k}^2-2v_{11}\sum_kv_{23k}^2\\
				=&~2v_{11}v_{22}(v_{11}+v_{22})-2v_{22}\big((v_{113}-v_3)^2+v_{123}^2\big)\\
				&-2v_{11}\big((v_{113}+v_3)^2+v_{123}^2\big).
			\end{aligned}
		\end{equation}
		Since $\sum_iv_iv_{ij}=0$ for every $j$, we also calculate
		\begin{align*}
			\Delta\bigg(\sum_{i,j}v_iv_jv_{ij}\bigg)
			=&~\sum_{i,j}v_iv_j\sum_kv_{ijkk}+4\sum_{i,j,k}v_iv_{jk}v_{ijk}+2\sum_{i,j,k}v_{ik}v_{jk}v_{ij}\\
			&+4\sum_j\bigg(\sum_iv_iv_{ij}\bigg)^2-4\sum_{i,j}v_iv_jv_{ij}\\
			=&~v_3^2\sum_kv_{33kk}+4v_3(v_{11}v_{131}+v_{22}v_{232})+2(v_{11}^3+v_{22}^3)\\
			=&-2v_3^2(v_{11}+v_{22})+4v_3(v_{11}-v_{22})v_{113}+2(v_{11}^3+v_{22}^3).
		\end{align*}
		Thus,
		\begin{equation}\label{3-dim-Degeneration rank two}
			\begin{aligned}
				&~\Delta\bigg(\det(v_{ij})-\sum_{i,j}v_iv_jv_{ij}\bigg)\\
				=&-2(v_{11}+v_{22})(v_{113}^2+v_{123}^2)-8v_3(v_{11}-v_{22})v_{113}-2(v_{11}+v_{22})(v_{11}-v_{22})^2\\
				=&-2(v_{11}+v_{22})\bigg(\bigg(v_{113}+\frac{2v_3(v_{11}-v_{22})}{v_{11}+v_{22}}\bigg)^2+v_{123}^2\bigg)\\
				&+\frac{2(v_{11}-v_{22})^2}{v_{11}+v_{22}}\big(4v_3^2-(v_{11}+v_{22})^2\big).
			\end{aligned}
		\end{equation}
		Recalling that $\lambda>1$ by Lemma \ref{Lemma-Eigenvalue lower bound}, we rewrite \eqref{Equation-v} as $$\lambda=v_{11}+v_{22}-v_3^2$$ and arrive at $$4v_3^2-(v_{11}+v_{22})^2=-(v_{11}+v_{22}-2)^2-4(\lambda-1)<0.$$ Hence, the right-hand side of \eqref{3-dim-Degeneration rank two} is nonpositive and the two minimum conditions yield
		\begin{align*}
			0
			\leq&~\Delta\bigg(\bigg(\max_{\partial\Omega_{t_\ast}}(\kappa_{1,t_\ast}\kappa_{2,t_\ast})^{-1}\bigg)\det(v_{ij})-\sum_{i,j}v_iv_jv_{ij}\bigg)\\
			=&~\Delta\bigg(\det(v_{ij})-\sum_{i,j}v_iv_jv_{ij}\bigg)-\bigg(1-\max_{\partial\Omega_{t_\ast}}(\kappa_{1,t_\ast}\kappa_{2,t_\ast})^{-1}\bigg)\Delta\det(v_{ij})\\
			\leq&-\bigg(1-\max_{\partial\Omega_{t_\ast}}(\kappa_{1,t_\ast}\kappa_{2,t_\ast})^{-1}\bigg)\Delta\det(v_{ij})\leq 0.
		\end{align*}
		By strict horo-convexity, $1-\max_{\partial\Omega_{t_\ast}}(\kappa_{1,t_\ast}\kappa_{2,t_\ast})^{-1}>0$. Hence, $$\Delta\det(v_{ij})=0\quad\mathrm{and}\ \Delta\bigg(\det(v_{ij})-\sum_{i,j}v_iv_jv_{ij}\bigg)=0.$$ It follows from \eqref{3-dim-Degeneration rank two} that $v_{11}=v_{22}$ and $v_{113}=v_{123}=0$. Relation \eqref{3-dim-Laplacian-det-rank two} then reduces to $0=4v_{11}(v_{11}^2-v_3^2)$. Since $v_{11}>0$, this gives $v_3^2=v_{11}^2$. Therefore, $$\lambda=2v_{11}-v_{11}^2\leq 1,$$ in contradiction with Lemma \ref{Lemma-Eigenvalue lower bound}.
		
		\textbf{Case 2.} $\nabla^2v$ has rank one. Take $v_{11}>0$ and $v_{22}=v_{33}=0$. Then, $v_1=0$ by \eqref{3-dim-Inequality-phi-t_ast}. As in Case 1, we first calculate $\Delta\det(v_{ij})$. Here the nonnegative functions $v_{22}$, $v_{33}$, $v_{22}+v_{33}+2v_{23}$ and $v_{22}+v_{33}-2v_{23}$ all vanish at $x_0$. Thus, $v_{22k}=v_{33k}=v_{23k}=0$. Since $\sum_iv_iv_{ik}=0$ for every $k$, differentiation of \eqref{Equation-v} implies $$v_{11k}=-v_{22k}-v_{33k}=0.$$ It follows from \eqref{Formula-Third derivative commutation} that $v_{211}=-v_2$ and $v_{311}=-v_3$. Consequently,
		\begin{equation}\label{3-dim-Laplacian-det-rank one}
			\Delta\det(v_{ij})=2v_{11}\sum_k(v_{22k}v_{33k}-v_{23k}^2)=0.
		\end{equation}
		For $i,j\in\{2,3\}$, all $v_{ijk}$ vanish and \eqref{3-dim-4th derivative term} becomes $$\sum_kv_{ijkk}=2v_iv_j+2\lambda\delta_{ij}.$$ We therefore obtain
		\begin{align*}
			\Delta\bigg(\sum_{i,j}v_iv_jv_{ij}\bigg)
			=&\sum_{i,j=2}^3v_iv_j\sum_kv_{ijkk}+4v_{11}(v_2v_{211}+v_3v_{311})+2v_{11}^3\\
			=&~2(v_2^2+v_3^2)^2+2\lambda(v_2^2+v_3^2)-4v_{11}(v_2^2+v_3^2)+2v_{11}^3\\
			=&~2v_{11}\big(v_{11}^2-(v_2^2+v_3^2)\big),
		\end{align*}
		where the last equality uses $\lambda=v_{11}-v_2^2-v_3^2$. By strict horo-convexity and \eqref{3-dim-Inequality-phi-t_ast}, the function $\det(v_{ij})-\sum_{i,j}v_iv_jv_{ij}$ also attains its minimum value $0$ at $x_0$. Together with \eqref{3-dim-Laplacian-det-rank one}, this leads to $$2v_{11}\big(v_{11}^2-(v_2^2+v_3^2)\big)\leq 0.$$ As $v_{11}>0$, we conclude that $v_2^2+v_3^2\geq v_{11}^2$ and therefore $$\lambda=v_{11}-v_2^2-v_3^2\leq v_{11}-v_{11}^2\leq\frac{1}{4},$$ again contradicting Lemma \ref{Lemma-Eigenvalue lower bound}.
		
		We have now proved that $\nabla^2v_{t_\ast}>0$ in $\Omega_{t_\ast}$, so $\phi_{t_\ast}$ is well defined. By \eqref{3-dim-Inequality-phi-t_ast}, $$\phi_{t_\ast}\leq\max_{\partial\Omega_{t_\ast}}(\kappa_{1,t_\ast}\kappa_{2,t_\ast})^{-1}<1$$ and hence both conditions hold at $t=t_\ast$.
		
		As in the proof of Theorem \ref{Theorem-2-dimensional}, continuity, the boundary convexity estimate in Proposition \ref{Proposition-Boundary convexity estimate} and the boundary behavior of $\phi_t$ stated in \eqref{Boundary behavior-higher dim} ensure that both conditions remain valid for times slightly below $t_\ast$ if $t_\ast>0$. This contradicts the definition of $t_\ast$. Hence, $t_\ast=0$ and $\nabla^2v>0$ in $\Omega$.
		
		Finally, the argument at the end of Section \ref{Subsection-Proof of Theorem 1.1} shows that $u$ has a unique critical point, its maximum point. Every $\Sigma_t$ with $0<t<\max_\Omega u$ is a smooth closed surface. Its strict convexity follows from \eqref{Equality-Second fundamental form}. Moreover, \eqref{Estimate-Auxiliary function-3-dim} and \eqref{Estimate-Boundary-3-dim} show that $$\frac{1}{K}\leq\phi\leq\max_{\partial\Omega}(\kappa_1\kappa_2)^{-1}.$$ Equivalently, $$K\geq\min_{\partial\Omega}(\kappa_1\kappa_2)>1.$$ This completes the proof in dimension three.
	\end{proof}
	
	\section{The higher-dimensional case}\label{Section-n}
	
	In this section, we prove Theorem \ref{Theorem-Higher-dimensional} for $n\geq 4$. Let $u$ be the first Dirichlet eigenfunction satisfying \eqref{Equation-Dirichlet problem} on a bounded, smooth and strictly horo-convex domain $\Omega\subset\mathbb{H}^n$. The proof relies on a minimum principle for $e^{-v}(1-\phi)$, where $v=-\log u$ and $\phi$ is the auxiliary function defined in \eqref{Auxiliary function-Higher dimensional}. We then follow the same domain deformation strategy as in dimensions two and three.
	
	Throughout this section, Latin indices $i,j,k,l,m,p,q$ range from $1$ to $n$ and $a$ from $1$ to $2$, while Greek indices $\alpha,\beta$ range from $3$ to $n-1$.
	
	\subsection{A minimum principle}\label{Subsection-Minimum principle}
	
	The following proposition will provide a uniform upper bound for $\phi$ along the domain deformation.
	\begin{proposition}\label{Proposition-Maximum principle-n-dim}
		Suppose that $\nabla^2v>0$ and $\phi<1$ in $\Omega$. Then, $e^{-v}(1-\phi)$ has no local minimum in $\Omega$.
	\end{proposition}
	\begin{proof}
		Suppose to the contrary that $e^{-v}(1-\phi)$ attains a local minimum at some $x_0\in\Omega$. Then, there exists a neighborhood $\mathcal{O}\subset\Omega$ of $x_0$ such that $$e^{-v(x)}\big(1-\phi(x)\big)\geq e^{-v(x_0)}\big(1-\phi(x_0)\big)\quad\mathrm{for\ all}\ x\in\mathcal{O}.$$ If $\nabla v(x_0)=0$, then $\phi(x_0)=0$, while $\nabla^2v>0$ implies that $v$ has a strict local minimum at $x_0$. Taking $\mathcal{O}$ sufficiently small, we therefore have $$e^{-v(x)}\big(1-\phi(x)\big)\leq e^{-v(x)}<e^{-v(x_0)}\quad\mathrm{for}\ x\in\mathcal{O}\setminus\{x_0\},$$ which is a contradiction. Hence, $\nabla v(x_0)\neq 0$ and $0<\phi(x_0)<1$.
		
		We now choose suitable coordinates at $x_0$ to simplify the subsequent calculations. We first take an orthonormal basis whose last vector is $\nabla v/|\nabla v|$ and rotate the tangent directions to diagonalize the tangential submatrix $(v^{ij})_{1\leq i,j\leq n-1}$ of the inverse Hessian. Choose unit eigenvectors corresponding to the two largest eigenvalues as the first two basis vectors $\xi,\eta$. This pair attains the maximum in \eqref{Auxiliary function-Higher dimensional}. With this choice, the tangential components of $(\nabla^2v)^{-1}\xi$ and $(\nabla^2v)^{-1}\eta$ are parallel to $\xi,\eta$, respectively, while their normal components are parallel to $\nabla v$. Thus, both vectors belong to the three-dimensional space spanned by $\xi,\eta$ and $\nabla v$. The gradient does not belong to the plane spanned by $(\nabla^2v)^{-1}\xi$ and $(\nabla^2v)^{-1}\eta$. Otherwise, $(\nabla^2v)\nabla v$ would be a linear combination of $\xi,\eta$ and hence orthogonal to $\nabla v$, contradicting the positive definiteness of the Hessian and $\nabla v(x_0)\neq 0$. Then, we choose an orthonormal basis of this plane as the first two coordinate directions and take the last coordinate direction along the component of $\nabla v$ orthogonal to it. These three directions span the same space as $\xi,\eta$ and $\nabla v$.  Taking the remaining coordinate directions in its orthogonal complement, we obtain orthonormal coordinates at $x_0$ such that
		\begin{equation}\label{n-dim-Coordinate condition}
			v_n>0\quad \mathrm{and}\ v_\alpha=0.
		\end{equation}
		Moreover, the first two coordinate directions are linear combinations of $(\nabla^2v)^{-1}\xi$ and $(\nabla^2v)^{-1}\eta$. Applying the Hessian to either direction therefore yields a linear combination of $\xi,\eta$. Since the remaining coordinate directions are orthogonal to $\xi,\eta$, we have
		\begin{equation}\label{n-dim-Coordinate condition'}
			v_{1\alpha}=0\quad\mathrm{and}\ v_{2\alpha}=0.
		\end{equation}
		Finally, we rotate the first two coordinate directions and the remaining $n-3$ coordinate directions separately to diagonalize the corresponding Hessian submatrices. The conditions \eqref{n-dim-Coordinate condition} and \eqref{n-dim-Coordinate condition'} remain unchanged and the rotated frame also satisfies
		\begin{equation}\label{n-dim-Coordinate condition''}
			v_{12}=0\quad \mathrm{and}\ v_{\alpha\beta}=0\ (\alpha\neq\beta).
		\end{equation}
		
		All calculations below are performed at $x_0$ and the coordinate directions are extended by parallel translation along radial geodesics. We first compute the value of $\phi$. The vectors with components $v_{i1}$ and $v_{i2}$ span the maximizing tangent plane and are therefore orthogonal to $\nabla v$. Using \eqref{n-dim-Coordinate condition} and \eqref{n-dim-Coordinate condition''}, we obtain
		\begin{equation}\label{n-dim-Coordinate condition'''}
			v_1=-\frac{v_{1n}}{v_{11}}v_n\quad\mathrm{and}\ v_2=-\frac{v_{2n}}{v_{22}}v_n.
		\end{equation}
		From \eqref{n-dim-Coordinate condition}-\eqref{n-dim-Coordinate condition'''}, the definition \eqref{Auxiliary function-Higher dimensional} takes the form
		\begin{equation}\label{n-dim-Coordinate condition''''}
			\phi=\frac{v_n^2}{v_{11}v_{22}}.
		\end{equation}
		
		To differentiate $\phi$, we use the construction described in Section \ref{Subsection-Auxiliary functions}. At $x_0$, the tangential projections of the first two coordinate directions span the maximizing tangent plane. Hence, the smooth function in \eqref{Auxiliary function-Higher dimensional'} agrees with the original auxiliary function \eqref{Auxiliary function-Higher dimensional} at $x_0$ and is bounded above by it nearby. We continue to denote this smooth function by $\phi$. Then, $e^{-v}(1-\phi)$ still attains a local minimum at $x_0$, so
		\begin{equation}\label{n-dim-1st condition'}
			\phi_k=-(1-\phi)v_k.
		\end{equation}
		Moreover,
		\begin{equation}\label{n-dim-2nd condition}
			\Delta\big(e^{-v}(1-\phi)\big)\geq 0.
		\end{equation}
		
		Differentiating \eqref{Auxiliary function-Higher dimensional'} using $$(v^{ij})_k=-\sum_{l,m}v^{il}v^{jm}v_{lmk}$$ and simplifying by \eqref{n-dim-Coordinate condition}-\eqref{n-dim-Coordinate condition'''}, we obtain
		\begin{equation}\label{n-dim-1st condition}
			\begin{aligned}
				\phi_k
				=&~\frac{2v_n}{v_{11}v_{22}}v_{nk}-\frac{v_n^2}{v_{11}^2v_{22}}v_{11k}-\frac{v_n^2}{v_{11}v_{22}^2}v_{22k}\\
				=&~\phi\bigg(\frac{2v_{nk}}{v_n}-\frac{v_{11k}}{v_{11}}-\frac{v_{22k}}{v_{22}}\bigg).
			\end{aligned}
		\end{equation}
		For the second derivatives, we also use the identity $$(v^{ij})_{kk}=2\sum_{l,m,p,q}v^{il}v^{mp}v^{qj}v_{lmk}v_{pqk}-\sum_{l,m}v^{il}v^{mj}v_{lmkk}.$$ Differentiating \eqref{Auxiliary function-Higher dimensional'} once more and summing over $k$, we find
		\begin{equation}\label{n-dim-Equality-F1}
			\begin{aligned}
				\frac{\Delta\phi}{\phi}
				=&~\frac{2}{v_n^2}\sum_{i=3}^n\sum_k(v_{ik}^2+v_iv_{ikk})-\frac{4}{v_n^2}\sum_kv_{nk}^2-\frac{1}{v_{11}}\sum_kv_{11kk}-\frac{1}{v_{22}}\sum_kv_{22kk}\\
				&+\sum_k\bigg(\frac{v_{11k}^2}{v_{11}^2}+\frac{v_{22k}^2}{v_{22}^2}+\frac{2v_{12k}^2}{v_{11}v_{22}}\bigg)+\frac{|\nabla\phi|^2}{\phi^2}\\
				&+2\sum_{a=1}^2\sum_k\frac{1}{v_{aa}}\bigg(\frac{v^{nn}}{v_n^2}\bigg(\sum_iv_iv_{aik}+\sum_iv_{ai}v_{ik}\bigg)^2\\
				&\hspace{2.8cm}+\frac{2}{v_n}\sum_{\alpha=3}^{n-1}v^{n\alpha}\bigg(\sum_iv_iv_{aik}+\sum_iv_{ai}v_{ik}\bigg)\bigg(v_{a\alpha k}-\frac{v_{an}v_{\alpha k}}{v_n}\bigg)\\
				&\hspace{2.8cm}+\sum_{\alpha,\beta=3}^{n-1}v^{\alpha\beta}\bigg(v_{a\alpha k}-\frac{v_{an}v_{\alpha k}}{v_n}\bigg)\bigg(v_{a\beta k}-\frac{v_{an}v_{\beta k}}{v_n}\bigg)\bigg).
			\end{aligned}
		\end{equation}
		
		To evaluate the fourth-order derivatives in \eqref{n-dim-Equality-F1}, we first differentiate equation \eqref{Equation-v} and employ the commutation formula \eqref{Formula-Third derivative commutation} to obtain $$\sum_kv_{ikk}=(\Delta v)_i-(n-1)v_i=2\sum_kv_kv_{ik}-(n-1)v_i.$$ The second derivatives of \eqref{Equation-v} can be written as
		\begin{align*}
			(\Delta v)_{ij}
			=&~2\sum_kv_{ki}v_{kj}+2\sum_kv_kv_{kij}\\
			=&~2\sum_kv_{ki}v_{kj}+2\sum_kv_kv_{ijk}+2v_iv_j-2\delta_{ij}\sum_kv_k^2,
		\end{align*}
		where we have used \eqref{Formula-Ricci commutation} in the last equality. It follows from \eqref{Formula-Fourth derivative commutation} that
		\begin{equation}\label{n-dim-4th derivative term}
			\begin{aligned}
				\sum_kv_{ijkk}
				=&~(\Delta v)_{ij}-2nv_{ij}+2\delta_{ij}\Delta v\\
				=&~2\sum_kv_kv_{ijk}+2\sum_kv_{ki}v_{kj}+2v_iv_j-2nv_{ij}+2\lambda\delta_{ij}.
			\end{aligned}
		\end{equation}
		Using \eqref{n-dim-Coordinate condition'}, \eqref{n-dim-Coordinate condition''} and \eqref{n-dim-4th derivative term}, we calculate
		\begin{align*}
			&-\frac{1}{v_{11}}\sum_kv_{11kk}-\frac{1}{v_{22}}\sum_kv_{22kk}\\
			=&-2\sum_kv_k\bigg(\frac{v_{11k}}{v_{11}}+\frac{v_{22k}}{v_{22}}\bigg)-2\bigg(v_{11}+v_{22}+\frac{v_{1n}^2}{v_{11}}+\frac{v_{2n}^2}{v_{22}}\bigg)\\
			&-2\bigg(\frac{v_1^2}{v_{11}}+\frac{v_2^2}{v_{22}}\bigg)+4n-2\lambda\bigg(\frac{1}{v_{11}}+\frac{1}{v_{22}}\bigg).
		\end{align*}
		Note that $$\sum_{i=3}^n\sum_kv_{ik}^2-2\sum_kv_{nk}^2=\sum_{\alpha=3}^{n-1}(v_{\alpha\alpha}^2+v_{n\alpha}^2)-\sum_kv_{nk}^2.$$ Thus, 
		\begin{equation}\label{n-dim-Equality-F1'}
			\begin{aligned}
				\frac{\Delta\phi}{\phi}
				=&~4n-2\bigg(v_{11}+v_{22}+\frac{v_{1n}^2}{v_{11}}+\frac{v_{2n}^2}{v_{22}}\bigg)-2\bigg(\frac{v_1^2}{v_{11}}+\frac{v_2^2}{v_{22}}\bigg)-2\lambda\bigg(\frac{1}{v_{11}}+\frac{1}{v_{22}}\bigg)\\
				&+\frac{2}{v_n^2}\sum_{i=3}^n\sum_kv_iv_{ikk}-2\sum_kv_k\bigg(\frac{v_{11k}}{v_{11}}+\frac{v_{22k}}{v_{22}}\bigg)\\
				&+\sum_k\bigg(\frac{v_{11k}^2}{v_{11}^2}+\frac{v_{22k}^2}{v_{22}^2}+\frac{2v_{12k}^2}{v_{11}v_{22}}\bigg)+\frac{|\nabla\phi|^2}{\phi^2}\\
				&+\frac{2}{v_n^2}\sum_{\alpha=3}^{n-1}(v_{\alpha\alpha}^2+v_{n\alpha}^2)-\frac{2}{v_n^2}\sum_kv_{nk}^2\\
				&+2\sum_{a=1}^2\sum_k\frac{1}{v_{aa}}\bigg(\frac{v^{nn}}{v_n^2}\bigg(\sum_iv_iv_{aik}+\sum_iv_{ai}v_{ik}\bigg)^2\\
				&\hspace{2.8cm}+\frac{2}{v_n}\sum_{\alpha=3}^{n-1}v^{n\alpha}\bigg(\sum_iv_iv_{aik}+\sum_iv_{ai}v_{ik}\bigg)\bigg(v_{a\alpha k}-\frac{v_{an}v_{\alpha k}}{v_n}\bigg)\\
				&\hspace{2.8cm}+\sum_{\alpha,\beta=3}^{n-1}v^{\alpha\beta}\bigg(v_{a\alpha k}-\frac{v_{an}v_{\alpha k}}{v_n}\bigg)\bigg(v_{a\beta k}-\frac{v_{an}v_{\beta k}}{v_n}\bigg)\bigg).
			\end{aligned}
		\end{equation}
		
		We next handle the third-order derivative terms. For the linear terms in the third-order derivatives, we have
		\begin{equation}\label{n-dim-3rd derivative identity}
			\begin{aligned}
				&~\frac{2}{v_n^2}\sum_{i=3}^n\sum_kv_iv_{ikk}-2\sum_kv_k\bigg(\frac{v_{11k}}{v_{11}}+\frac{v_{22k}}{v_{22}}\bigg)\\
				=&~\frac{4}{v_n}\sum_kv_kv_{nk}-2(n-1)-2\sum_kv_k\bigg(\frac{v_{11k}}{v_{11}}+\frac{v_{22k}}{v_{22}}\bigg)\\
				=&~2\sum_kv_k\frac{\phi_k}{\phi}-2(n-1),
			\end{aligned}
		\end{equation}
		where the last equality follows from \eqref{n-dim-1st condition}. For the quadratic terms, we first note from $v_{12}=0$ and the positive definiteness of $\nabla^2v$ that
		\begin{equation}\label{n-dim-Condition-positive definite}
			v_{nn}-\frac{v_{1n}^2}{v_{11}}-\frac{v_{2n}^2}{v_{22}}>0.
		\end{equation}
		Using \eqref{n-dim-Coordinate condition'} and \eqref{n-dim-Coordinate condition''}, the identities $\sum_jv_{ij}v^{jk}=\delta_{ik}$ imply $$v^{n\alpha}=-\frac{\sum_{\beta=3}^{n-1}v^{\alpha\beta}v_{\beta n}}{v_{nn}-\frac{v_{1n}^2}{v_{11}}-\frac{v_{2n}^2}{v_{22}}}$$ and $$v^{nn}=\frac{1}{v_{nn}-\frac{v_{1n}^2}{v_{11}}-\frac{v_{2n}^2}{v_{22}}}+\frac{\sum_{\alpha,\beta=3}^{n-1}v^{\alpha\beta}v_{\alpha n}v_{\beta n}}{\big(v_{nn}-\frac{v_{1n}^2}{v_{11}}-\frac{v_{2n}^2}{v_{22}}\big)^2}.$$ Consequently, the last term in \eqref{n-dim-Equality-F1'} can be written as
		\begin{equation}\label{n-dim-3rd derivative identity1}
			\begin{aligned}
				&2\sum_{a=1}^2\sum_k\frac{1}{v_{aa}}\bigg(\frac{v^{nn}}{v_n^2}\bigg(\sum_iv_iv_{aik}+\sum_iv_{ai}v_{ik}\bigg)^2\\
				&\hspace{2.8cm}+\frac{2}{v_n}\sum_{\alpha=3}^{n-1}v^{n\alpha}\bigg(\sum_iv_iv_{aik}+\sum_iv_{ai}v_{ik}\bigg)\bigg(v_{a\alpha k}-\frac{v_{an}v_{\alpha k}}{v_n}\bigg)\\
				&\hspace{2.8cm}+\sum_{\alpha,\beta=3}^{n-1}v^{\alpha\beta}\bigg(v_{a\alpha k}-\frac{v_{an}v_{\alpha k}}{v_n}\bigg)\bigg(v_{a\beta k}-\frac{v_{an}v_{\beta k}}{v_n}\bigg)\bigg)\\
				=&~\frac{2}{v_n^2\big(v_{nn}-\frac{v_{1n}^2}{v_{11}}-\frac{v_{2n}^2}{v_{22}}\big)}\sum_{a=1}^2\sum_k\frac{1}{v_{aa}}\bigg(\sum_iv_iv_{aik}+\sum_iv_{ai}v_{ik}\bigg)^2\\
				&+2\sum_{a=1}^2\sum_k\frac{1}{v_{aa}}\sum_{\alpha,\beta=3}^{n-1}v^{\alpha\beta}\\
				&~\cdot\bigg(v_{a\alpha k}-\frac{v_{an}v_{\alpha k}}{v_n}-\frac{v_{\alpha n}}{v_n\big(v_{nn}-\frac{v_{1n}^2}{v_{11}}-\frac{v_{2n}^2}{v_{22}}\big)}\bigg(\sum_iv_iv_{aik}+\sum_iv_{ai}v_{ik}\bigg)\bigg)\\
				&\cdot~\bigg(v_{a\beta k}-\frac{v_{an}v_{\beta k}}{v_n}-\frac{v_{\beta n}}{v_n\big(v_{nn}-\frac{v_{1n}^2}{v_{11}}-\frac{v_{2n}^2}{v_{22}}\big)}\bigg(\sum_iv_iv_{aik}+\sum_iv_{ai}v_{ik}\bigg)\bigg).
			\end{aligned}
		\end{equation}
		The second term on the right-hand side is nonnegative, since $(v^{\alpha\beta})_{\alpha,\beta=3}^{n-1}$ is positive definite. For the terms involving $v_{11k}^2$, $v_{22k}^2$, $v_{nk}^2$ and $|\nabla\phi|^2$, we apply \eqref{n-dim-Coordinate condition'''}, \eqref{n-dim-1st condition'} and \eqref{n-dim-1st condition} to write
		\begin{equation}\label{n-dim-3rd derivaitve identity2}
			\begin{aligned}
				&\sum_k\bigg(\frac{v_{11k}^2}{v_{11}^2}+\frac{v_{22k}^2}{v_{22}^2}\bigg)+\frac{|\nabla\phi|^2}{\phi^2}-\frac{2}{v_n^2}\sum_kv_{nk}^2\\
				=&~2\sum_k\bigg(\frac{v_{11k}^2}{v_{11}^2}+\frac{v_{22k}^2}{v_{22}^2}\bigg)+\frac{2}{v_{11}v_{22}}\sum_kv_{11k}v_{22k}\\
				&-\frac{4}{v_n}\sum_kv_{nk}\bigg(\frac{v_{11k}}{v_{11}}+\frac{v_{22k}}{v_{22}}\bigg)+\frac{2}{v_n^2}\sum_kv_{nk}^2\\
				=&~\frac{1}{2}\sum_k\bigg(\frac{v_{11k}}{v_{11}}-\frac{v_{22k}}{v_{22}}\bigg)^2+\frac{3}{2}\sum_k\bigg(\frac{v_{11k}}{v_{11}}+\frac{v_{22k}}{v_{22}}\bigg)^2\\
				&-\frac{4}{v_n}\sum_kv_{nk}\bigg(\frac{v_{11k}}{v_{11}}+\frac{v_{22k}}{v_{22}}\bigg)+\frac{2}{v_n^2}\sum_kv_{nk}^2\\
				=&~\frac{1}{2}\sum_k\bigg(\frac{v_{11k}}{v_{11}}-\frac{v_{22k}}{v_{22}}\bigg)^2+\frac{3}{2}\sum_k\bigg(\frac{2v_{nk}}{3v_n}+\frac{1-\phi}{\phi}v_k\bigg)^2-\frac{2}{3v_n^2}\sum_kv_{nk}^2\\
				=&~\frac{1}{2}\sum_k\bigg(\frac{v_{11k}}{v_{11}}-\frac{v_{22k}}{v_{22}}\bigg)^2+\frac{2(1-\phi)}{\phi}\bigg(v_{nn}-\frac{v_{1n}^2}{v_{11}}-\frac{v_{2n}^2}{v_{22}}\bigg)+\frac{3(1-\phi)^2}{2\phi^2}|\nabla v|^2.
			\end{aligned}
		\end{equation}
		The last two terms are nonnegative by \eqref{n-dim-Condition-positive definite} and $0<\phi(x_0)<1$. Substituting \eqref{n-dim-3rd derivative identity}-\eqref{n-dim-3rd derivaitve identity2} into \eqref{n-dim-Equality-F1'} and dropping the nonnegative terms including those containing $v_{n\alpha}^2$, we obtain
		\begin{equation}\label{n-dim-Inequality-F1}
			\begin{aligned}
				\frac{\Delta\phi-2\sum_kv_k\phi_k}{\phi}
				\geq&~2n+2-2\bigg(v_{11}+v_{22}+\frac{v_{1n}^2}{v_{11}}+\frac{v_{2n}^2}{v_{22}}\bigg)-2\bigg(\frac{v_1^2}{v_{11}}+\frac{v_2^2}{v_{22}}\bigg)\\
				&-2\lambda\bigg(\frac{1}{v_{11}}+\frac{1}{v_{22}}\bigg)+\frac{1}{2}\sum_k\bigg(\frac{v_{11k}}{v_{11}}-\frac{v_{22k}}{v_{22}}\bigg)^2+\frac{2}{v_{11}v_{22}}\sum_kv_{12k}^2\\
				&+\frac{2}{v_n^2}\sum_{\alpha=3}^{n-1}v_{\alpha\alpha}^2\\
				&+\frac{2}{v_n^2\big(v_{nn}-\frac{v_{1n}^2}{v_{11}}-\frac{v_{2n}^2}{v_{22}}\big)}\sum_{a=1}^2\sum_k\frac{1}{v_{aa}}\bigg(\sum_iv_iv_{aik}+\sum_iv_{ai}v_{ik}\bigg)^2.
			\end{aligned}
		\end{equation}
		
		We further estimate the third-order derivative terms in \eqref{n-dim-Inequality-F1} as in the three-dimensional proof. The ingredient is to collect the terms involving $\frac{v_nv_{111}}{v_{11}}$, $\frac{v_nv_{112}}{v_{11}}$ and $\frac{v_nv_{11n}}{v_{11}}$. By \eqref{n-dim-1st condition'} and \eqref{n-dim-1st condition}, we express $v_{22k}$ in terms of $v_{11k}$ as
		\begin{equation}\label{n-dim-3rd derivative relation}
			v_{22k}=v_{22}\bigg(\frac{2v_{nk}}{v_n}-\frac{v_{11k}}{v_{11}}+\frac{1-\phi}{\phi}v_k\bigg).
		\end{equation}
		Thus, we retain only the terms with $k=1,2,n$ and use \eqref{n-dim-Coordinate condition'''} and \eqref{n-dim-Coordinate condition''''} to obtain
		\begin{equation}\label{n-dim-Estimate-3rd derivaitve1}
			\begin{aligned}
				\frac{1}{2}\sum_k\bigg(\frac{v_{11k}}{v_{11}}-\frac{v_{22k}}{v_{22}}\bigg)^2
				=&~\frac{2}{v_n^2}\sum_k\bigg(\frac{v_nv_{11k}}{v_{11}}-v_{nk}-\frac{1-\phi}{2\phi}v_nv_k\bigg)^2\\
				\geq&~\frac{2}{v_n^2}\bigg(\bigg(\frac{v_nv_{111}}{v_{11}}-\bigg(1-\frac{1-\phi}{2}v_{22}\bigg)v_{1n}\bigg)^2\\
				&\hspace{0.8cm}+\bigg(\frac{v_nv_{112}}{v_{11}}-\bigg(1-\frac{1-\phi}{2}v_{11}\bigg)v_{2n}\bigg)^2\\
				&\hspace{0.8cm}+\bigg(\frac{v_nv_{11n}}{v_{11}}-v_{nn}-\frac{1-\phi}{2}v_{11}v_{22}\bigg)^2\bigg).
			\end{aligned}
		\end{equation}
		On the other hand, \eqref{Formula-Ricci commutation}, \eqref{Formula-Third derivative commutation}, \eqref{n-dim-Coordinate condition'''} and the expression for $v_{22k}$ imply $$v_{121}=v_{112}+\frac{v_nv_{2n}}{v_{22}}$$ and $$v_{122}=v_{22}\bigg(\frac{2v_{1n}}{v_n}-\frac{v_{111}}{v_{11}}-\frac{1-\phi}{\phi}\frac{v_nv_{1n}}{v_{11}}\bigg)+\frac{v_nv_{1n}}{v_{11}}.$$ After the terms with $k=3,4,\cdots,n-1$ are omitted, \eqref{n-dim-Coordinate condition''''} yields
		\begin{equation}\label{n-dim-Estimate-3rd derivaitve2}
			\begin{aligned}
				\frac{2}{v_{11}v_{22}}\sum_kv_{12k}^2
				\geq&~\frac{2}{v_{11}v_{22}}(v_{121}^2+v_{122}^2+v_{12n}^2)\\
				=&~\frac{2}{v_n^2}\bigg(\frac{v_{11}}{v_{22}}\bigg(\frac{v_nv_{112}}{v_{11}}+\phi v_{2n}\bigg)^2\\
				&\hspace{0.8cm}+\frac{v_{22}}{v_{11}}\bigg(\frac{v_nv_{111}}{v_{11}}-(2+\phi-(1-\phi)v_{22})v_{1n}\bigg)^2+\phi v_{12n}^2\bigg).
			\end{aligned}
		\end{equation}
		For the last term in \eqref{n-dim-Inequality-F1}, we restrict the sum to $k=a$ and deduce from Young's inequality that
		\begin{equation}\label{n-dim-Estimate-3rd derivaitve3}
			\begin{aligned}
				&~\frac{2}{v_n^2\big(v_{nn}-\frac{v_{1n}^2}{v_{11}}-\frac{v_{2n}^2}{v_{22}}\big)}\sum_{a=1}^2\sum_k\frac{1}{v_{aa}}\bigg(\sum_iv_iv_{aik}+\sum_iv_{ai}v_{ik}\bigg)^2\\
				\geq&~\frac{2}{v_n^2\big(v_{nn}-\frac{v_{1n}^2}{v_{11}}-\frac{v_{2n}^2}{v_{22}}\big)}\sum_{a=1}^2\frac{1}{v_{aa}}\bigg(\sum_iv_iv_{aia}+\sum_iv_{ai}^2\bigg)^2\\
				\geq&~\frac{4}{v_n^2}\sum_{a=1}^2\bigg(\sum_iv_iv_{aia}+\sum_iv_{ai}^2\bigg)-\frac{2(v_{11}+v_{22})}{v_n^2}\bigg(v_{nn}-\frac{v_{1n}^2}{v_{11}}-\frac{v_{2n}^2}{v_{22}}\bigg).
			\end{aligned}
		\end{equation}
		Here, we have used \eqref{n-dim-Condition-positive definite}. To calculate the terms linear in the third-order derivatives in \eqref{n-dim-Estimate-3rd derivaitve3}, we express them in terms of $v_{111}$, $v_{112}$ and $v_{11n}$ as in dimension three. Taking $k=2,n$ in \eqref{n-dim-3rd derivative relation} and using \eqref{n-dim-Coordinate condition'''} and the commutation formula \eqref{Formula-Third derivative commutation}, we have
		\begin{align*}
			v_{222}=&~v_{22}\bigg(\frac{2v_{2n}}{v_n}-\frac{v_{112}}{v_{11}}-\frac{1-\phi}{\phi}\frac{v_nv_{2n}}{v_{22}}\bigg),\\
			v_{2n2}=&~v_{22}\bigg(\frac{2v_{nn}}{v_n}-\frac{v_{11n}}{v_{11}}+\frac{1-\phi}{\phi}v_n\bigg)-v_n\\
			\intertext{and}
			v_{1n1}=&~v_{11n}-v_n.
		\end{align*} 
		It follows from these identities and the expressions for $v_{121}$ and $v_{122}$ that
		\begin{equation}\label{n-dim-3rd derivative identity3}
			\begin{aligned}
				&\sum_{a=1}^2\bigg(\sum_iv_iv_{aia}+\sum_iv_{ai}^2\bigg)\\
				=&~v_1(v_{111}+v_{122})+v_2(v_{121}+v_{222})+v_n(v_{1n1}+v_{2n2})\\
				&+v_{11}^2+v_{22}^2+v_{1n}^2+v_{2n}^2\\
				=&-\frac{v_nv_{1n}}{v_{11}}\bigg(v_{111}+v_{22}\bigg(\frac{2v_{1n}}{v_n}-\frac{v_{111}}{v_{11}}-\frac{1-\phi}{\phi}\frac{v_nv_{1n}}{v_{11}}\bigg)+\frac{v_nv_{1n}}{v_{11}}\bigg)\\
				&-\frac{v_nv_{2n}}{v_{22}}\bigg(v_{112}+\frac{v_nv_{2n}}{v_{22}}+v_{22}\bigg(\frac{2v_{2n}}{v_n}-\frac{v_{112}}{v_{11}}-\frac{1-\phi}{\phi}\frac{v_nv_{2n}}{v_{22}}\bigg)\bigg)\\
				&+v_n\bigg(v_{11n}+v_{22}\bigg(\frac{2v_{nn}}{v_n}-\frac{v_{11n}}{v_{11}}+\frac{1-\phi}{\phi}v_n\bigg)-2v_n\bigg)+v_{11}^2+v_{22}^2+v_{1n}^2+v_{2n}^2\\
				=&-\frac{v_nv_{1n}(v_{11}-v_{22})}{v_{11}^2}v_{111}-\frac{v_nv_{2n}(v_{11}-v_{22})}{v_{11}v_{22}}v_{112}+\frac{v_n(v_{11}-v_{22})}{v_{11}}v_{11n}\\
				&+2v_{22}v_{nn}+v_{11}^2+v_{22}^2-2v_n^2+\bigg(1-\frac{v_{22}}{v_{11}}(2+\phi)\bigg)v_{1n}^2-\bigg(1+\frac{v_{11}}{v_{22}}\phi\bigg)v_{2n}^2\\
				&+\frac{v_{22}(1-\phi)}{\phi}|\nabla v|^2,
			\end{aligned}
		\end{equation}
		in which \eqref{n-dim-Coordinate condition}-\eqref{n-dim-Coordinate condition'''} were used in the first two equalities and \eqref{n-dim-Coordinate condition''''} in the last equality.
		
		We now substitute \eqref{n-dim-Estimate-3rd derivaitve1}-\eqref{n-dim-3rd derivative identity3} into \eqref{n-dim-Inequality-F1} to obtain
		\begin{align*}
			&~\Delta\phi-2\sum_kv_k\phi_k\\
			\geq&~\phi\bigg(2n+2-2\bigg(v_{11}+v_{22}+\frac{v_{1n}^2}{v_{11}}+\frac{v_{2n}^2}{v_{22}}\bigg)-2\bigg(\frac{v_1^2}{v_{11}}+\frac{v_2^2}{v_{22}}\bigg)-2\lambda\bigg(\frac{1}{v_{11}}+\frac{1}{v_{22}}\bigg)\bigg)\\
			&+\frac{2}{v_{11}v_{22}}\bigg(\frac{v_{11}+v_{22}}{v_{11}}\bigg(\frac{v_nv_{111}}{v_{11}}\bigg)^2+\frac{v_{11}+v_{22}}{v_{22}}\bigg(\frac{v_nv_{112}}{v_{11}}\bigg)^2+\bigg(\frac{v_nv_{11n}}{v_{11}}\bigg)^2\\
			&\hspace{1.8cm}-\frac{2v_{1n}}{v_{11}}\bigg(2v_{11}+(1+\phi)v_{22}-\frac{1-\phi}{2}v_{22}(v_{11}+2v_{22})\bigg)\frac{v_nv_{111}}{v_{11}}\\
			&\hspace{1.8cm}-\frac{(1-\phi)v_{11}(2-v_{22})v_{2n}}{v_{22}}\frac{v_nv_{112}}{v_{11}}\\
			&\hspace{1.8cm}+2\bigg(v_{11}-v_{22}-v_{nn}-\frac{1-\phi}{2}v_{11}v_{22}\bigg)\frac{v_nv_{11n}}{v_{11}}\\
			&\hspace{1.8cm}+\phi v_{12n}^2+\sum_{\alpha=3}^{n-1}v_{\alpha\alpha}^2+2v_{11}^2+2v_{22}^2+v_{nn}^2-4v_n^2+(3v_{22}-v_{11})v_{nn}\\
			&\hspace{1.8cm}+(1-\phi)v_{11}v_{22}(v_{nn}+2v_{22})+\frac{1}{4}(1-\phi)^2v_{11}^2v_{22}^2\\
			&\hspace{1.8cm}+\bigg(4-(1-\phi)v_{22}+\frac{v_{22}}{v_{11}}\big(1+\phi-(1-\phi)v_{22}\big)^2+\frac{1}{4}(1-\phi)^2v_{22}^2\bigg)v_{1n}^2\\
			&\hspace{1.8cm}+\bigg((1-\phi)v_{11}+\frac{v_{11}}{v_{22}}(1-\phi)^2+\frac{1}{4}(1-\phi)^2v_{11}^2\bigg)v_{2n}^2\bigg).
		\end{align*}
		As in the three-dimensional proof, for the terms with $v_{111}$, $v_{112}$ and $v_{11n}$, one has
		\begin{align*}
			&~\frac{v_{11}+v_{22}}{v_{11}}\bigg(\frac{v_nv_{111}}{v_{11}}\bigg)^2+\frac{v_{11}+v_{22}}{v_{22}}\bigg(\frac{v_nv_{112}}{v_{11}}\bigg)^2+\bigg(\frac{v_nv_{11n}}{v_{11}}\bigg)^2\\
			&-\frac{2v_{1n}}{v_{11}}\bigg(2v_{11}+(1+\phi)v_{22}-\frac{1-\phi}{2}v_{22}(v_{11}+2v_{22})\bigg)\frac{v_nv_{111}}{v_{11}}\\
			&-\frac{(1-\phi)v_{11}(2-v_{22})v_{2n}}{v_{22}}\frac{v_nv_{112}}{v_{11}}\\
			&+2\bigg(v_{11}-v_{22}-v_{nn}-\frac{1-\phi}{2}v_{11}v_{22}\bigg)\frac{v_nv_{11n}}{v_{11}}\\
			=&~\frac{v_{11}+v_{22}}{v_{11}}\bigg(\frac{v_nv_{111}}{v_{11}}-\frac{2v_{11}+(1+\phi)v_{22}-\frac{1-\phi}{2}v_{22}(v_{11}+2v_{22})}{v_{11}+v_{22}}v_{1n}\bigg)^2\\
			&+\frac{v_{11}+v_{22}}{v_{22}}\bigg(\frac{v_nv_{112}}{v_{11}}-\frac{(1-\phi)v_{11}(2-v_{22})}{2(v_{11}+v_{22})}v_{2n}\bigg)^2\\
			&+\bigg(\frac{v_nv_{11n}}{v_{11}}+v_{11}-v_{22}-v_{nn}-\frac{1-\phi}{2}v_{11}v_{22}\bigg)^2\\
			&-\frac{\big(2v_{11}+(1+\phi)v_{22}-\frac{1-\phi}{2}v_{22}(v_{11}+2v_{22})\big)^2}{v_{11}(v_{11}+v_{22})}v_{1n}^2-\frac{(1-\phi)^2v_{11}^2(2-v_{22})^2}{4v_{22}(v_{11}+v_{22})}v_{2n}^2\\
			&-\bigg(v_{11}-v_{22}-v_{nn}-\frac{1-\phi}{2}v_{11}v_{22}\bigg)^2\\
			\geq&-\frac{\big(2v_{11}+(1+\phi)v_{22}-\frac{1-\phi}{2}v_{22}(v_{11}+2v_{22})\big)^2}{v_{11}(v_{11}+v_{22})}v_{1n}^2-\frac{(1-\phi)^2v_{11}^2(2-v_{22})^2}{4v_{22}(v_{11}+v_{22})}v_{2n}^2\\
			&-\bigg(v_{11}-v_{22}-v_{nn}-\frac{1-\phi}{2}v_{11}v_{22}\bigg)^2.
		\end{align*}
		Dropping the nonnegative term involving $v_{12n}^2$ gives
		\begin{equation}\label{n-dim-Inequality-F2}
			\begin{aligned}
				\Delta\phi-2\sum_kv_k\phi_k
				\geq&~2(n-3)\phi+\frac{2}{v_{11}v_{22}}\sum_{\alpha=3}^{n-1}v_{\alpha\alpha}^2\\
				&+\frac{2(v_{11}+v_{22})}{v_{11}v_{22}}\bigg((1-\phi)\Delta v-\sum_{\alpha=3}^{n-1}v_{\alpha\alpha}\bigg)\\
				&+2(1-\phi)^2\bigg(v_{11}+v_{22}+\frac{v_{1n}^2}{v_{11}}+\frac{v_{2n}^2}{v_{22}}\bigg)\\
				&+\frac{(1-\phi)^2}{2(v_{11}+v_{22})}\bigg(\frac{(v_{22}+2)^2}{v_{11}}v_{1n}^2+\frac{(v_{11}+2)^2}{v_{22}}v_{2n}^2\bigg),
			\end{aligned}
		\end{equation}
		where \eqref{Equation-v}, \eqref{n-dim-Coordinate condition'''} and \eqref{n-dim-Coordinate condition''''} have been used. Since $\nabla^2v>0$ and $0<\phi(x_0)<1$, the terms with $(1-\phi)\Delta v$ and $(1-\phi)^2$ are nonnegative. Thus, \eqref{n-dim-Inequality-F2} yields
		\begin{equation}\label{n-dim-Inequality-F3}
			\begin{aligned}
				&~\Delta\phi-2\sum_kv_k\phi_k\\
				\geq&~2(n-3)\phi+\frac{2}{v_{11}v_{22}}\sum_{\alpha=3}^{n-1}\big(v_{\alpha\alpha}^2-(v_{11}+v_{22})v_{\alpha\alpha}\big)\\
				=&~\frac{2}{v_{11}v_{22}}\sum_{\alpha=3}^{n-1}(v_{\alpha\alpha}-v_{11})(v_{\alpha\alpha}-v_{22})-2(n-3)(1-\phi).
			\end{aligned}
		\end{equation}
		
		We next estimate $(v_{\alpha\alpha}-v_{11})(v_{\alpha\alpha}-v_{22})$ for $\alpha=3,4,\cdots,n-1$. Recall the choice of $\xi,\eta$ as the eigenvectors corresponding to the largest tangential eigenvalues of the inverse Hessian. Thus, the coordinate construction implies $$v^{\alpha\alpha}\leq\frac{\sum_{i,j}v^{ij}v_{ia}v_{ja}}{\sum_iv_{ia}^2}=\frac{v_{aa}}{v_{aa}^2+v_{an}^2}\leq\frac{1}{v_{aa}},\quad a=1,2.$$ Here, the equality is obtained from \eqref{n-dim-Coordinate condition'}, \eqref{n-dim-Coordinate condition''} and $\sum_jv_{ij}v^{jk}=\delta_{ik}$. From \eqref{n-dim-Coordinate condition'} and \eqref{n-dim-Coordinate condition''}, $$v_{\alpha\alpha}v^{\alpha\alpha}+v_{\alpha n}v^{n\alpha}=1,\quad \mathrm{and}\ v_{\alpha\alpha}v^{n\alpha}+v_{\alpha n}v^{nn}=0.$$ Thus, $$v^{\alpha\alpha}=\frac{1}{v_{\alpha\alpha}}+\frac{v_{\alpha n}^2}{v_{\alpha\alpha}^2}v^{nn}\geq\frac{1}{v_{\alpha\alpha}}.$$  Consequently, $v_{\alpha\alpha}\geq\max\{v_{11},v_{22}\}$ and therefore $$(v_{\alpha\alpha}-v_{11})(v_{\alpha\alpha}-v_{22})\geq 0.$$ Returning to \eqref{n-dim-Inequality-F3}, we conclude that
		\begin{equation}\label{n-dim-Inequality-F4}
			\Delta\phi-2\sum_kv_k\phi_k\geq-2(n-3)(1-\phi).
		\end{equation}
		
		Finally, we observe that $$\lambda-2(n-3)>\frac{(n-1)^2}{4}-2(n-3)=\frac{(n-5)^2}{4}\geq 0$$ by Lemma \ref{Lemma-Eigenvalue lower bound}. Equation \eqref{Equation-v} and estimate \eqref{n-dim-Inequality-F4} imply
		\begin{align*}
			\Delta\big(e^{-v}(1-\phi)\big)
			=&~e^{-v}\bigg(-\lambda(1-\phi)-\Delta\phi+2\sum_kv_k\phi_k\bigg)\\
			\leq&-e^{-v}\big(\lambda-2(n-3)\big)(1-\phi)\\
			<&~0,
		\end{align*}
		since $0<\phi(x_0)<1$. This contradicts \eqref{n-dim-2nd condition} and the proof is complete.
	\end{proof}
	
	\subsection{Proof of Theorem \ref{Theorem-Higher-dimensional} for $n\geq 4$}\label{Subsection-Proof of Theorem 1.2-n}
	
	We follow the continuity argument in Sections \ref{Subsection-Proof of Theorem 1.1} and \ref{Subsection-Proof of Theorem 1.2-n=3}. Here, the minimum principle in Proposition \ref{Proposition-Maximum principle-n-dim} provides a uniform upper bound for $\phi_t$ along the domain deformation, which allows us to exclude degeneration of the Hessian.
	
	\begin{proof}
		Let $\{\Omega_t\}_{t\geq 0}$ and $v_t=-\log u_t$ be as in Section \ref{Subsection-Domain deformation}, with $\lambda_t$ the first Dirichlet eigenvalue of $\Omega_t$. Let $\phi_t$ denote the auxiliary function defined by \eqref{Auxiliary function-Higher dimensional} for $v_t$. By Proposition \ref{Proposition-Domain deformation}, there exists $T>0$ such that $\nabla^2v_t>0$ and $\phi_t<1$ throughout $\Omega_t$ whenever $t\geq T$. Let $t_\ast\in[0,T)$ be the smallest time such that both conditions hold for every $t>t_\ast$. As in Sections \ref{Subsection-Proof of Theorem 1.1} and \ref{Subsection-Proof of Theorem 1.2-n=3}, it suffices to verify both conditions at $t=t_\ast$ and then show that $t_\ast=0$.
		
		We begin with a uniform estimate for $\phi_t$. The boundary convexity estimate and the continuous boundary extension of $\phi_t$ are uniform for $t$ near $t_\ast$. Since $\Omega_{t_\ast}$ is strictly horo-convex, there exists $\varepsilon\in(0,1)$ independent of $t$ such that $$\phi_t(x)\leq 1-\varepsilon\quad\mathrm{on}\ \{x\in\Omega_t:0<e^{-v_t(x)}\leq\varepsilon\}$$ whenever $t>t_\ast$ is sufficiently close to $t_\ast$. On $\{x\in\Omega_t:e^{-v_t(x)}=\varepsilon\}$, we therefore have $e^{-v_t}(1-\phi_t)\geq\varepsilon^2$. By Proposition \ref{Proposition-Maximum principle-n-dim}, $$e^{-v_t}(1-\phi_t)\geq\varepsilon^2\quad\mathrm{in}\ \{x\in\Omega_t:e^{-v_t(x)}\geq\varepsilon\}.$$ Since the normalized eigenfunctions $u_t=e^{-v_t}$ depend smoothly on $t$, they are uniformly bounded above for $t$ close to $t_\ast$. Consequently, the preceding estimates yield a constant $C\in(0,1)$ independent of $t$ such that
		\begin{equation}\label{n-dim-Estimate-phi_t}
			\phi_t\leq C\quad\mathrm{in}\ \Omega_t.
		\end{equation}
		
		We now estimate $\phi_t$ in \eqref{Auxiliary function-Higher dimensional} from below. At a fixed regular point, choose any two coordinate directions and label them $1,2$. Rotate these directions within their plane so that $(v_t)_{12}=0$ and rotate the remaining directions to arrange $(v_t)_\alpha=0$, for $\alpha=3,4,\cdots,n-1$. Thus, the gradient belongs to the three-dimensional space spanned by the first, second and last directions. Within this space, the plane orthogonal to the gradient is tangent to the level hypersurface. The maximum in \eqref{Auxiliary function-Higher dimensional} is therefore at least the value corresponding to this plane. Hence, by the block inverse formula and the equivalence of \eqref{Auxiliary function-3-dim} and \eqref{Auxiliary function-3-dim'}, we have
		\begin{align*}
			\phi_t
			\geq&~\frac{\sum_{l,m\in\{1,2,n\}}(v_t)_l(v_t)_m(v_t)_{lm}}{\det\big((v_t)_{lm}\big)_{l,m\in\{1,2,n\}}}\\
			=&~\frac{(v_t)_n^2}{(v_t)_{11}(v_t)_{22}}\\
			&+\frac{(v_t)_{11}}{\det\big((v_t)_{lm}\big)_{l,m\in\{1,2,n\}}}\bigg((v_t)_1+\frac{(v_t)_{1n}}{(v_t)_{11}}(v_t)_n\bigg)^2\\
			&+\frac{(v_t)_{22}}{\det\big((v_t)_{lm}\big)_{l,m\in\{1,2,n\}}}\bigg((v_t)_2+\frac{(v_t)_{2n}}{(v_t)_{22}}(v_t)_n\bigg)^2\\
			\geq&~\frac{(v_t)_n^2}{(v_t)_{11}(v_t)_{22}},
		\end{align*}
		where $(v_t)_{12}=0$ has been used in the equality. We may therefore write the estimate in the original coordinates since these rotations preserve both $(v_t)_{11}(v_t)_{22}-(v_t)_{12}^2$ and $\sum_{k=3}^n(v_t)_k^2$. As the first two directions were arbitrary, \eqref{n-dim-Estimate-phi_t} implies
		\begin{equation}\label{n-dim-Estimate-gradient}
			\sum_{k\neq i,j}(v_t)_k^2\leq C\big((v_t)_{ii}(v_t)_{jj}-(v_t)_{ij}^2\big)
		\end{equation}
		for all $i\neq j$ in every orthonormal frame.
		
		For the calculations at $t=t_\ast$, write $v=v_{t_\ast}$ and $\lambda=\lambda_{t_\ast}$. By continuity, $\nabla^2v\geq 0$ and \eqref{n-dim-Estimate-gradient} yields
		\begin{equation}\label{n-dim-Estimate-gradient'}
			\sum_{k\neq i,j}v_k^2\leq C(v_{ii}v_{jj}-v_{ij}^2).
		\end{equation}
		Suppose that $\nabla^2v$ degenerates at some $x_0\in\Omega_{t_\ast}$. Choose normal orthonormal coordinates at $x_0$ such that $\nabla^2v(x_0)$ is diagonal. Since equation \eqref{Equation-v} guarantees $\Delta v>0$, its rank lies between $1$ and $n-1$. All calculations below are carried out at $x_0$. The coordinate directions are extended as in Section \ref{Subsection-Minimum principle}.
		
		\textbf{Case 1.} $\nabla^2v$ has rank $n-1$. We may assume that $v_{ii}>0$, $i=1,2,\cdots,n-1$ and $v_{nn}=0$. Taking the indices $n,i$ in \eqref{n-dim-Estimate-gradient'} with $i=1,2,\cdots,n-1$ shows that $v_1=v_2=\cdots=v_{n-1}=0$. Moreover, $v_{nn}$ attains its minimum value $0$ at $x_0$, so $v_{nnk}=0$ and $v_{nnkk}\geq 0$. Thus, \eqref{n-dim-4th derivative term} and equation \eqref{Equation-v} reduce to $$\sum_kv_{nnkk}=2\lambda+2v_n^2=2\sum_{j=1}^{n-1}v_{jj}.$$ For distinct $i,j\in\{1,2,\cdots,n-1\}$, inequality \eqref{n-dim-Estimate-gradient'} shows that $$C(v_{nn}v_{ii}-v_{ni}^2)-\sum_{k\neq n,i}v_k^2$$ attains its minimum value $0$ at $x_0$. Its second derivative in the $x_j$-direction therefore satisfies
		\begin{equation}\label{n-dim-Estimate-4th}
			0\leq C(v_{ii}v_{nnjj}-2v_{nij}^2)-2v_{jj}^2.
		\end{equation}
		Take $i=j+1$ for $j=1,2,\cdots,n-2$ and take $i=1$ when $j=n-1$. Then, by \eqref{n-dim-Estimate-gradient'}, \eqref{n-dim-Estimate-4th} and the Cauchy-Schwarz inequality,
		\begin{align*}
			\sum_{j=1}^{n-1}v_{jj}
			=&~\frac{1}{2}\sum_{k=1}^nv_{nnkk}\geq\frac{1}{2}\sum_{j=1}^{n-1}v_{nnjj}\\
			\geq&~\frac{1}{C}\bigg(\sum_{j=1}^{n-2}\frac{v_{jj}^2}{v_{j+1,j+1}}+\frac{v_{n-1,n-1}^2}{v_{11}}\bigg)\\
			\geq&~\frac{1}{C}\frac{\big(\sum_{j=1}^{n-1}v_{jj}\big)^2}{\sum_{j=1}^{n-1}v_{jj}}=\frac{1}{C}\sum_{j=1}^{n-1}v_{jj}.
		\end{align*}
		This contradicts $C<1$ and $\sum_{j=1}^{n-1}v_{jj}>0$.
		
		\textbf{Case 2.} $\nabla^2v$ has rank at most $n-2$. In this case, take $v_{11}>0$ and $v_{n-1,n-1}=v_{nn}=0$. As in Case 2 of Section \ref{Subsection-Proof of Theorem 1.2-n=3}, the first derivatives of $v_{n-1,n-1}$, $v_{nn}$ and $v_{n-1,n}$ vanish at $x_0$. By \eqref{n-dim-Estimate-gradient'}, we also have $v_k=0$ for $k=1,2,\cdots,n-2$. Consequently, the nonnegative function $$C(v_{n-1,n-1}v_{nn}-v_{n-1,n}^2)-\sum_{k=1}^{n-2}v_k^2$$ vanishes at $x_0$. Its second derivative in the $x_1$-direction yields $0\leq -2v_{11}^2<0$, again a contradiction.
		
		We have now proved that $\nabla^2v_{t_\ast}>0$ throughout $\Omega_{t_\ast}$, so $\phi_{t_\ast}$ is well defined. By continuity, \eqref{n-dim-Estimate-phi_t} implies $$\phi_{t_\ast}\leq C<1.$$ Thus, both conditions remain valid at $t=t_\ast$.
		
		As in Section \ref{Subsection-Proof of Theorem 1.1}, continuity and the boundary convexity estimate ensure that both conditions continue to hold for times slightly below $t_\ast$ if $t_\ast>0$. This contradicts the definition of $t_\ast$. Therefore, $t_\ast=0$, and $\nabla^2v>0$ and $\phi<1$ throughout $\Omega$.
		
		Finally, the argument at the end of Section \ref{Subsection-Proof of Theorem 1.1} shows that $u$ has a unique critical point, its maximum point. Every $\Sigma_t$, $0<t<\max_\Omega u$, is smooth, closed and strictly convex. Let $\kappa_1\leq\kappa_2\leq\cdots\leq\kappa_{n-1}$ denote its principal curvatures. The comparison following \eqref{Auxiliary function-Higher dimensional} yields $$\frac{1}{\kappa_1\kappa_2}\leq\phi<1.$$ We conclude that $\kappa_1\kappa_2>1$. This also implies $\kappa_2>1$ and hence $$K=(\kappa_1\kappa_2)\prod_{\alpha=3}^{n-1}\kappa_\alpha>1.$$
	\end{proof}
	
	We conclude this section by weakening the strict horo-convexity assumption on $\Omega$ in Theorem \ref{Theorem-Higher-dimensional} to horo-convexity.
	
	\begin{proof}[Proof of Corollary \ref{Corollary}]
		Let $\Omega_\varepsilon$, $u_\varepsilon$ and $v_\varepsilon$ be as in the proof of Corollary \ref{Corollary-n=2}. Theorem \ref{Theorem-Higher-dimensional} and the estimates in Sections \ref{Subsection-Proof of Theorem 1.2-n=3} and \ref{Subsection-Proof of Theorem 1.2-n} give $\nabla^2v_\varepsilon>0$ and $\phi_\varepsilon<1$ throughout $\Omega_\varepsilon$. Local smooth convergence therefore implies $\nabla^2v\geq 0$ in $\Omega$ and $\phi\leq 1$ wherever $\nabla^2v>0$. 
		
		We begin with a uniform estimate for $\phi_\varepsilon$. By Proposition \ref{Proposition-Boundary convexity estimate} and Hopf's lemma, $\nabla^2v>0$ and $\nabla v\neq 0$ in a connected neighborhood of $\partial\Omega$. We first show that $\phi<1$ in this neighborhood. A geodesic ball containing $\Omega$ touches $\partial\Omega$ at a point where all boundary principal curvatures exceed $1$, so Proposition \ref{Proposition-Boundary behavior} ensures that $\phi<1$ near this point. A similar calculation retaining the gradient terms in the proofs of Propositions \ref{Proposition-Maximum principle-3-dim} and \ref{Proposition-Maximum principle-n-dim} shows $$\Delta\phi\geq -C(1-\phi+|\nabla\phi|),\quad 0<\phi\leq 1,$$ where $C$ is locally bounded. For $n\geq 4$, we use the smooth function in \eqref{Auxiliary function-Higher dimensional'} as in the proof of Proposition \ref{Proposition-Maximum principle-n-dim}. Since $\phi<1$ near the boundary point chosen above, the strong maximum principle then establishes $\phi<1$ throughout the connected boundary neighborhood.
		
		We now use the inequality $\phi<1$ near $\partial\Omega$ to estimate $\phi_\varepsilon$ on a fixed level set. Choose $t>0$ sufficiently small that $\{0<u\leq t\}$ lies in the preceding neighborhood. Compactness and local smooth convergence then ensure that $$\phi_\varepsilon\leq C<1\quad\mathrm{on}\ \Sigma_t,$$ uniformly for small $\varepsilon>0$.
		
		For $n=3$, Proposition \ref{Proposition-Maximum principle-3-dim} extends this bound to $\{u>t\}$. For $n\geq 4$, we instead apply Proposition \ref{Proposition-Maximum principle-n-dim} to $e^{-v_\varepsilon}(1-\phi_\varepsilon)$. Since $u_\varepsilon\rightarrow t>0$ uniformly on $\Sigma_t$ and $\phi_\varepsilon\leq C<1$ there, the quantity $e^{-v_\varepsilon}(1-\phi_\varepsilon)$ has a uniform positive lower bound on $\Sigma_t$. Proposition \ref{Proposition-Maximum principle-n-dim} gives the same lower bound in $\{u>t\}$. Together with the uniform upper bound for $u_\varepsilon=e^{-v_\varepsilon}$, this yields, after increasing $C$ if necessary, $$\phi_\varepsilon\leq C<1\quad\mathrm{in}\ \{u>t\},$$ with $C$ independent of $\varepsilon$.
		
		It remains to exclude degeneration of $\nabla^2v$ in $\{u>t\}$. By continuity, the corresponding estimates in Sections \ref{Subsection-Proof of Theorem 1.2-n=3} and \ref{Subsection-Proof of Theorem 1.2-n} now hold for $v$ with $C<1$, so the degeneration arguments in those sections apply. Hence, $\nabla^2v>0$ in $\{u>t\}$. Together with Proposition \ref{Proposition-Boundary convexity estimate}, this proves $$\nabla^2v>0\quad\mathrm{in}\ \Omega.$$
	\end{proof}

	\section*{AI use disclosure}
	
	This paper was written by the authors and was not generated by artificial intelligence. The ideas and proof strategies are the authors' own. AI tools were used primarily to explore auxiliary functions, especially in higher dimensions, and to assist with proof details and language editing. All calculations and proof details have been verified by the authors, who take full responsibility for the manuscript.
	
	\bibliographystyle{plain}
	
	\bibliography{mybibliography-hyperbolic}

@article {Acker-PP1981,
	AUTHOR = {Acker, A. and Payne, L. E. and Philippin, G.},
	TITLE = {On the convexity of level lines of the fundamental mode in the
	clamped membrane problem, and the existence of convex
	solutions in a related free boundary problem},
	JOURNAL = {Z. Angew. Math. Phys.},
	FJOURNAL = {Zeitschrift f\"ur Angewandte Mathematik und Physik. ZAMP.
	Journal of Applied Mathematics and Physics. Journal de
	Math\'ematiques et de Physique Appliqu\'ees},
	VOLUME = {32},
	YEAR = {1981},
	NUMBER = {6},
	PAGES = {683--694},
	ISSN = {0044-2275,1420-9039},
	MRCLASS = {35P99 (35R35 53C45 73K15)},
	MRNUMBER = {648766},
	MRREVIEWER = {John\ M.\ Ball},
	DOI = {10.1007/BF00946979},
	URL = {https://doi.org/10.1007/BF00946979},
}

@article {Alvarez-LL1997,
	AUTHOR = {Alvarez, O. and Lasry, J.-M. and Lions, P.-L.},
	TITLE = {Convex viscosity solutions and state constraints},
	JOURNAL = {J. Math. Pures Appl. (9)},
	FJOURNAL = {Journal de Math\'{e}matiques Pures et Appliqu\'{e}es. Neuvi\`eme S\'{e}rie},
	VOLUME = {76},
	YEAR = {1997},
	NUMBER = {3},
	PAGES = {265--288},
	ISSN = {0021-7824},
	MRCLASS = {35J60 (35D10 49L25)},
	MRNUMBER = {1441987},
	MRREVIEWER = {Katsuyuki Ishii},
	DOI = {10.1016/S0021-7824(97)89952-7},
	URL = {https://doi.org/10.1016/S0021-7824(97)89952-7},
}

@article {Andrews-Clutterbuck2011,
	AUTHOR = {Andrews, Ben and Clutterbuck, Julie},
	TITLE = {Proof of the fundamental gap conjecture},
	JOURNAL = {J. Amer. Math. Soc.},
	FJOURNAL = {Journal of the American Mathematical Society},
	VOLUME = {24},
	YEAR = {2011},
	NUMBER = {3},
	PAGES = {899--916},
	ISSN = {0894-0347,1088-6834},
	MRCLASS = {35J10 (35P15)},
	MRNUMBER = {2784332},
	MRREVIEWER = {J.\ B.\ Kennedy},
	DOI = {10.1090/S0894-0347-2011-00699-1},
	URL = {https://doi.org/10.1090/S0894-0347-2011-00699-1},
}

@article {Bourni-CNSWW2022,
	AUTHOR = {Bourni, Theodora and Clutterbuck, Julie and Nguyen, Xuan Hien
		and Stancu, Alina and Wei, Guofang and Wheeler,
		Valentina-Mira},
	TITLE = {The vanishing of the fundamental gap of convex domains in
		{$\Bbb {H}^n$}},
	JOURNAL = {Ann. Henri Poincar\'e},
	FJOURNAL = {Annales Henri Poincar\'e. A Journal of Theoretical and
		Mathematical Physics},
	VOLUME = {23},
	YEAR = {2022},
	NUMBER = {2},
	PAGES = {595--614},
	ISSN = {1424-0637,1424-0661},
	MRCLASS = {58J50},
	MRNUMBER = {4386445},
	MRREVIEWER = {Shoo\ Seto},
	DOI = {10.1007/s00023-021-01096-3},
	URL = {https://doi.org/10.1007/s00023-021-01096-3},
}

@article {Brascamp-Lieb1976,
	AUTHOR = {Brascamp, Herm Jan and Lieb, Elliott H.},
	TITLE = {On extensions of the {B}runn-{M}inkowski and
		{P}r\'{e}kopa-{L}eindler theorems, including inequalities for log
		concave functions, and with an application to the diffusion
		equation},
	JOURNAL = {J. Funct. Anal.},
	VOLUME = {22},
	YEAR = {1976},
	NUMBER = {4},
	PAGES = {366--389},
	MRCLASS = {26A87 (26A51)},
	MRNUMBER = {0450480},
	MRREVIEWER = {L. Leindler},
	DOI = {10.1016/0022-1236(76)90004-5},
	URL = {https://doi.org/10.1016/0022-1236(76)90004-5},
}

@article {Caffarelli-Friedman1985,
	AUTHOR = {Caffarelli, Luis A. and Friedman, Avner},
	TITLE = {Convexity of solutions of semilinear elliptic equations},
	JOURNAL = {Duke Math. J.},
	FJOURNAL = {Duke Mathematical Journal},
	VOLUME = {52},
	YEAR = {1985},
	NUMBER = {2},
	PAGES = {431--456},
	ISSN = {0012-7094,1547-7398},
	MRCLASS = {35B99 (35J99)},
	MRNUMBER = {792181},
	MRREVIEWER = {L.\ E.\ Payne},
	DOI = {10.1215/S0012-7094-85-05221-4},
	URL = {https://doi.org/10.1215/S0012-7094-85-05221-4},
}

@article {Cho2026,
	AUTHOR = {Cho, Gunhee},
	TITLE = {Log-{C}oncavity and {F}undamental {G}aps on {C}omplex {P}rojective {D}omains},
	JOURNAL = {arXiv:2610.35431},
	FJOURNAL = {},
	VOLUME = {},
	YEAR = {2026},
	NOTE = {Preprint},
	NUMBER = {},
	PAGES = {},
	ISSN = {},
	MRCLASS = {},
	MRNUMBER = {},
	DOI = {},
	URL = {https://doi.org/10.48550/arXiv.2610.35431},
}

@article {Dai-ENW2026,
	AUTHOR = {Dai, Xianzhe and Ennis, John M. and Nguyen, Xuan Hien and Wei, Guofang},
	TITLE = {Log-concavity and level-set horoconvexity of the first eigenfunction on horoconvex domains in the hyperbolic plane},
	JOURNAL = {arXiv:2607.27120},
	FJOURNAL = {},
	VOLUME = {},
	YEAR = {2026},
	NOTE = {Preprint},
	NUMBER = {},
	PAGES = {},
	ISSN = {},
	MRCLASS = {},
	MRNUMBER = {},
	DOI = {},
	URL = {https://doi.org/10.48550/arXiv.2607.27120},
}

@article {Dai-SW2021,
	AUTHOR = {Dai, Xianzhe and Seto, Shoo and Wei, Guofang},
	TITLE = {Fundamental gap estimate for convex domains on sphere---the
		case {$n = 2$}},
	JOURNAL = {Comm. Anal. Geom.},
	FJOURNAL = {Communications in Analysis and Geometry},
	VOLUME = {29},
	YEAR = {2021},
	NUMBER = {5},
	PAGES = {1095--1125},
	ISSN = {1019-8385,1944-9992},
	MRCLASS = {58J50},
	MRNUMBER = {4349140},
	DOI = {10.4310/CAG.2021.v29.n5.a3},
	URL = {https://doi.org/10.4310/CAG.2021.v29.n5.a3},
}

@article {He-WZ2020,
	AUTHOR = {He, Chenxu and Wei, Guofang and Zhang, Qi S.},
	TITLE = {Fundamental gap of convex domains in the spheres},
	JOURNAL = {Amer. J. Math.},
	FJOURNAL = {American Journal of Mathematics},
	VOLUME = {142},
	YEAR = {2020},
	NUMBER = {4},
	PAGES = {1161--1191},
	ISSN = {0002-9327,1080-6377},
	MRCLASS = {58J50 (31B10 35J61 35R01)},
	MRNUMBER = {4124117},
	MRREVIEWER = {Alberto\ G.\ Setti},
	DOI = {10.1353/ajm.2020.0033},
	URL = {https://doi.org/10.1353/ajm.2020.0033},
}

@article {Hsu2026,
	AUTHOR = {Hsu, Jia-Lin},
	TITLE = {Strict {L}og-{C}oncavity of {F}irst {D}irichlet {E}igenfunctions under {S}ymmetry in $\mathbb{CP}^n$},
	JOURNAL = {arXiv:2610.00716},
	FJOURNAL = {},
	VOLUME = {},
	YEAR = {2026},
	NOTE = {Preprint},
	NUMBER = {},
	PAGES = {},
	ISSN = {},
	MRCLASS = {},
	MRNUMBER = {},
	DOI = {},
	URL = {https://doi.org/10.48550/arXiv.2610.00716},
}

@article{Hu-LW2022,
	AUTHOR = {Hu, Yingxiang and Li, Haizhong and Wei, Yong},
	TITLE = {Locally constrained curvature flows and geometric inequalities in hyperbolic space},
	FJOURNAL = {Mathematische Annalen},
	JOURNAL = {Math. Ann.},
	ISSN = {0025-5831},
	VOLUME = {382},
	NUMBER = {3-4},
	PAGES = {1425--1474},
	YEAR = {2022},
	DOI = {10.1007/s00208-020-02076-4},
	MRCLASS = {53E10,53E99,53C40,52A39},
}

@article {Jia-MS2023,
	AUTHOR = {Jia, Xiaohan and Ma, Xi-Nan and Shi, Shujun},
	TITLE = {Convexity estimates for {G}reen's function and the first
	eigenfunction of {L}aplace operator},
	JOURNAL = {Potential Anal.},
	FJOURNAL = {Potential Analysis. An International Journal Devoted to the
	Interactions between Potential Theory, Probability Theory,
	Geometry and Functional Analysis},
	VOLUME = {59},
	YEAR = {2023},
	NUMBER = {4},
	PAGES = {1525--1545},
	ISSN = {0926-2601,1572-929X},
	MRCLASS = {35J05 (31B35 35B45 35J08 35P15)},
	MRNUMBER = {4684367},
	DOI = {10.1007/s11118-022-10023-y},
	URL = {https://doi.org/10.1007/s11118-022-10023-y},
}

@book {Kawohl1985,
	AUTHOR = {Kawohl, Bernhard},
	TITLE = {Rearrangements and convexity of level sets in {PDE}},
	SERIES = {Lecture Notes in Mathematics},
	VOLUME = {1150},
	PUBLISHER = {Springer-Verlag, Berlin},
	YEAR = {1985},
	PAGES = {iv+136},
	ISBN = {3-540-15693-3},
	MRCLASS = {35-02 (35B50 35J60 49A50)},
	MRNUMBER = {810619},
	MRREVIEWER = {Michael\ Wiegner},
	DOI = {10.1007/BFb0075060},
	URL = {https://doi.org/10.1007/BFb0075060},
}

@article {Kennington1985,
	AUTHOR = {Kennington, Alan U.},
	TITLE = {Power concavity and boundary value problems},
	JOURNAL = {Indiana Univ. Math. J.},
	FJOURNAL = {Indiana University Mathematics Journal},
	VOLUME = {34},
	YEAR = {1985},
	NUMBER = {3},
	PAGES = {687--704},
	ISSN = {0022-2518},
	MRCLASS = {35J65 (35B99)},
	MRNUMBER = {794582},
	MRREVIEWER = {Maria J. Esteban},
	DOI = {10.1512/iumj.1985.34.34036},
	URL = {https://doi.org/10.1512/iumj.1985.34.34036},
}

@article {Khan-NTW2025,
	AUTHOR = {Khan, Gabriel and Nguyen, Xuan Hien and Tuerkoen, Malik and
		Wei, Guofang},
	TITLE = {Log-concavity and fundamental gaps on surfaces of positive
		curvature},
	JOURNAL = {Comm. Anal. Geom.},
	FJOURNAL = {Communications in Analysis and Geometry},
	VOLUME = {33},
	YEAR = {2025},
	NUMBER = {1},
	PAGES = {239--260},
	ISSN = {1019-8385,1944-9992},
	MRCLASS = {58J50 (53C21)},
	MRNUMBER = {4870313},
	MRREVIEWER = {Shoo\ Seto},
	DOI = {10.4310/cag.250221035055},
	URL = {https://doi.org/10.4310/cag.250221035055},
}

@article {Korevaar1983b,
	AUTHOR = {Korevaar, Nicholas J.},
	TITLE = {Convex solutions to nonlinear elliptic and parabolic boundary
		value problems},
	JOURNAL = {Indiana Univ. Math. J.},
	FJOURNAL = {Indiana University Mathematics Journal},
	VOLUME = {32},
	YEAR = {1983},
	NUMBER = {4},
	PAGES = {603--614},
	ISSN = {0022-2518},
	MRCLASS = {35J65 (35B50 35K60)},
	MRNUMBER = {703287},
	DOI = {10.1512/iumj.1983.32.32042},
	URL = {https://doi.org/10.1512/iumj.1983.32.32042},
}

@incollection {Korevaar1987,
	AUTHOR = {Korevaar, Nicholas J.},
	TITLE = {Convexity properties of solutions to elliptic {PDE}s},
	BOOKTITLE = {Variational methods for free surface interfaces ({M}enlo
		{P}ark, {C}alif., 1985)},
	PAGES = {115--121},
	PUBLISHER = {Springer, New York},
	YEAR = {1987},
	ISBN = {0-387-96396-0},
	MRCLASS = {35J15},
	MRNUMBER = {872895},
}

@article {Korevaar-Lewis1987,
	AUTHOR = {Korevaar, Nicholas J. and Lewis, John L.},
	TITLE = {Convex solutions of certain elliptic equations have constant
		rank {H}essians},
	JOURNAL = {Arch. Rational Mech. Anal.},
	FJOURNAL = {Archive for Rational Mechanics and Analysis},
	VOLUME = {97},
	YEAR = {1987},
	NUMBER = {1},
	PAGES = {19--32},
	ISSN = {0003-9527},
	MRCLASS = {35J60},
	MRNUMBER = {856307},
	MRREVIEWER = {D. E. Edmunds},
	DOI = {10.1007/BF00279844},
	URL = {https://doi.org/10.1007/BF00279844},
}

@article {Lee-Wang1987,
	AUTHOR = {Lee, Yng Ing and Wang, Ai Nung},
	TITLE = {Estimate of {$\lambda_2-\lambda_1$} on spheres},
	JOURNAL = {Chinese J. Math.},
	FJOURNAL = {Chinese Journal of Mathematics},
	VOLUME = {15},
	YEAR = {1987},
	NUMBER = {2},
	PAGES = {95--97},
	ISSN = {0379-7570},
	MRCLASS = {58G25},
	MRNUMBER = {909159},
	MRREVIEWER = {G\'erard\ Besson},
}

@article {Ma-SY2012,
	AUTHOR = {Ma, Xi-Nan and Shi, Shujun and Ye, Yu},
	TITLE = {The convexity estimates for the solutions of two elliptic
	equations},
	JOURNAL = {Comm. Partial Differential Equations},
	FJOURNAL = {Communications in Partial Differential Equations},
	VOLUME = {37},
	YEAR = {2012},
	NUMBER = {12},
	PAGES = {2116--2137},
	ISSN = {0360-5302,1532-4133},
	MRCLASS = {35J05 (35B30 35J25 35P05)},
	MRNUMBER = {3005538},
	MRREVIEWER = {Ian\ Geoffrey\ Wood},
	DOI = {10.1080/03605302.2012.727129},
	URL = {https://doi.org/10.1080/03605302.2012.727129},
}

@article {Ma-Zhang2014,
	AUTHOR = {Ma, Xi-Nan and Zhang, Yongbing},
	TITLE = {The convexity and the {G}aussian curvature estimates for the
	level sets of harmonic functions on convex rings in space
	forms},
	JOURNAL = {J. Geom. Anal.},
	FJOURNAL = {Journal of Geometric Analysis},
	VOLUME = {24},
	YEAR = {2014},
	NUMBER = {1},
	PAGES = {337--374},
	ISSN = {1050-6926,1559-002X},
	MRCLASS = {35J05 (31C05 52A20 53C21)},
	MRNUMBER = {3145927},
	DOI = {10.1007/s12220-012-9339-8},
	URL = {https://doi.org/10.1007/s12220-012-9339-8},
}

@article {Makar-Limanov1971,
	AUTHOR = {Makar-Limanov, L. G.},
	TITLE = {The solution of the {D}irichlet problem for the equation
		{$\Delta u=-1$} in a convex region},
	JOURNAL = {Mat. Zametki},
	FJOURNAL = {Akademiya Nauk SSSR. Matematicheskie Zametki},
	VOLUME = {9},
	YEAR = {1971},
	PAGES = {89--92},
	ISSN = {0025-567X},
	MRCLASS = {31.20},
	MRNUMBER = {279321},
	MRREVIEWER = {M. A. Atakhodzhaev},
}

@article {Singer-WYY1985,
	AUTHOR = {Singer, I. M. and Wong, Bun and Yau, Shing-Tung and Yau,
		Stephen S.-T.},
	TITLE = {An estimate of the gap of the first two eigenvalues in the
		{S}chr\"{o}dinger operator},
	JOURNAL = {Ann. Scuola Norm. Sup. Pisa Cl. Sci. (4)},
	FJOURNAL = {Annali della Scuola Normale Superiore di Pisa. Classe di
		Scienze. Serie IV},
	VOLUME = {12},
	YEAR = {1985},
	NUMBER = {2},
	PAGES = {319--333},
	ISSN = {0391-173X},
	MRCLASS = {35P15 (35J10 58G25)},
	MRNUMBER = {829055},
	MRREVIEWER = {David Gurarie},
	URL = {http://www.numdam.org/item?id=ASNSP_1985_4_12_2_319_0},
}

@article {Seto-WW2019,
	AUTHOR = {Seto, Shoo and Wang, Lili and Wei, Guofang},
	TITLE = {Sharp fundamental gap estimate on convex domains of sphere},
	JOURNAL = {J. Differential Geom.},
	FJOURNAL = {Journal of Differential Geometry},
	VOLUME = {112},
	YEAR = {2019},
	NUMBER = {2},
	PAGES = {347--389},
	ISSN = {0022-040X,1945-743X},
	MRCLASS = {58J50},
	MRNUMBER = {3960269},
	MRREVIEWER = {He-Jun\ Sun},
	DOI = {10.4310/jdg/1559786428},
	URL = {https://doi.org/10.4310/jdg/1559786428},
}

@article {Shih1989,
	AUTHOR = {Shih, Ying},
	TITLE = {A counterexample to the convexity property of the first
		eigenfunction on a convex domain of negative curvature},
	JOURNAL = {Comm. Partial Differential Equations},
	FJOURNAL = {Communications in Partial Differential Equations},
	VOLUME = {14},
	YEAR = {1989},
	NUMBER = {7},
	PAGES = {867--876},
	ISSN = {0360-5302,1532-4133},
	MRCLASS = {58G25 (35P99)},
	MRNUMBER = {1003017},
	MRREVIEWER = {P.\ G\"unther},
	DOI = {10.1080/03605308908820634},
	URL = {https://doi.org/10.1080/03605308908820634},
}

@article {Wei-Xiao2025,
	AUTHOR = {Wei, Guofang and Xiao, Ling},
	TITLE = {Super {L}og-concavity of the {F}irst {E}igenfunctions for {H}oro-convex {D}omains in {H}yperbolic {S}pace},
	JOURNAL = {Anal. PDE},
	FJOURNAL = {Analysis \& PDE},
	VOLUME = {},
	YEAR = {to appear},
	NOTE = {arXiv:2510.13072},
	NUMBER = {},
	PAGES = {},
	ISSN = {},
	MRCLASS = {},
	MRNUMBER = {},
	DOI = {},
	URL = {https://doi.org/10.48550/arXiv.2510.13072},
}

@article {Xia2026,
	AUTHOR = {Xia, Yusen},
	TITLE = {Log-{C}oncavity of {F}irst {D}irichlet {E}igenfunctions on $\mathbb{CP}^2$},
	JOURNAL = {arXiv:2609.22659},
	FJOURNAL = {},
	VOLUME = {},
	YEAR = {2026},
	NOTE = {Preprint},
	NUMBER = {},
	PAGES = {},
	ISSN = {},
	MRCLASS = {},
	MRNUMBER = {},
	DOI = {},
	URL = {https://doi.org/10.48550/arXiv.2609.22659},
}

@article {Zhang-Zhou2026,
	AUTHOR = {Zhang, Wei and Zhou, Qi},
	TITLE = {Power convexity of the torsion function on horo-convex domains in hyperbolic space},
	JOURNAL = {arXiv:2609.02516},
	FJOURNAL = {},
	VOLUME = {},
	YEAR = {2026},
	NOTE = {Preprint},
	NUMBER = {},
	PAGES = {},
	ISSN = {},
	MRCLASS = {},
	MRNUMBER = {},
	DOI = {},
	URL = {https://doi.org/10.48550/arXiv.2609.02516},
}

\end{document}